\documentclass[12pt,reqno]{amsart}
\usepackage{amsmath,amsfonts,amsthm,amsopn,amssymb}
\usepackage{cite,marginnote}
\usepackage{color,enumitem,graphicx}
\usepackage[colorlinks=true,urlcolor=blue,
citecolor=red,linkcolor=blue,linktocpage,pdfpagelabels,
bookmarksnumbered,bookmarksopen]{hyperref}
\usepackage[english]{babel}

\usepackage[left=2.9cm,right=2.9cm,top=2.8cm,bottom=2.8cm]{geometry}

\numberwithin{equation}{section}

\makeindex

\newtheorem{Proposition}{Proposition}[section]

\newtheorem{Remark}{Remark}[section]

\newtheorem{theo}{Theorem}[section]

\newtheorem{lemma}{Lemma}[section]
\newtheorem{iteration lemma}{iteration Lemma}[section]

\newcommand{\bt}{\begin{theorem}}
\newcommand{\et}{\end{theorem}}
\newcommand{\bl}{\begin{lemma}}
\newcommand{\el}{\end{lemma}}
\newcommand{\bd}{\begin{definition}}
\newcommand{\ed}{\end{definition}}
\newcommand{\bc}{\begin{corollary}}
\newcommand{\ec}{\end{corollary}}
\newcommand{\bp}{\begin{proof}}
\newcommand{\ep}{\end{proof}}
\newcommand{\bx}{\begin{example}}
\newcommand{\ex}{\end{example}}
\newcommand{\bi}{\begin{exercise}}
\newcommand{\ei}{\end{exercise}}
\newcommand{\bo}{\begin{proposition}}
\newcommand{\eo}{\end{proposition}}
\newcommand{\br}{\begin{remark}}
\newcommand{\er}{\end{remark}}
\newcommand{\beq}{\begin{equation}}
\newcommand{\eeq}{\end{equation}}
\newcommand{\ba}{\begin{align}}
\newcommand{\ea}{\end{align}}
\newcommand{\bn}{\begin{enumerate}}
\newcommand{\en}{\end{enumerate}}
\newcommand{\bg}{\begin{align*}}
\newcommand{\bcs}{\begin{cases}}
\newcommand{\ecs}{\end{cases}}

\newcommand{\bean}{\begin{eqnarray*}}
\newcommand{\eean}{\end{eqnarray*}}

\def\bd{\mathrm{bd}\,}

\title[Normalized solutions for critical Choquard systems]{Normalized solutions for critical Choquard systems}

\author[H. Zhang]{Hui Zhang}
\author[J. J. Zhang]{Jianjun Zhang}
\author[X.~X.~Zhong]{Xuexiu Zhong}

\address[H.\ Zhang]{\newline\indent Department of  Mathematics,
Jinling Institute of Technology,
\newline\indent
Nanjing 211169, China
\newline\indent and
\newline\indent Department of Mathematics,
Nanjing University,
\newline\indent
Nanjing 210093, China}
\email{\href{mailto:huihz0517@126.com}{huihz0517@126.com}}

\address[J. J. \ Zhang]{\newline\indent College of Mathematics and Statistics, Chongqing Jiaotong University,
\newline\indent
Chongqing 400074, China}
\email{\href{mailto:zhangjianjun09@tsinghua.org.cn}{zhangjianjun09@tsinghua.org.cn}}

\address[X.~X.~Zhong]{\newline\indent South China Research Center for Applied Mathematics and Interdisciplinary Studies
\newline\indent
South China Normal University
\newline\indent
Guangzhou 510631, China}
\email{\href{mailto:zhongxuexiu1989@163.com}{zhongxuexiu1989@163.com}}

\thanks{Hui Zhang was supported by China Postdoctoral Science Foundation (No.2021M691527). Xuexiu Zhong was supported by National Natural Science Foundation of China (No.11801581), Guangdong Basic and Applied Basic Research Foundation (2021A1515010034), Guangzhou Basic and Applied Basic Research Foundation(No.202102020225). Jianjun Zhang was supported by National Natural Science Foundation of China (No.11871123).}

\subjclass[2010]{35J20; 35J50; 35B33.}
\keywords{Choquard system; Normalized solution; Critical exponent; Variational method.}

\begin{document}

\begin{abstract}
In this paper, we consider the critical Choquard system with prescribed mass
\begin{equation*}
\aligned
\left\{ \begin{array}{lll}
-\Delta u+\lambda_1u=(I_\mu\ast |u|^{2^*_\mu})|u|^{2^*_\mu-2}u+\nu p(I_\mu\ast |v|^q)|u|^{p-2}u\ & \text{in}\quad \mathbb{R}^N,\\
-\Delta v+\lambda_2v=(I_\mu\ast |v|^{2^*_\mu})|v|^{2^*_\mu-2}v+\nu q(I_\mu\ast |u|^p)|v|^{q-2}v\ & \text{in}\quad \mathbb{R}^N,\\
\int_{\mathbb{R}^N}u^2=a^2,\quad\int_{\mathbb{R}^N}v^2=b^2,
\end{array}\right.\endaligned
\end{equation*}
where $N\geq3$, $0<\mu<N$, $\nu\in\mathbb{R}$, $I_\mu:\mathbb{R}^N\rightarrow\mathbb{R}$ is a Riesz potential, and $2_{\mu,*}:=\frac{2N-\mu}{N}<p,q<\frac{2N-\mu}{N-2}:=2^*_{\mu},$
with $2_{\mu,*}, 2^*_{\mu}$ called the lower and upper critical exponent in the  sense of Hardy-Littlewood-Sobolev inequality respectively. When $\nu<0$, we prove that no normalized ground state exists. When $\nu>0$, we study the existence, non-existence and asymptotic behavior of normalized solutions by distinguishing three cases: $L^2$-subcritical case: $p+q<4+\frac{4-2\mu}{N}$; $L^2$-critical case: $p+q=4+\frac{4-2\mu}{N}$; $L^2$-supercritical case: $p+q>4+\frac{4-2\mu}{N}$.
In particular,  in $L^2$-subcritical case, and either $N\in\{3,4\}$ or $N\geq5$ with $(\frac N2-1)p+\frac {N}{2}q\leq 2N-\mu$ and $(\frac N2-1)q+\frac {N}{2}p\leq 2N-\mu$, we prove that there exists $\nu_0>0$ such that the system has a positive radial normalized ground state for $0<\nu<\nu_0$. In $L^2$-critical case and $N\in\{3,4\}$, we show there is $\nu'_0>0$ such that  the system has a positive radial normalized ground state for $0<\nu<\nu'_0$. In $L^2$-supercritical case  and $N\in\{3,4\}$, there are two thresholds $\nu_2\geq\nu_1\geq0$ such that
a positive radial
normalized solution exists if $\nu>\nu_2$, and no normalized ground state exists for $\nu<\nu_1$. Furthermore, we give the  concrete ranges of $p$ and $q$ for $\nu_2=\nu_1=0$ and $\nu_2>0$.
\end{abstract}

\maketitle

\section{Introduction and main results}
\renewcommand{\theequation}{1.\arabic{equation}}
In the present paper we are interested in the critical Choquard system with prescribed mass
\begin{equation}\label{0.0}
\aligned
\left\{ \begin{array}{lll}
-\Delta u+\lambda_1u=(I_\mu\ast |u|^{2^*_\mu})|u|^{2^*_\mu-2}u+\nu p(I_\mu\ast |v|^q)|u|^{p-2}u\ & \text{in}\quad \mathbb{R}^N,\\
-\Delta v+\lambda_2v=(I_\mu\ast |v|^{2^*_\mu})|v|^{2^*_\mu-2}v+\nu q(I_\mu\ast |u|^p)|v|^{q-2}v\ & \text{in}\quad \mathbb{R}^N,\\
\int_{\mathbb{R}^N}u^2=a^2,\quad\int_{\mathbb{R}^N}v^2=b^2,
\end{array}\right.\endaligned
\end{equation}
where $a,b>0$, $N\geq3$, $\nu\in\mathbb{R}$, $\ast$ is the convolution product on $\mathbb{R}^N$, $p,q$ satisfy
\vskip 0.1 true cm
 \noindent(H$_0$) $2_{\mu,*}:=\frac{2N-\mu}{N}<p,q<2^{*}_\mu:=\frac{2N-\mu}{N-2}$,
 \vskip 0.1 true cm
 \noindent and $I_\mu:\mathbb{R}^N\rightarrow\mathbb{R}$ is a Riesz potential of order $\mu\in(0,N)$
defined by
\begin{equation}\label{1.1}I_\mu(x)=\frac{A_{N,\mu}}{|x|^\mu}, \ \text{for each}\ x\in\mathbb{R}^N\backslash\{0\},\ \text{where}\quad A_{N,\mu}=\frac{\Gamma(\frac\mu2)}{\Gamma(\frac{N-\mu}2)\pi^{\frac N2}2^{N-\mu}},\end{equation}
with $\Gamma(s)=\int^{+\infty}_0x^{s-1}e^{-x}dx$. In some references,  see \cite{Jianjunzhang,MS} for instance, the Riesz potential is defined by a different form $I_\mu(x)=\frac{A_{N,\mu}}{|x|^{N-\mu}}$.

System (\ref{0.0}) comes from the study of two-component coupled Schr\"{o}dinger system
\begin{equation}\label{1.0}
\aligned
\left\{ \begin{array}{lll}
-i\partial_t\Psi_1=\Delta\Psi_1+(K(x)\ast |\Psi_1|^{r_1})|\Psi_1|^{r_1-2}\Psi_1+\nu p(K(x)\ast |\Psi_2|^{q})|\Psi_1|^{p-2}\Psi_1,  \\
-i\partial_t\Psi_2=\Delta\Psi_2+(K(x)\ast |\Psi_2|^{r_2})|\Psi_2|^{r_2-2}\Psi_2+\nu q (K(x)\ast |\Psi_1|^{p})|\Psi_2|^{q-2}\Psi_2,
\end{array}\right.\endaligned
\end{equation}
which describes several physical phenomenon, such as the Bose-Einstein condensates with multiple states, or the propagation of mutually
incoherent waves packets in nonlinear optics, see \cite{Timm}. In system (\ref{1.0}), $(t,x)\in\mathbb{R}^+\times\mathbb{R}^N$, $\Psi_1$, $\Psi_2$: $\mathbb{R}^+\times\mathbb{R}^N\rightarrow\mathbb{C}$, $K$ is a nonnegative response function which possesses information about the
self-interaction between the particles. The coupling constant $\nu>0$ corresponds to the attraction
and $\nu<0$ to the repulsion between two components in the system. Physically, an important and well-known feature of (\ref{1.0}) is conservation of masses: the $L^2$-norms $|\Psi_1(t,\cdot)|_2$, $|\Psi_2(t,\cdot)|_2$ of solutions are independent of $t\in\mathbb{R}$. Moreover, these norms have a
clear physical meaning, for example, they represent the number of
particles of each component in Bose-Einstein condensates, or the power supply in the
nonlinear optics framework.

An important topic of (\ref{1.0}) is to look for  standing wave solutions $\Psi_1(t,x)=e^{-i\lambda_1t}u(x)$ and $\Psi_2(t,x)=e^{-i\lambda_2t}v(x)$, then $(\Psi_1,\Psi_2)$ solves (\ref{1.0}) if and only if $(u,v)$ is a solution of the system
\begin{equation}\label{1.0.1}
\aligned
\left\{ \begin{array}{lll}
-\Delta u+\lambda_1u={(K(x)\ast |u|^{r_1})|u|^{r_1-2}u}+\nu p(K(x)\ast |v|^q)|u|^{p-2}u\ & \text{in}\quad \mathbb{R}^N,\\
-\Delta v+\lambda_2v=(K(x)\ast |v|^{r_2})|v|^{r_2-2}v+\nu q(K(x)\ast |u|^p)|v|^{q-2}v\ & \text{in}\quad \mathbb{R}^N,
\end{array}\right.\endaligned
\end{equation}
there are two different ways to deal with (\ref{1.0.1}):

\indent(i) the frequencies $\lambda_1$ and $\lambda_2$ are prescribed (fixed frequency problem);\\
\indent(ii) $L^2$-norms $|u|_2,|v|_2$ are prescribed (fixed
mass problem).\\
In this paper, we focus on the case (ii). In this case, the real parameters $\lambda_1$ and $\lambda_2$ appear as Lagrange multipliers, and the solution of (\ref{1.0.1}) with prescribed mass is called a normalized solution.

If the response function is a delta function, i.e. $K(x)=\delta(x)$, the nonlinear response is local and the problem (\ref{1.0.1}) with prescribed mass turns out to be
\begin{equation}\label{1.0.3}
\aligned
\left\{ \begin{array}{lll}
-\Delta u+\lambda_1u=|u|^{2r_1-2}u+\nu p|v|^q|u|^{p-2}u\ & \text{in}\quad \mathbb{R}^N,\\
-\Delta v+\lambda_2v=|v|^{2r_2-2}v+\nu q|u|^p|v|^{q-2}v\ & \text{in}\quad \mathbb{R}^N,\\
\int_{\mathbb{R}^N}u^2=a^2,\ \int_{\mathbb{R}^N}v^2=b^2.
\end{array}\right.\endaligned
\end{equation}
Regarding the existence and qualitative properties of normalized solutions of (\ref{1.0.3}) we refer the reader to \cite{Bartsch-Jeanjean-Soave,
Bartsch-soave-CVPDE,Bartsch-zhong-zou,Gou-jeanjean-Nonl} for Sobolev subcritical case,
and \cite{Bartsch-Li-Zou} for Sobolev critical case. In \cite{Bartsch-Li-Zou}, Bartsch et al.  considered the nonlinear Schr\"{o}dinger system (\ref{1.0.3}) with $2r_1=2r_2=2^*:=2N/(N-2)$
by assuming $N\in\{3,4\}$, $p,q>1$ and $p+q<2^*.$
They proved that a normalized ground
state does not exist for $\nu<0$. When $\nu>0$ and $p+q\leq2+\frac{4}{N}$, the system
has a normalized ground state for $0<\nu<\nu_0$ with $\nu_0$ explicitly
given. In the case $p+q>2+\frac{4}{N}$, they proved the existence of a threshold $\nu_1\geq0$ such that a
normalized ground state exists if and only if $\nu>\nu_1$. Moreover, the asymptotic behavior of normalized ground
states was investigated  as $\nu\rightarrow0^+$ or $\nu\rightarrow+\infty$.

If the response function $K(x)$ is the Reisz potential  function $I_\mu(x)$ given in (\ref{1.1}), then (\ref{1.0.1}) with prescribed mass is changed into
\begin{equation}\label{1.0.6}
\aligned
\left\{ \begin{array}{lll}
-\Delta u+\lambda_1u=(I_\mu\ast |u|^{r_1})|u|^{r_1-2}u+\nu p(I_\mu\ast |v|^q)|u|^{p-2}u\ & \text{in}\quad \mathbb{R}^N,\\
-\Delta v+\lambda_2v=(I_\mu\ast |v|^{r_2})|v|^{r_2-2}v+\nu q(I_\mu\ast |u|^p)|v|^{q-2}v\ & \text{in}\quad \mathbb{R}^N,\\
\int_{\mathbb{R}^N}u^2=a^2,\ \int_{\mathbb{R}^N}v^2=b^2.
\end{array}\right.\endaligned
\end{equation}
 To state the previous results of (\ref{1.0.6}), we first recall the Hardy-Littlewood-Sobolev inequality to explain the meaning of ``critical'' for the nonlocal Choquard system (\ref{1.0.6}).
\begin{Proposition} \label{p1.1} (Hardy-Littlewood-Sobolev inequality) (\cite[Theorem 4.3]{LL})
Let $s,r>1$ and $0<\mu<N$ with $\frac{1}{s}+\frac{\mu}{N}+\frac{1}{r}=2$, $f\in L^s(\mathbb{R}^N)$
and $h\in L^r(\mathbb{R}^N)$. There exists a sharp constant $C(s,N,\mu,r)$, independent of $f$ and $h$, such that
\begin{equation}\label{1.2}\int_{\mathbb{R}^N}\int_{\mathbb{R}^N}\frac{f(x)h(y)}{|x-y|^\mu}dydx\leq C(s,N,\mu,r)|f|_s|h|_r,\end{equation}
where $|\cdot|_q$ is the $L^q(\mathbb{R}^N)$ norm for $q\in[1,+\infty]$. If $s=r=2N/({2N-\mu})$, then
\begin{equation}\label{1.3}C(s,N,\mu,r)=C(N,\mu)=\pi^{\frac{\mu}{2}}
\frac{\Gamma(\frac{N-\mu}{2})}{\Gamma(N-\frac{\mu}{2})}\Bigl(\frac{\Gamma(\frac N2)}{\Gamma(N)}\Bigr)^{-1+\frac{\mu}{N}}.\end{equation}
\end{Proposition}

By Proposition \ref{p1.1}, for any $u\in H^1(\mathbb{R}^N)$, the term
$\int_{\mathbb{R}^N}(I_\mu*|u|^r)|u|^r$
is well defined if
$2_{\mu,*}\leq r\leq 2^*_\mu.$
Thus $2_{\mu,*}$ and $2^*_\mu$ are called the lower and upper  critical exponent in the sense of Hardy-Littlewood-Sobolev inequality respectively. Moreover, Choquard system (\ref{1.0.6}) is usually said to be of critical growth if the nonlinearities is of upper critical growth.

There are some results on normalized solutions of a single Choquard equation with prescribed mass
\begin{equation*}\label{1.0.8}\aligned
\left\{ \begin{array}{lll}-\Delta u+\lambda u=\theta(I_\mu\ast|u|^{r})|u|^{r-2}u+\nu(I_\mu\ast|u|^{p})|u|^{p-2}u \ \ \text{in}\ \mathbb{R}^N,\\
\int_{\mathbb{R}^N}u^2=a^2.
\end{array}\right.\endaligned\end{equation*}
 For instance, see \cite{Li-Ye,Ye,Bartsch-Liu-Liu} for subcritical case, and \cite{Yeyang} for critical case in the sense of Hardy-Littlewood-Sobolev inequality. On the contrary, the results on normalized solutions of the coupled system (\ref{1.0.6}) are few, and as we know, the nonlinearities of (\ref{1.0.6}) in the existing related results \cite{WangYang,wang,gengwang} are all of subcritical growth in the sense of Hardy-Littlewood-Sobolev inequality, i.e. $p,q,r_1,r_2<2^*_\mu$.
 So it is natural to ask whether normalized solutions exist for the critical Choquard system (\ref{0.0})? In the paper, we aim to study the existence, nonexistence and asymptotic behavior of normalized
solutions of (\ref{0.0}).

Denote $H:=H^1(\mathbb{R}^N)\times H^1(\mathbb{R}^N)$ and $|\cdot|_q$ as the $L^q(\mathbb{R}^N)$ norm for $q\in[1,+\infty]$.
Solutions of (\ref{0.0}) are critical points of the functional
$$\aligned J_\nu(u,v)=&\frac12(|\nabla u|^2_2+|\nabla v|^2_2)-\frac{1}{22^*_\mu}\int_{\mathbb{R}^N}\bigl[
(I_\mu*|u|^{2^*_\mu})|u|^{2^*_\mu}+(I_\mu*|v|^{2^*_\mu})
|v|^{2^*_\mu}\bigr]\\&-\nu\int_{\mathbb{R}^N}
(I_\mu*|u|^{p})|v|^{q},\endaligned$$
constrained to the $L^2$-torus
$$\mathcal{S}(a,b):=\{(u,v)\in H:|u|^2_2=a^2,\
|v|^2_2=b^2\}.$$
Before proceeding, we recall the nonlocal type Gagliardo-Nirenberg inequality.

\begin{Proposition}\label{Pro:20230817-p1}(\cite{MS}) Let $N\geq1$, $\mu\in (0,N)$ and $2_{\mu,*}<p<2^*_\mu$, then
\begin{equation}\label{1.11}\int_{\mathbb{R}^N}(I_\mu\ast|u|^p)|u|^p\leq
\frac{p}{|Q_p|^{2p-2}_2}|\nabla u|^{N(p-2)+\mu}_2
|u|^{2N-\mu-(N-2)p}_2,
\end{equation}
where equality holds for $u=Q_p$, where $Q_p$ is  a nontrivial solution of
$$-[N(p-2)+\mu]\Delta Q_p+[2N-\mu-(N-2)p]Q_p=2(I_\mu\ast|Q_p|^p)|Q_p|^{p-2}Q_p.$$
\end{Proposition}

Let $C_{N,p}=p|Q_p|^{2-2p}_2$, then inequality (\ref{1.11}) can be rewritten as
\begin{equation}\label{1.12}\int_{\mathbb{R}^N}(I_\mu\ast|u|^p)|u|^p\leq
C_{N,p}|\nabla u|^{2\gamma_p}_2
|u|^{2p-2\gamma_p}_2,
\end{equation}
where \begin{equation}\label{1.7.0}\gamma_p:=\frac{N(p-2)+\mu}{2}.\end{equation}
To study the system, we need a vector-valued version of the Gagliardo-Nirenberg
inequality. For $2_{\mu,*}<p,q<2^{*}_\mu$,
by (3.3) in \cite{Gvans} and (\ref{1.12}) we have
$$\aligned
\int_{\mathbb{R}^N}(I_\mu*|u|^p)|v|^q&\leq\Bigl(
\int_{\mathbb{R}^N}(I_\mu*|u|^p)|u|^p\Bigr)^{\frac12}\Bigl(
\int_{\mathbb{R}^N}(I_\mu*|v|^q)|v|^q\Bigr)^{\frac12}\\
&\leq C_{N,p}C_{N,q}|\nabla u|^{\gamma_p}_2|\nabla v|^{\gamma_p}_2|u|^{p-\gamma_p}_2|v|^{q-\gamma_q}_2\\
&\leq C_{N,p}C_{N,q}(|\nabla u|^2_2+|\nabla v|^2_2)^{\frac{\gamma_p+\gamma_q}{2}}(|u|^2_2+|v|^2_2)^{\frac{p+q-\gamma_p-\gamma_q}
{2}}.\endaligned$$
Define
$$C^{-1}_{N,p,q}=\inf_{u,v\in H^1(\mathbb{R}^N)\backslash\{0\}}\frac{(|\nabla u|^2_2+|\nabla v|^2_2)^{\frac{\gamma_p+\gamma_q}{2}}
(|u|^2_2+|v|^2_2)^{p+q-\gamma_p-\gamma_q}}
{\int_{\mathbb{R}^N}(I_\mu*|u|^p)|v|^q}.$$
Then $+\infty>C^{-1}_{N,p,q}\geq(C_{N,p}C_{N,q})^{-1}$. From
the above definition we know
\begin{equation}\label{1.14}\int_{\mathbb{R}^N}(I_\mu*|u|^p)|v|^q\leq
C_{N,p,q}(|\nabla u|^2_2+|\nabla v|^2_2)^{\frac{\gamma_p+\gamma_q}{2}}(|u|^2_2+|v|^2_2)
^{\frac{p+q-\gamma_p-\gamma_q}{2}}.
\end{equation}
Clearly, $p+q=4+\frac{4-2\mu}{N}$ is the $L^2$-critical exponent of (\ref{0.0}). In addition, $p+q>4+\frac{4-2\mu}{N}$ and  $p+q<4+\frac{4-2\mu}{N}$ are the $L^2$-supcritical and $L^2$-subcritical exponents of (\ref{0.0}) respectively.

 The solutions of (\ref{0.0}) satisfy the Pohozaev identity
\begin{equation}\label{1.8.0}\aligned P_\nu(u,v)=&|\nabla u|^2_2+|\nabla v|^2_2-\int_{\mathbb{R}^N}\bigl[
(I_\mu*|u|^{2^*_\mu})|u|^{2^*_\mu}+(I_\mu*|v|^{2^*_\mu})
|v|^{2^*_\mu}\bigr]\\&-\nu(\gamma_p+\gamma_q)\int_{\mathbb{R}^N}
(I_\mu*|u|^{p})|v|^{q}=0.\endaligned\end{equation}
Set
\begin{equation*}\label{1.7}
\mathcal{P}_\nu(a,b):=\{(u,v)\in \mathcal{S}(a,b):P_\nu(u,v)=0\},\
m_\nu(a,b):=\inf_{(u,v)\in\mathcal{P}_\nu(a,b) }J_\nu(u,v).\end{equation*}
 A normalized ground state of (\ref{0.0}) is a solution $(u,v)\in \mathcal{P}_\nu(a,b)$ of (\ref{0.0}) that achieves
$m_\nu(a,b)$. Furthermore,
denote $H_{rad}:=\{(u,v)\in H: u,v \text{ are radial}\}$ and
\begin{equation*}\label{1.0.10}m_{r,\nu}(a,b):=\inf_{(u,v)\in \mathcal{P}_\nu(a,b)\cap H_{rad}}J_\nu(u,v).\end{equation*}
 Define $L^2$-invariant scaling $s\star u(x):=e^{\frac{N}{2}s}u(e^sx)$ and $s\star (u,v):=(s\star u,s\star v)$. For $\nu\in\mathbb{R}$ and $(u,v)\in \mathcal{S}(a,b)$, we consider the map $\Psi^\nu_{(u,v)}:\mathbb{R}\rightarrow\mathbb{R}$ defined by
$$\aligned\Psi^\nu_{(u,v)}(t):=J_\nu(t\star(u,v))=&\frac12e^{2t}(|\nabla u|^2_2+|\nabla v|^2_2)-\frac{1}{22^*_\mu}e^{22^*_\mu t}\int_{\mathbb{R}^N}[(I_\mu\ast|u|^{2^*_\mu})
|u|^{2^*_\mu}\\&+(I_\mu\ast|v|^{2^*_\mu})|v|^{2^*_\mu}]-\nu e^{(\gamma_p+\gamma_q)t}\int_{\mathbb{R}^N}(I_\mu\ast|u|^{p})
|v|^{q}.\endaligned$$
An easy computation shows that
$$\mathcal{P}_\nu(a,b)=\bigl\{(u,v)\in \mathcal{S}(a,b):(\Psi^\nu_{(u,v)})'(0)=0\bigr\}.$$
We divide $\mathcal{P}_\nu(a,b)$ into three disjoint sets as follows
$$\mathcal{P}^+_\nu(a,b):=\{(u,v)\in \mathcal{P}_\nu(a,b):(\Psi^\nu_{(u,v)})''(0)>0\},$$
$$\mathcal{P}^0_\nu(a,b):=\{(u,v)\in \mathcal{P}_\nu(a,b):(\Psi^\nu_{(u,v)})''(0)=0\},$$
$$\mathcal{P}^-_\nu(a,b):=\{(u,v)\in \mathcal{P}_\nu(a,b):(\Psi^\nu_{(u,v)})''(0)<0\}.$$
The level $m_\nu(a,b)$ or $m_{r,\nu}(a,b)$ plays a key role in looking for normalized solutions and is involved in a best constant
\begin{equation}\label{1.10} S_{H,L}:=\inf_{u\in D^{1,2}(\mathbb{R}^N)\backslash\{0\}}\frac{|\nabla u|^2_2}{\Bigl(\int_{\mathbb{R}^N}(I_\mu\ast|u|^{2^*_\mu})|u|^{2^*_\mu}\Bigr)^{\frac{1}{2^*_\mu}}}
.\end{equation}

Now we state the main results. Firstly, we consider the  attractive case $\nu>0$.
\begin{theo}\label{th1.1} Let (H$_0$) hold and assume $p+q<4+\frac{4-2\mu}{N}$.\\

\noindent(i) If $N\in\{3,4\}$, then there exists $\nu_0>0$  explicitly given in (\ref{2.13}) below such that, system (\ref{0.0}) admits a normalized ground state $(\bar{u},\bar{v})$ when $0<\nu<\nu_0$. \\

\noindent(ii) If $N\geq5$ with $(\frac N2-1)p+\frac {N}{2}q\leq 2N-\mu$ and $(\frac N2-1)q+\frac {N}{2}p\leq 2N-\mu$, then there exists $\nu_0>0$ such that, system (\ref{0.0}) admits a normalized ground state $(\bar{u},\bar{v})$ when $0<\nu<\nu_0$.\\

\noindent(iii) $(\bar{u},\bar{v})$ obtained in (i) or (ii) is a local
minimizer of $J_\nu|_{\mathcal{S}(a,b)}$ satisfying
$$J_\nu(\bar{u},\bar{v})=m_\nu(a,b)=m_{r,\nu}(a,b)=\inf_{ \mathcal{P}^+_\nu(a,b)}J_\nu<0,\quad \text{when}\ 0<\nu<\nu_0.$$
\end{theo}

\begin{theo}\label{th1.2} Let (H$_0$) hold. Assume $p+q=4+\frac{4-2\mu}{N}$ with $N\in\{3,4\}$, and suppose \begin{equation}\label{1.11.0}\aligned \min\{p,q\}
\left\{ \begin{array}{lll}
>\frac{6-\mu}{2} \text{ with } 2<\mu<3,\ & \text{if}\quad N=3,\\
<2,   \ & \text{if}\quad N=4.
\end{array}\right.\endaligned
\end{equation}

\noindent(i) There exists $\nu'_0$ explicitly given in (\ref{niu0}) below, such that  system (\ref{0.0}) admits a normalized ground state $(\bar{u},\bar{v})$ provided
$0<\nu<\nu'_0$.\\

\noindent(ii) $(\bar{u},\bar{v})$
is of mountain pass type as a critical point of $J_\nu|_{\mathcal{S}(a,b)}$ satisfying
$$J_\nu(\bar{u},\bar{v})=m_\nu(a,b)=m_{r,\nu}(a,b)=\inf_{ \mathcal{P}^-_\nu(a,b)}J_\nu\in
(0,\frac{2^*_\mu-1}{22^*_\mu}S^{\frac{2^*_\mu}{2^*_\mu-1}}_{H,L}),\quad \text{when}\ 0<\nu<\nu'_0.$$
\end{theo}

\begin{theo}\label{th1.3}
 Let (H$_0$) hold and assume $p+q>4+\frac{4-2\mu}{N}$ with $N\in\{3,4\}$,  additionally suppose
\begin{equation}\label{bu1}\min\{p,q\}>\frac{6-\mu}{2},\quad \text{if}\quad N=3.\end{equation}
(i) There exist $\nu_2\geq\nu_1\geq0$ such that system (\ref{0.0}) admits a normalized solution $(\bar{{u}},\bar{{v}})$ if $\nu>\nu_2$, and no normalized ground state exists if $\nu<\nu_1$. Moreover,
\begin{equation*}\aligned
\left\{ \begin{array}{lll}
\nu_2=\nu_1=0,\ & \text{if}\quad N=3,\ \min\{p,q\}<2\ \text{and}\ 2<\mu<3,\\
\nu_2=\nu_1=0,\ & \text{if}\quad N=4\ \text{and}\ \min\{p,q\}<2,\\
\nu_2>0,   \ & \text{if}\ \min\{p,q\}\geq2+\frac{4-\mu}{N}.
\end{array}\right.\endaligned
\end{equation*}
\noindent(ii) $(\bar{{u}},\bar{{v}})$ is of mountain pass type as a critical point of $J_\nu|_{\mathcal{S}(a,b)}$  satisfying $$J_\nu(\bar{{u}},\bar{{v}})=m_{r,\nu}(a,b)\in(0,
\frac{2^*_\mu-1}{22^*_\mu}S^{\frac{2^*_\mu}{2^*_\mu-1}}_{H,L}),\quad \text{when}\ \nu>\nu_2.$$
\end{theo}
\begin{theo}\label{th1.4}Both components of normalized solutions obtained in Theorems \ref{th1.1}-\ref{th1.3} are positive and radially symmetric. The corresponding Lagrange multipliers $\lambda_1,\lambda_2$ are positive.
\end{theo}
\begin{Remark} (i) As we shall see in Section 3, when $N=3$, it is necessary to assume that $\min\{p,q\}>\frac{6-\mu}{2}$ to ensure the least energy $m_\nu(a,b)$ or $m_{r,\nu}(a,b)$ lies in a suitable interval. Moreover, we would like to point out that, in Theorem \ref{th1.3}, $\min\{p,q\}\geq2+\frac{4-\mu}{N}$ with $N=3$ implies $\min\{p,q\}>\frac{6-\mu}{2}$ and so in this case we may only assume $\min\{p,q\}\geq2+\frac{4-\mu}{N}$.\\
\noindent(ii) $m_\nu(a,b)$ and $m_{r,\nu}(a,b)$ in Theorem \ref{th1.3} also satisfy
\begin{equation*}\aligned
\left\{ \begin{array}{lll}
 m_\nu(a,b)=m_{r,\nu}(a,b)=\frac{2^*_\mu-1}{22^*_\mu}S^{\frac{2^*_\mu}{2^*_\mu-1}}_{H,L},  \ & \text{when}\quad 0<\nu\leq\nu_1,\\
 m_\nu(a,b)\leq m_{r,\nu}(a,b)=\frac{2^*_\mu-1}{22^*_\mu}S^{\frac{2^*_\mu}{2^*_\mu-1}}_{H,L},   \ & \text{when}\quad \nu_1\leq\nu\leq\nu_2,\\
 m_\nu(a,b)\leq m_{r,\nu}(a,b)\in(0,\frac{2^*_\mu-1}{22^*_\mu}S^{\frac{2^*_\mu}{2^*_\mu-1}}_{H,L}),   \ & \text{when}\quad \nu>\nu_2,
\end{array}\right.\endaligned
\end{equation*}
and only the last holds if $\nu_1=\nu_2=0$. In addition, $m_\nu(a,b)\leq m_{r,\nu}(a,b)\rightarrow0^+$ when $\nu\rightarrow+\infty.$
\end{Remark}

Next we investigate the asymptotic behavior of normalized solutions of (\ref{0.0}) as $\nu\rightarrow0^+$ or $\nu\rightarrow+\infty$. One type of asymptotic behavior is related to the minimizer of $S_{H,L}$. For this goal, we recall that the best Sobolev embedding constant
$S=\inf_{u\in D^{1,2}(\mathbb{R}^N)\backslash\{0\}}\frac{|\nabla u|^2_2}{|u|^{2}_{2^*}},$
is achieved by the family
\begin{equation}\label{1.16}U_{\epsilon,\xi}(x)=[N(N-2)]^{\frac{N-2}{4}}\epsilon^{\frac{N-2}{2}}
(\epsilon^2+|x-\xi|^2)^{-\frac{N-2}{2}},\ \ x,\xi\in\mathbb{R}^N, \ \epsilon>0,\end{equation}
and the minimizer of $S_{H,L}$ is given by the following lemma.
\begin{lemma} (\cite[Lemma 1.2]{DuYang})\label{l1.1} The constant $S_{H,L}$ defined in (\ref{1.10}) satisfies $S_{H,L}=S\bigl(C(N,\mu)\bigr)^{-\frac{1}{2^*_\mu}},$ where $S$ is the best Sobolev constant and $C(N,\mu)$ is given in (\ref{1.3}). Moreover,
\begin{equation*}\tilde{U}_{\epsilon,\xi}(x)=S^{\frac{(N-\mu)(2-N)}{4(N-\mu+2)}}
\bigl(C(N,\mu)\bigr)^{\frac{2-\mu}{2(N-\mu+2)}}U_{\epsilon,\xi}(x),\end{equation*}
is the unique family of positive minimizers for $S_{H,L}$ that satisfies
\begin{equation*} -\Delta u=(I_\mu*|u|^{2^*_\mu})|u|^{2^*_\mu-2}u,\quad\text{in}\ \ \mathbb{R}^N.\end{equation*}
\end{lemma}

Another type of	asymptotic behavior is related to the solutions of the system
\begin{equation}\label{limit}
\aligned
\left\{ \begin{array}{lll}
-\Delta u+\lambda_1u= p(I_\mu\ast |v|^q)|u|^{p-2}u\ & \text{in}\quad \mathbb{R}^N,\\
-\Delta v+\lambda_2v= q(I_\mu\ast |u|^p)|v|^{q-2}v\ & \text{in}\quad \mathbb{R}^N,\\
\int_{\mathbb{R}^N}u^2=a^2,\quad\int_{\mathbb{R}^N}v^2=b^2.
\end{array}\right.\endaligned
\end{equation}
The functional of (\ref{limit}) is
\begin{equation}\label{limitfuctional}\tilde{J}(u,v)=\frac12(|\nabla u|^2_2+|\nabla v|^2_2)-\int_{\mathbb{R}^N}
(I_\mu*|u|^{p})|v|^{q},\end{equation}
which is well defined in $H$,
and constrained to the $L^2$-torus
$\mathcal{S}(a,b)$.
The solutions of (\ref{limit}) satisfy the Pohozaev identity
\begin{equation*}\label{limitpoho}\tilde{P}(u,v)=|\nabla u|^2_2+|\nabla v|^2_2-(\gamma_p+\gamma_q)\int_{\mathbb{R}^N}
(I_\mu*|u|^{p})|v|^{q}=0.\end{equation*} A normalized
ground state of (\ref{limit}) is a solution $(u,v)\in \mathcal{\tilde{P}}(a,b)$ of (\ref{limit}) satisfies  $\tilde{J}(u,v)=\tilde{m}(a,b)$ where
\begin{equation}\label{1.15}
\tilde{m}(a,b)=\inf_{(u,v)\in\mathcal{\tilde{P}}(a,b) }\tilde{J}(u,v),\ \text{with} \ \mathcal{\tilde{P}}(a,b):=\{(u,v)\in \mathcal{S}(a,b):\tilde{P}(u,v)=0\}. \end{equation}
\begin{theo}\label{th1.5}  Let (H$_0$) hold and assume $p+q\neq 4+\frac{4-2\mu}{N}$. In addition, assume either $N\in\{3,4\}$ or $N\geq5$ with $(\frac N2-1)p+\frac {N}{2}q\leq 2N-\mu$ and $(\frac N2-1)q+\frac {N}{2}p\leq 2N-\mu$, then system (\ref{limit}) admits a  normalized ground state $(u,v)$ with $\lambda_1,\lambda_2>0$. Moreover, $u,v$ are positive and radially symmetric. \end{theo}

Now we state the results of asymptotic behavior of normalized solutions of (\ref{0.0}).

\begin{theo}\label{th1.6} Let $\mathcal{L}(a,b)$ be the set of normalized ground states of (\ref{limit}) with $\lambda_1,\lambda_2>0$. Then the following results hold.\\
\noindent(i) If $p+q<4+\frac{4-2\mu}{N}$ and $(u_\nu,v_\nu)$ is a family of positive radial normalized ground states of (\ref{0.0}) from Theorem \ref{th1.1}  with $\nu\rightarrow0^+$,
then there exist $t_\nu\sim(\gamma_p+\gamma_q-2)^{-1}\ln\nu$ and $(\tilde{u},\tilde{v})\in\mathcal{L}(a,b)\cap H_{rad}$ such that
$$t_\nu\star (u_\nu,v_\nu)\rightarrow(\tilde{u},\tilde{v})\ \text{in } H,\ \quad\text{as}\ \nu\rightarrow0^+.$$
 \noindent(ii) If $p+q>4+\frac{4-2\mu}{N}$ and $(u_\nu,v_\nu)$ is a family of  positive radial normalized solutions of (\ref{0.0}) from Theorem \ref{th1.3} with $\nu\rightarrow+\infty$,
then there exist $t_\nu\sim(\gamma_p+\gamma_q-2)^{-1}\ln\nu$ and $(\tilde{u},\tilde{v})\in\mathcal{L}(a,b)\cap H_{rad}$ such that
$$t_\nu\star (u_\nu,v_\nu)\rightarrow(\tilde{u},\tilde{v})\ \text{in } H, \quad\text{as}\ \nu\rightarrow+\infty.$$
\end{theo}
\begin{theo}\label{th1.7}
\noindent(i) If $p+q=4+\frac{4-2\mu}{N}$ and $(u_\nu,v_\nu)$ is a family of positive radial normalized ground states of (\ref{0.0}) from Theorem \ref{th1.2}  with $\nu\rightarrow0^+$, then there exists $r_\nu>0$ such that
$$\bigl(u_\nu,r^{\frac{N-2}{2}}_\nu v_\nu(r_\nu x)\bigr)\rightarrow(0,\widetilde{U}_{\epsilon,0})\quad \text{or}\quad \bigl(r^{\frac{N-2}{2}}_\nu u_\nu(r_\nu x), v_\nu\bigr)\rightarrow(\widetilde{U}_{\epsilon,0},0),$$
in $D^{1,2}(\mathbb{R}^N)\times D^{1,2}(\mathbb{R}^N)$ as $\nu\rightarrow0^+$ for some $\epsilon>0$.\\
\noindent(ii) If $p+q>4+\frac{4-2\mu}{N}$ and $(u_\nu,v_\nu)$ is a family of positive radial normalized solutions of (\ref{0.0}) from Theorem \ref{th1.3} with $\nu\rightarrow0^+$, then there exists $r_\nu>0$ such that
$$\bigl(u_\nu,r^{\frac{N-2}{2}}_\nu v_\nu(r_\nu x)\bigr)\rightarrow(0,\widetilde{U}_{\epsilon,0})\quad \text{or}\quad \bigl(r^{\frac{N-2}{2}}_\nu u_\nu(r_\nu x), v_\nu\bigr)\rightarrow(\widetilde{U}_{\epsilon,0},0),$$
in $D^{1,2}(\mathbb{R}^N)\times D^{1,2}(\mathbb{R}^N)$ as $\nu\rightarrow0^+$  for some $\epsilon>0$.
\end{theo}

Finally, we consider the repulsive case $\nu<0$.
\begin{theo}\label{th1.8} Let (H$_0$)  hold and assume $\nu\leq0$. Then $m_\nu(a,b)=
\frac{2^*_\mu-1}{22^*_\mu}S^{\frac{2^*_\mu}{2^*_\mu-1}}_{H,L}$ and system (\ref{0.0}) has no normalized ground state.\end{theo}
\begin{Remark}
The proofs in this paper can also be adapted to another form of critical Choquard system
\begin{equation*}
\aligned
\left\{ \begin{array}{lll}
-\Delta u+\lambda_1u=\theta(|x|^{-4}\ast |u|^{2})u+\nu(|x|^{-\alpha}\ast v^2)u+\mu(|x|^{-\beta}\ast v^2)u\ & \text{in}\quad \mathbb{R}^N,\\
-\Delta v+\lambda_2v=\theta(|x|^{-4}\ast |v|^{2})v+\nu(|x|^{-\alpha}\ast u^2)v+\mu(|x|^{-\beta}\ast u^2)v\ & \text{in}\quad \mathbb{R}^N,\\
\int_{\mathbb{R}^N}u^2=a^2,\quad\int_{\mathbb{R}^N}v^2=b^2,
\end{array}\right.\endaligned
\end{equation*} where $0<\alpha,\beta<4$, $\theta\geq0$, $\nu,\mu\in\mathbb{R}$, and similar existence, nonexistence and asymptotic results as Theorems \ref{th1.1}-\ref{th1.8} can be obtained. Here we would like to point out, in the sense of Hardy-Littlewood-Sobolev inequality, $\int_{\mathbb{R}^N}(|x|^{-\alpha}\ast|u|^{2})u^2$ is well defined when $0\leq\alpha\leq4$, $\alpha=0$ and $\alpha=4$ are the lower and upper critical
exponents respectively.
\end{Remark}

Before going to prove our main results, some remarks on these results are in order: (i) the restore of the compactness of Palais-Smale sequence is usually strongly relied on a Liouville type result, which is only valid when $N\leq4$. In this paper, we are allowed to consider the case $N\geq5$ in the $L^2$-subcritical case. On the other hand, when using a Liouville type result for $N\in\{3,4\}$,  the appearance of nonlocal terms leads that  the
integrability of solutions should be improved.  We shall use the fact that the least energy  in the $L^2$-subcritical case is less than zero and some tricks in \cite{Gou-jeanjean-Nonl} to deal with the case $N\geq5$. When $N\in\{3,4\}$, by virtue of Kato inequality and the arguments in \cite{MS} for a single Choquard equation,  we develop the results about improved integrability of solutions for a nonlocal system. (ii) For the
$L^2$-supercritical case, due to the non-homogeneous and non-locality of nonlinear terms, it is
difficult to use the Riesz's rearrangement inequality (see \cite{LL}) to prove that, the least energy restricted on the radial subspace is equal to that on the whole space. Hence, it is necessary to distinguish them and introduce another threshold for the least energy restricted on the radial subspace, we then prove that it is achieved by a normalized solution. (iii) The Hardy-Littlewood-Sobolev critical term and nonlocal coupling term lead that, the estimates of the functional level and the threshold turn out to be more complex. We shall estimate carefully and give the concrete ranges of $p$ and $q$ for the thresholds. (iv) It is not easy to show the nonexistence of solutions in the $L^2$-supercritical case. We shall borrow the idea in \cite{Chen-Zou} and use the minimax characterization of the least energy to show the nonexistence of normalized solutions. (v) Some new properties of two limit problems are needed to analyze the asymptotic behavior of normalized solutions. We shall study two limit problems and in particular give a nonlocal critical compactness lemma, which can be viewed as the counterpart of \cite[Theorem 1.41]{WM} about a local critical compactness lemma.

{\bf Notations}: $C,C_0,\tilde{C},C',C_1, C_2...$  may denote different positive constants. For $1\leq p\leq+\infty$, the norm in $L^p(\Omega)$ is denoted by $|\cdot|_{p,\Omega}$ when $\Omega$ is  a proper subset of $\mathbb{R}^N$, by $|\cdot|_{p}$ when $\Omega=\mathbb{R}^N$. For any $r>0$ and $z\in\mathbb{R}^N$, $B_r(z)$ denotes the ball of radius $r$ centered at $z$ and $B_r=B_r(0)$. $\int_{\mathbb{R}^N}f(x)dx$ is denoted by $\int_{\mathbb{R}^N}f(x)$. $\mathcal{S}_r(a,b)=\mathcal{S}(a,b)\cap H_{rad}$. For simplicity, we write $A_\nu\thicksim B_\nu$, $A_\nu \lesssim B_\nu$ and $A_\nu \gtrsim B_\nu$ if there exist $C_1, C_2>0$ such that, respectively, $C_1B_\nu\leq A_\nu\leq C_2B_\nu$, $A_\nu\leq C_2B_\nu$ and $A_\nu\geq C_1B_\nu$ as $\nu\rightarrow0^+$ or $\nu\rightarrow+\infty$.

The paper is organized as
follows. In Sections 2 and 3, for the case $\nu>0$, we introduce the regularity and compactness results, and then give the existence and nonexistence results of normalized solutions. In Section 4, we study two limit problems and investigate the asymptotic behavior of normalized solutions. In Section 5 we give the nonexistence result in the case $\nu<0$.

\section{The regularity and compactness results in the case $\nu>0$}
\renewcommand{\theequation}{2.\arabic{equation}}
In this section we fix $\nu>0$. To give the compactness lemma, we need a regularity result. Different from \cite{Bartsch-Li-Zou}, the nonlocal term leads that the integrability of solutions of system (\ref{0.0}) should be improved.  Firstly let us recall an important inequality for
nonlocal nonlinearities by Moroz and Van Schaftingen \cite{MS}.
\begin{lemma}\label{l2.0}(\cite[Lemma 3.2]{MS})
Let $N\geq2$, $\mu\in (0,N)$ and $\theta\in(0,2)$. If $H,K\in L^{\frac{2N}{N-\mu+2}}(\mathbb{R}^N)+L^{\frac{2N}{N-\mu}}(\mathbb{R}^N)$ and $1-\frac{\mu}{N}<\theta<1+\frac{\mu}{N}$, then for every $\epsilon>0$, there exists $C_{\epsilon,\theta}\in\mathbb{R}$ such that for every $u\in H^1(\mathbb{R}^N)$,
$$\int_{\mathbb{R}^N}\bigl(I_\mu\ast (H|u|^\theta)\bigr)K|u|^{2-\theta}\leq \epsilon^2\int_{\mathbb{R}^N}|\nabla u|^2+C_{\epsilon,\theta}\int_{\mathbb{R}^N}|u|^2.$$
\end{lemma}

Now we give the result about improved integrability of solutions of a nonlocal system.
\begin{lemma}\label{l2.0.0} Let $N\geq2$ and $\mu_1,\mu_2\in\mathbb{R}$. If $H_{ij}, K_{ij}\in L^{\frac{2N}{N-\mu+2}}(\mathbb{R}^N)+L^{\frac{2N}{N-\mu}}(\mathbb{R}^N)$, $i,j=1,2$,
and $(u,v)\in H$ solves
\begin{equation}\label{2.0.0}
\aligned
\left\{ \begin{array}{lll}
-\Delta u+\mu_1u=(I_\mu\ast H_{11}u)K_{11}+(I_\mu\ast H_{12}v)K_{12}\ & \text{in}\quad \mathbb{R}^N,\\
-\Delta v+\mu_2v=(I_\mu\ast H_{21}v)K_{21}+(I_\mu\ast H_{22}u)K_{22}\ & \text{in}\quad \mathbb{R}^N,
\end{array}\right.\endaligned
\end{equation}
then $u,v\in L^p(\mathbb{R}^N)$ for every $p\in[2,\frac{2N^2}{(N-2)(N-\mu)})$.
\end{lemma}
{\bf Proof}: Using Lemma \ref{l2.0} with $\theta=1$, there exists $\lambda>0$ such that
$$\int_{\mathbb{R}^N}(I_\mu\ast|H_{i1}\phi|)|K_{i1}\phi|\leq\frac14|\nabla\phi|^2_2+
\frac{\lambda}{4}|\phi|^2_2,\quad\forall \phi\in H^1(\mathbb{R}^N),\ i=1,2.$$
By (3.3) in \cite{Gvans} and Lemma \ref{l2.0}, for any $(\phi,\psi)\in H$ we have
$$\aligned
\int_{\mathbb{R}^N}(I_\mu\ast|H_{i2}\phi|)|K_{i2}\psi|&\leq
\Bigl(\int_{\mathbb{R}^N}(I_\mu\ast|H_{i2}\phi|)|H_{i2}\phi|\Bigr)^{\frac12}
\Bigl(\int_{\mathbb{R}^N}(I_\mu\ast|K_{i2}\psi|)|K_{i2}\psi|\Bigr)^{\frac12}\\
&\leq \bigl(\frac14|\nabla\phi|^2_2+
\frac{\lambda}{4}|\phi|^2_2\bigr)^{\frac12}\bigl(\frac14|\nabla\psi|^2_2+
\frac{\lambda}{4}|\psi|^2_2\bigr)^{\frac12}\\&\leq\frac{1}{4}
(|\nabla\phi|^2_2+|\nabla\psi|^2_2)+\frac{\lambda}{4}(|\phi|^2_2+|\psi|^2_2),\ i=1,2.
\endaligned$$
For $i,j\in\{1,2\}$, choose sequences $\{H^{(k)}_{ij}\}_{k\in \mathbb{N}}$, $\{K^{(k)}_{ij}\}_{k\in \mathbb{N}}$  in $L^{\frac{2N}{N-\mu}}(\mathbb{R}^N)$, such that $|H^{(k)}_{ij}|\leq |H_{ij}|$, $|K^{(k)}_{ij}|\leq |K_{ij}|$ and $H^{(k)}_{ij}\rightarrow H_{ij}$ and $K^{(k)}_{ij}\rightarrow K_{ij}$ a.e. in $\mathbb{R}^N$. For each $k\in\mathbb{N}$, the form $a_k:H\times H\rightarrow\mathbb{R}$ defined for $(\phi_1,\phi_2), (\psi_1,\psi_2)\in H$ by
$$\aligned
a_k((\phi_1,\phi_2), (\psi_1,\psi_2))=&
\int_{\mathbb{R}^N}(\nabla\phi_1\nabla\psi_1+\lambda\phi_1\psi_1)
-\Sigma^{2}_{i=1}\int_{\mathbb{R}^N}(I_\mu\ast|H^{(k)}_{1i}\phi_i|)|K^{(k)}_{1i}\psi_1|\\
&+\int_{\mathbb{R}^N}(\nabla\phi_2\nabla\psi_2+\lambda\phi_2\psi_2)
-\Sigma^{2}_{i=1}\int_{\mathbb{R}^N}(I_\mu\ast|H^{(k)}_{2i}\phi_{3-i}|)|K^{(k)}_{2i}
\psi_2|,
\endaligned$$
is bilinear and coercive. Applying the Lax-Milgram theorem, there exists
a unique solution $(u_k,v_k)\in H$ of
\begin{equation}\label{2.0.9}
\aligned
\left\{ \begin{array}{lll}
-\Delta u_k+\lambda u_k=\bigl(I_\mu\ast (H^{(k)}_{11}u_k)\bigr)K^{(k)}_{11}+\bigl(I_\mu\ast (H^{(k)}_{12}v_k)\bigr)K^{(k)}_{12}+(\lambda-\mu_1)u\ & \text{in}\quad \mathbb{R}^N,\\
-\Delta v_k+\lambda v_k=\bigl(I_\mu\ast (H^{(k)}_{21}v_k)\bigr)K^{(k)}_{21}+\bigl(I_\mu\ast (H^{(k)}_{22}u_k)\bigr)K^{(k)}_{22}+(\lambda-\mu_2)v\ & \text{in}\quad \mathbb{R}^N,
\end{array}\right.\endaligned
\end{equation}
where $(u,v)$ is the given solution of (\ref{2.0.0}). It can be proved that $(u_k,v_k)\rightharpoonup (u,v)$ in $H$ as $k\rightarrow+\infty$.

Testing the first equality of (\ref{2.0.9}) with $sgn u_k$, where $sgn$ is the sign function,
we obtain
$$sgn(u_k)(-\Delta u_k+\lambda u_k)=sgn(u_k)\Bigl[\bigl(I_\mu\ast (H^{(k)}_{11}u_k)\bigr)K^{(k)}_{11}+\bigl(I_\mu\ast (H^{(k)}_{12}v_k)\bigr)K^{(k)}_{12}+(\lambda-\mu_1)u\Bigr].$$
By \cite[Lemma A]{Kato}, we get $- sgn(u_k)\Delta u_k\geq-\Delta|u_k|$ in the sense of distribution. Thus \begin{equation}\label{2.0.1}-\Delta |u_k|+\lambda |u_k|\leq(I_\mu\ast |H^{(k)}_{11}||u_k|)|K^{(k)}_{11}|+(I_\mu\ast |H^{(k)}_{12}||v_k|)|K^{(k)}_{12}|+|\lambda-\mu_1||u|.\end{equation}
Similarly
\begin{equation}\label{2.0.2}-\Delta |v_k|+\lambda |v_k|\leq(I_\mu\ast |H^{(k)}_{21}||v_k|)|K^{(k)}_{21}|+(I_\mu\ast |H^{(k)}_{22}||u_k|)K^{(k)}_{22}+|\lambda-\mu_2||v|.\end{equation}
Adding (\ref{2.0.1}) and (\ref{2.0.2}) and letting $$H^{(k)}=\Sigma^{2}_{i,j=1}|H^{(k)}_{ij}|,\ K^{(k)}=\Sigma^{2}_{i,j=1}|K^{(k)}_{ij}|, \ w_k=|u_k|+|v_k|,$$
there holds
\begin{equation}\label{2.0.3}-\Delta w_k+\lambda w_k\leq \bigl(I_\mu\ast(H^{(k)}w_k)\bigr)K^{(k)}+\tilde{C}(|u|+|v|).\end{equation}
for some $\tilde{C}>0$. For $l>0$, define the truncation $w_{k,l}:\mathbb{R}^N\rightarrow\mathbb{R}$ by $w_{k,l}=\min\{w_k,l\}$.
Since $w^{p-1}_{k,l}\in H^1(\mathbb{R}^N)$, we can take it as a test function of (\ref{2.0.3})  and obtain
$$\aligned&
\frac{4(p-1)}{p^2}|\nabla(w_{k,l})^{\frac p2}|^2_2+\lambda|w_{k,l}|^p_p \leq\int_{\mathbb{R}^N}\bigl[
\bigl(I_\mu\ast(H^{(k)}w_k)\bigr)K^{(k)}+
\tilde{C}(|u|+|v|)\bigr](w_{k,l})^{p-1}.
\endaligned$$
Denote $H_0=\Sigma^{2}_{i,j=1}|H_{ij}|$ and $K_0=\Sigma^{2}_{i,j=1}|K_{ij}|$. Then $H^{(k)}\leq H_0$ and $K^{(k)}\leq K_0$. If $p<\frac{2N}{N-\mu}$, by Lemma \ref{l2.0} with $\theta=\frac{2}{p}$, there exists $C>0$ such that
$$\aligned\int_{\mathbb{R}^N}
\bigl(I_\mu\ast(H_0w_{k,l})\bigr)K_0(w_{k,l})^{p-1}\leq
\frac{2(p-1)}{p^2}|\nabla(w_{k,l})^{\frac p2}|^2_2+C|w_{k,l}|^p_p.
\endaligned$$
Therefore
\begin{equation}\label{2.0.4}\frac{2(p-1)}{p^2}
|\nabla(w_{k,l})^{\frac p2}|^2_2\leq C'(|w_k|^p_p+|u|^p_p+|v|^p_p)+\int_{A_{k,l}}
\bigl(I_\mu\ast (K^{(k)}w^{p-1}_k)\bigr)H^{(k)}w_{k},\end{equation}
where $A_{k,l}=\{x\in\mathbb{R}^N:w_k(x)\geq l\}.$
Below taking similar arguments in \cite[Proposition 3.1]{MS} for a scalar Choquard equation, using $p<\frac{2N}{N-\mu}$ we get $\limsup_{k\rightarrow\infty}|w_k|^{p}_{\frac{Np}{N-2}}\leq\limsup_{k\rightarrow\infty}|w_k|^{p}_p$.
Then by iterating over $p$ a finite number of times we cover the range $p\in[2,\frac{2N^2}{(N-2)(N-\mu)})$.\ \ \ \ $\Box$

\begin{lemma}\label{l2.0.1} Let $\mu_1,\mu_2\in\mathbb{R}$. If $(u,v)\in H$ is a  weak solution of the system
\begin{equation}\label{2.0.5}
\aligned
\left\{ \begin{array}{lll}
-\Delta u+\mu_1 u=(I_\mu\ast |u|^{2^*_\mu})|u|^{2^*_\mu-2}u+p(I_\mu\ast |v|^q)|u|^{p-2}u\ & \text{in}\quad \mathbb{R}^N,\\
-\Delta v+\mu_2 v=(I_\mu\ast |v|^{2^*_\mu})|v|^{2^*_\mu-2}v+q(I_\mu\ast |u|^p)|v|^{q-2}v\ & \text{in}\quad \mathbb{R}^N,
\end{array}\right.\endaligned
\end{equation}
then $(u,v)$ is a smooth solution.
\end{lemma}
{\bf Proof}: Rewrite system (\ref{2.0.5}) in the form of (\ref{2.0.0})
with $H_{11}=K_{11}=|u|^{2^*_\mu-2}u$, $H_{21}=K_{21}=|v|^{2^*_\mu-2}v$, $H_{12}=\frac1qK_{22}=|v|^{q-2}v$, $H_{22}=\frac1pK_{12}=|u|^{p-2}u.$ Moreover, for $i,j=1,2$, $H_{ij}, K_{ij}\in L^{\frac{2N}{N-\mu}}(\mathbb{R}^N)+L^{\frac{2N}{N-\mu+2}}(\mathbb{R}^N)$.
By Lemma \ref{l2.0.0} we have $u,v\in L^s(\mathbb{R}^N)$ for all $s\in [2,\frac{2N^2}{(N-\mu)(N-2)})$. We use an equivalent form of Hardy-Littlewood-Sobolev inequality, which says that for any $r>\frac{N}{\mu}$, if $w\in L^{\frac{Nr}{(N-\mu)r+N}}(\mathbb{R}^N)$ there holds
\begin{equation}\label{2.0.6}|I_\mu\ast w|_r\leq C(r,N,\mu)|w|_{\frac{Nr}{(N-\mu)r+N}}.\end{equation}
Since $|u|^p\leq C(|u|^{2^*_\mu}+|u|^{2_{\mu,*}})$, and $u\in L^s(\mathbb{R}^N)$ for all $s\in [2,\frac{2N^2}{(N-\mu)(N-2)})$, we infer
$|u|^p, |u|^{2^*_\mu}\in L^{\frac{Nr}{(N-\mu)r+N}}$ for all $r>\frac{2N}{\mu}$.  Using (\ref{2.0.6}) and  letting $r\rightarrow+\infty$ there holds
$$\bigl|I_\mu\ast |u|^p\bigr|_\infty\leq C|u|^p_{\frac{N}{(N-\mu)}p}<+\infty,\quad\bigl|I_\mu\ast |u|^{2^*_\mu}\bigr|_\infty\leq C|u|^{2^*_\mu}_{\frac{N}{(N-\mu)}2^*_\mu}<+\infty.$$
Then from the first equation of system (\ref{2.0.5}) we know
 \begin{equation*}
-\Delta u+\mu_1u\leq C(|u|^{2_{\mu,*}-1}+|u|^{2^*_\mu-1}), \quad \text{in}\ \mathbb{R}^N.
\end{equation*}
Using the classical bootstrap method and regularity theory for Sobolev subcritical local problems in bounded domains, see \cite{Struwe}, we deduce that $u\in C^{2,\sigma}_{loc}$ for some $\sigma\in(0,1)$. Similarly, from the second equation of system (\ref{2.0.5}) we obtain that $v\in C^{2,\sigma}_{loc}$ for some $\sigma\in(0,1)$. Hence, $(u,v)$ is a smooth solution of (\ref{2.0.5}).
 \ \ \ \ \ $\Box$
\begin{lemma}\label{l3.4.0} Fix $\nu>0$.
Let $\{(u_n,v_n)\}\subset \mathcal{S}_r(a,b)$ be such that, for $n\rightarrow+\infty$,
\begin{equation}\label{2.3}J_\nu(u_n,v_n)\rightarrow c,\quad P_\nu(u_n,v_n)\rightarrow0,\end{equation}
for some fixed $c\in\mathbb{R}$.
Then $\{(u_n,v_n)\}$ is bounded in $H$.
\end{lemma}
{\bf Proof}: If
$4+\frac{4-2\mu}{N}<p+q<22^*_\mu$, then $\gamma_p+\gamma_q>2$. By (\ref{2.3}) we have
\begin{equation}\label{2.7.1}\aligned
c+o_n(1)\geq&\bigl(\frac{1}{2}-\frac{1}{22^*_\mu}\bigl)
\int_{\mathbb{R}^N}\bigl[(I_\mu\ast|u_n|^{2^*_\mu})|u_n|^{2^*_\mu}+
(I_\mu\ast|v_n|^{2^*_\mu})|v_n|^{2^*_\mu}\bigr]\\&+
\nu\frac{\gamma_p+\gamma_q-2}{2}
\int_{\mathbb{R}^N}(I_\mu\ast|u_n|^{p})|v_n|^{q}.
\endaligned\end{equation}
Since $\gamma_p+\gamma_q>2$, we get
$\int_{\mathbb{R}^N}(I_\mu\ast|u_n|^{2^*_\mu})|u_n|^{2^*_\mu}$, $
\int_{\mathbb{R}^N}(I_\mu\ast|v_n|^{2^*_\mu})|v_n|^{2^*_\mu}$ and $\int_{\mathbb{R}^N}(I_\mu\ast|u_n|^{p})|v_n|^{q}$ are bounded. Using
$P_\nu(u_n,v_n)\rightarrow0$, we have $\{|\nabla u_n|^2_2+|\nabla v_n|^2_2\}$ is bounded. Then $\{(u_n,v_n)\}$ is bounded in $H$.

If $22_{\mu,*}<p+q<4+\frac{4-2\mu}{N}$, then $0<\gamma_p+\gamma_q<2$.
From (\ref{1.14}) it follows that
\begin{equation}\label{2.7}\aligned
c+o_n(1)
\geq&\bigl(\frac{1}{2}-\frac{1}{22^*_\mu}\bigl)(|\nabla u_n|^2_2+|\nabla v_n|^2_2)\\&-\bigl(1-\frac{\gamma_p+\gamma_q}{22^*_\mu}\bigl)\nu C_{N,p,q}
(a^2+b^2)^{\frac{p+q-\gamma_p-\gamma_q}{2}}(|\nabla u_n|^2_2+|\nabla v_n|^2_2)^{\frac{\gamma_p+\gamma_q}{2}}.\endaligned
\end{equation}
Since $\gamma_p+\gamma_q<2$, we know $\{(u_n,v_n)\}$ is bounded in $H$.

If $p+q=4+\frac{4-2\mu}{N}$, then $\gamma_p+\gamma_q=2$.  For any $2_{\mu,*}<\bar{p}<p<2^*_\mu$, there exists some $\theta\in (0,1)$ such that $p=\bar{p}\theta +2^*_\mu(1-\theta)$, and then one easily has
$$\int_{\mathbb{R}^N}(I_\mu\ast|u_n|^{p})|u_n|^{p}\leq \left(\int_{\mathbb{R}^N}(I_\mu\ast|u_n|^{\bar{p}})|u_n|^{\bar{p}}\right)^\theta \left(\int_{\mathbb{R}^N}(I_\mu\ast|u_n|^{2^*_\mu})|u_n|^{2^*_\mu}\right)^{1-\theta}.$$
 Using (\ref{2.7.1}) with $\gamma_p+\gamma_q=2$, we get $\int_{\mathbb{R}^N}(I_\mu\ast|u_n|^{2^*_\mu})|u_n|^{2^*_\mu}$, $
\int_{\mathbb{R}^N}(I_\mu\ast|v_n|^{2^*_\mu})|v_n|^{2^*_\mu}$ are bounded.
So
$$\int_{\mathbb{R}^N}(I_\mu\ast|u_n|^{p})|u_n|^{p}\leq C_{\bar{p}} \left(\int_{\mathbb{R}^N}(I_\mu\ast|u_n|^{\bar{p}})|u_n|^{\bar{p}}\right)^\theta
\leq C_{N,a,\bar{p}} |\nabla u_n|_{2}^{2\gamma_{\bar{p}}\theta},$$
where we have used (\ref{1.12}).
Similarly, for any $2_{\mu,*}<\bar{q}<q<2^*_\mu$,
$$\int_{\mathbb{R}^N}(I_\mu\ast|v_n|^{q})|v_n|^{q}\leq C_{N,b,\bar{q}} |\nabla v_n|_{2}^{2{\gamma_{\bar{q}}}\theta'},$$
for some $\theta'<1$.
To show $\sup_{n}(|\nabla u_n|_{2}^{2}+|\nabla v_n|_{2}^{2})<+\infty$, we argue by contradiction and assume
$|\nabla u_n|_{2}^{2}+|\nabla v_n|_{2}^{2}\rightarrow +\infty$ after passing to a subsequence. In view of $\theta<1,\theta'<1$, we have
\begin{align*}
\int_{\mathbb{R}^N}(I_\mu\ast|u_n|^{p})|v_n|^{q}\leq & \left(\int_{\mathbb{R}^N}(I_\mu\ast|u_n|^{p})|u_n|^{p}\right)^{\frac{1}{2}}
\left(\int_{\mathbb{R}^N}(I_\mu\ast|v_n|^{q})|v_n|^{q}\right)^{\frac{1}{2}}\\
\leq & C_{N,a,b,\bar{p},\bar{q}} |\nabla u_n|_{2}^{\gamma_{\bar{p}}\theta} |\nabla v_n|_{2}^{\gamma_{\bar{q}}\theta'}\leq C_{N,a,b,\bar{p},\bar{q}} \left(|\nabla u_n|_{2}^{2}+|\nabla v_n|_{2}^{2}\right)^{\frac{\gamma_{\bar{p}}+\gamma_{\bar{q}}}{2}}.
\end{align*}
By the choice of $\bar{p}$ and $\bar{q}$, we get $\frac{\gamma_{\bar{p}}+\gamma_{\bar{q}}}{2}<\frac{\gamma_{p}+\gamma_{q}}{2}=1$.
Then
\begin{align*}
&c+o_n(1)+ \frac{1}{22^*_\mu}\int_{\mathbb{R}^N}\bigl[
(I_\mu*|u_n|^{2^*_\mu})|u_n|^{2^*_\mu}+(I_\mu*|v_n|^{2^*_\mu})
|v_n|^{2^*_\mu}\bigr]\\
=&\frac12(|\nabla u_n|^2_2+|\nabla v_n|^2_2)-\nu\int_{\mathbb{R}^N}
(I_\mu*|u_n|^{p})|v_n|^{q}\\
\geq& \frac12(|\nabla u_n|^2_2+|\nabla v_n|^2_2)-C_{N,\nu,a,b,\bar{p},\bar{q}} \left(|\nabla u_n|_{2}^{2}+|\nabla v_n|_{2}^{2}\right)^{\frac{\gamma_{\bar{p}}+\gamma_{\bar{q}}}{2}}\rightarrow +\infty,
\end{align*}
a contradiction. Hence, $\{(u_n,v_n)\}$
is bounded in $H$.\ \ \ \ \ $\Box$

\begin{Remark}\label{r1} In the case $p+q=4+\frac{4-2\mu}{N}$ in Lemma \ref{l3.4.0}, if additionally we assume $0<\nu<\nu'_0$ with $\nu'_0$ given in (\ref{niu0}) below, then the boundedness of $\{(u_n,v_n)\}$ in $H$ follows directly by using (\ref{2.7}) with $0<\nu<\nu'_0$. \end{Remark}

\begin{lemma}\label{l3.4.1}Assume $(u_n,v_n)\rightharpoonup(u,v)$ in $H$. Then
\begin{equation}\label{2.11}\int_{\mathbb{R}^N}
(I_\mu\ast|u_n|^{p})|v_n|^{q}=\int_{\mathbb{R}^N}
(I_\mu\ast|u|^{p})|v|^{q}+o_n(1).\end{equation}
\end{lemma}
{\bf Proof}:
Set $(\tilde{u}_n,\tilde{v}_n)=(u_n-u,v_n-v)$.
Note that
$$\aligned&\Bigl||u_n(x)|^p|v_n(y)|^q-|u(x)|^p|v(y)|^q\Bigr|\\
\leq&C|u_n(x)|^p|\tilde{v}_n(y)|[|v_n(y)|^{q-1}+|v(y)|^{q-1}]+C
|v(y)|^q|\tilde{u}_n(x)|[|u_n(x)|^{p-1}+|u(x)|^{p-1}],
\endaligned$$
for some $C>0$.
The Hardy-Littlewood-Sobolev inequality  implies
$$\aligned&\int_{\mathbb{R}^N}
[(I_\mu\ast|u_n|^{p})|v_n|^{q}-
(I_\mu\ast|u|^{p})|v|^{q}]\\ \leq& C|{u}_n|^p_{pr_0}|\tilde{v}_n|_{qr_0}
(|v_n|^{q-1}_{qr_0}+|v|^{q-1}_{qr_0})
+C|v|^q_{qr_0}|\tilde{u}_n|_{pr_0}
(|u_n|^{p-1}_{pr_0}+|u|^{p-1}_{pr_0})
,\endaligned$$
where $r_0=\frac{2N}{2N-\mu}$. Observe that $2_{\mu,*}<q,p<2^*_\mu$, we get $|\tilde{u}_n|_{pr_0}\rightarrow0$ and $|\tilde{v}_n|_{qr_0}\rightarrow0$. Hence, (\ref{2.11}) holds true.\ \ \ $\Box$

Now we recall a Liouville's type result.
\begin{lemma}\label{l2.5}(\cite[Lemma A.2]{IkomaANS})
Suppose $p\in(0,\frac{N}{N-2}]$ when $N\geq3$ and $p\in(0,\infty)$ when $N=1,2$. Let $u\in L^p(\mathbb{R}^N)$  be a smooth nonnegative function satisfying $-\Delta u\geq0$ in $\mathbb{R}^N$. Then $u\equiv0$.
\end{lemma}
We introduce a condition of $m_\nu(a,b)$ as follows
\begin{equation}\label{2.1}m_\nu(a_1,b_1)\leq m_\nu(a_2,b_2) \ \text{for any } 0<a_2\leq a_1, \ 0<b_2\leq b_1.\end{equation}
\begin{lemma}\label{l3.4} Assume $N\in\{3,4\}$ and (\ref{2.1}) is satisfied. Let $\{(u_n,v_n)\}\subset \mathcal{S}_r(a,b)$ satisfy the assumptions of Lemma \ref{l3.4.0} with
$u^-_n,v^-_n\rightarrow 0$ a.e.\ {in}\ $\mathbb{R}^N,$ and
\begin{equation}\label{2.5}
c\neq0\ \text{and}\ c<\frac{2^*_\mu-1}{22^*_\mu}S^{\frac{2^*_\mu}{2^*_\mu-1}}_{H,L}+
\min\{0,m_\nu(a,b)\},
\end{equation}
Moreover, assume
\begin{equation}\label{2.2}J'_\nu(u_n,v_n)+\lambda_{1,n}(u_n,0)+\lambda_{2,n}(0,v_n)\rightarrow0
\ \text{for some}\ \lambda_{1,n},\lambda_{2,n}\in\mathbb{R}.\end{equation}
Then there exist $(u,v)\in H_{rad}, u,v>0$ and $\lambda_1,\lambda_2>0$ such that $(u_n,v_n)\rightarrow(u,v)$ in $H$ and $(\lambda_{1,n},\lambda_{2,n})\rightarrow(\lambda_{1},\lambda_{2})$ in $\mathbb{R}^2$ after passing to a subsequence.
\end{lemma}
{\bf Proof}: From (\ref{2.2}) we obtain
\begin{equation}\label{3.6}\lambda_{1,n}=-a^{-2}\langle J'_\nu(u_n,v_n),(u_n,0)\rangle+o_n(1),\ \lambda_{2,n}=-b^{-2}\langle J'_\nu(u_n,v_n),(0,v_n)\rangle+o_n(1).\end{equation}
By Lemma \ref{l3.4.0}, we know $\{(u_n,v_n)\}$ is bounded in $H$.
Then $\{\lambda_{1,n}\}$ and $\{\lambda_{2,n}\}$ are bounded in $\mathbb{R}$ and so  there exist $(u,v)\in H_{rad}$, $\lambda_1,\lambda_2\in\mathbb{R}$ such that
$(u_n,v_n)\rightharpoonup (u,v)$ in $H$, $(u_n,v_n)\rightarrow(u,v)$ in $L^p(\mathbb{R}^N)\times L^p(\mathbb{R}^N)$, $(u_n,v_n)\rightarrow (u,v)$ a.e. in $\mathbb{R}^{2N}$ and $(\lambda_{1,n},\lambda_{2,n})\rightarrow(\lambda_1,\lambda_2)$ in $\mathbb{R}^2$. Furthermore, in view of (\ref{2.2}) and the fact that $u^-_n,v^-_n\rightarrow0$ a.e. in $\mathbb{R}^N$ we infer
\begin{equation}\label{2.8}
J'_\nu(u,v)+\lambda_1(u,0)+\lambda_2(0, v)=0,\
u\geq0,\ v\geq0,
\end{equation}
and so $P_\nu(u,v)=0$.

{\bf Claim 1}. $u\neq0$ and $v\neq0$. Otherwise, we may assume $u=0$. Then (\ref{2.8}) implies
\begin{equation*}-\Delta v+\lambda_2v=(I_\mu*|v|^{2^*_\mu})|v|^{2^*_\mu-2}v, \ v\geq0,\ v\in H^1(\mathbb{R}^N).\end{equation*}
If $v\neq0$, then $\lambda_2=0$. Thus, the maximum principle implies $v>0$. From Lemma \ref{l1.1} we know $v=\widetilde{U}_{\epsilon,\xi}$ for some $\epsilon>0$ and $\xi\in\mathbb{R}^N$. This is impossible since $\widetilde{U}_{\epsilon,\xi}\not\in L^2(\mathbb{R}^N)$ with $N\in\{3,4\}$. Therefore $v=0$.
We assume $|\nabla u_n|^2_2\rightarrow l_1\geq0$ and $|\nabla v_n|^2_2\rightarrow l_2\geq0$. If $l_1=l_2=0$, then by (\ref{1.14}) and the definition of $S_{H,L}$ we get
$J_\nu(u_n,v_n)\rightarrow0$, contradicts $c\neq0$. Hence,
$l_1+l_2>0$.
Since $v=0$, Lemma \ref{l3.4.1} implies
$\int_{\mathbb{R}^N}(I_\mu\ast|u_n|^{p})|v_n|^{q}\rightarrow0.$
From $P_\nu(u_n,v_n)\rightarrow0$ and the definition of $S_{H,L}$ we deduce
\begin{equation}\label{3.5}l_1+l_2\leq S^{-2^*_\mu}_{H,L}(l^{2^*_\mu}_1+l^{2^*_\mu}_2)\leq S^{-2^*_\mu}_{H,L}(l_1+l_2)^{2^*_\mu}.\end{equation}Then
 $l_1+l_2\geq S^{\frac{2^*_\mu}{2^*_\mu-1}}_{H,L}$.
Therefore, it is easy to see $c\geq \frac{2^*_\mu-1}{22^*_\mu} S^{\frac{2^*_\mu}{2^*_\mu-1}}_{H,L}$,
which contradicts (\ref{2.5}). Then $u\neq0$ and similarly $v\neq0$.

{\bf Claim 2}. $u_n\rightarrow u$ and $v_n\rightarrow v$ in $H^1(\mathbb{R}^N)$.

Indeed, firstly we show that $\lambda_1,\lambda_2>0$. If $\lambda_1\leq0$, then
$$-\Delta u=-\lambda_1 u+(I_\mu*|u|^{2^*_\mu})|u|^{2^*_\mu-2}u+
(I_\mu*|u|^{p})|v|^{q-2}v\geq0\ \ \text{in}\quad\mathbb{R}^N.$$
From Lemmas \ref{l2.5} and \ref{l2.0.1} we know $u\equiv0$ in $\mathbb{R}^N$, which is impossible. Hence $\lambda_1 >0$, and similarly $\lambda_2 >0$. Set $(\tilde{u}_n,\tilde{v}_n)=(u_n-u,v_n-v)$.
By the nonlocal Brezis-Lieb lemma we get
\begin{equation}\label{2.11.0}\int_{\mathbb{R}^N}(I_\mu\ast|u_n|^{2^*_\mu})|u_n|^{2^*_\mu}=
\int_{\mathbb{R}^N}(I_\mu\ast|u|^{2^*_\mu})|u|^{2^*_\mu}+
\int_{\mathbb{R}^N}(I_\mu\ast|\tilde{u}_n|^{2^*_\mu})
|\tilde{u}_n|^{2^*_\mu}+o_n(1),\end{equation}
and similar decomposition holds for $\{v_n\}$. Since $P_\nu(u_n,v_n)\rightarrow0$, $P(u,v)=0$,  by (\ref{2.11}) we have
\begin{equation}\label{2.12}\aligned|\nabla \tilde{u}_n|^2_2+|\nabla \tilde{v}_n|^2_2=\int_{\mathbb{R}^N}(I_\mu\ast|\tilde{u}_n|^{2^*_\mu})|\tilde{u}_n|^{2^*_\mu}
+\int_{\mathbb{R}^N}
(I_\mu\ast|\tilde{v}_n|^{2^*_\mu})|\tilde{v}_n|^{2^*_\mu}+o_n(1).
\endaligned\end{equation}
Taking same arguments as in Claim 1 we infer either $|\nabla \tilde{u}_n|^2_2+|\nabla \tilde{v}_n|^2_2\rightarrow0$ or $|\nabla \tilde{u}_n|^2_2+|\nabla \tilde{v}_n|^2_2\rightarrow l\geq S^{\frac{2^*_\mu}{2^*_\mu-1}}_{H,L}.$
If the latter case occurs, noting that $0<|u|_2\leq a$, $0<|v|_2\leq b$, from (\ref{2.12}) it follows that
$$\aligned
c+o_n(1)=&J_\nu(u_n,v_n)=J_\nu(u,v)+
\frac12(|\nabla \tilde{u}_n|^2_2+|\nabla \tilde{v}_n|^2_2)\\&-\frac1{22^*_\mu}\int_{\mathbb{R}^N}
\bigl[(I_\mu*|\tilde{u}_n|^{2^*_\mu})
|\tilde{u}_n|^{2^*_\mu}+(I_\mu*|\tilde{v}_n|^{2^*_\mu})
|\tilde{v}_n|^{2^*_\mu}\bigr]+o_n(1)\\
\geq&m_\nu(|u|_2,|v|_2)+
\bigl(\frac12-\frac1{22^*_\mu}\bigr)(|\nabla \tilde{u}_n|^2_2+|\nabla \tilde{v}_n|^2_2)+o_n(1)\\
\geq& m_\nu(a,b)+\frac{2^*_\mu-1}{22^*_\mu}S^{\frac{2^*_\mu}{2^*_\mu-1}}_{H,L}+o_n(1).
\endaligned$$
 This contradicts (\ref{2.5}). Hence, $|\nabla \tilde{u}_n|^2_2+|\nabla \tilde{v}_n|^2_2\rightarrow0$. Testing (\ref{2.2}) and (\ref{2.8}) with $(u_n,v_n)$ and $(u,v)$ respectively and subtracting,
together with (\ref{2.11}) and (\ref{2.11.0}) we deduce
\begin{equation*}
|\nabla \tilde{u}_n|^2_2+|\nabla \tilde{v}_n|^2_2+\lambda_1|\tilde{u}_n|^2_2+
\lambda_2|\tilde{v}_n|^2_2=\int_{\mathbb{R}^N}
\bigl[(I_\mu\ast|\tilde{u}_n|^{2^*_\mu})|\tilde{u}_n|^{2^*_\mu}
+(I_\mu\ast|\tilde{v}_n|^{2^*_\mu})|\tilde{v}_n|^{2^*_\mu}\bigr]+o_n(1).
\end{equation*}
Taking into account (\ref{2.12}), $|\nabla \tilde{u}_n|^2_2+|\nabla \tilde{v}_n|^2_2\rightarrow0$ and $\lambda_1,\lambda_2>0$, we know $\tilde{u}_n\rightarrow0$ and $\tilde{v}_n\rightarrow0$ in $H^1(\mathbb{R}^N)$. \ \ \ \ \ $\Box$
\section{Existence and nonexistence results in the case $\nu>0$}
\renewcommand{\theequation}{3.\arabic{equation}}
In this section we fix $\nu>0$ and give several existence and non-existence results by distinguishing three cases: $L^2$-subcritical case: $p+q<4+\frac{4-2\mu}{N}$; $L^2$-critical case: $p+q=4+\frac{4-2\mu}{N}$; $L^2$-supercritical case: $p+q>4+\frac{4-2\mu}{N}$.
\subsection{$L^2$-subcritical case}
\begin{lemma}\label{l2.2.0}Let $0<\nu<\nu_0$ with
\begin{equation}\label{2.13}\nu_0:=\frac{1}{\gamma_p+\gamma_q}
\frac{S^{\frac{(2-\gamma_p-\gamma_q)2^*_\mu}{22^*_\mu-2}}_{H,L}}{C_{N,p,q}}
\frac{(22^*_\mu-2)(2-\gamma_p-\gamma_q)
^{\frac{2-\gamma_p-\gamma_q}{22^*_\mu-2}}}{(22^*_\mu-\gamma_p-\gamma_q)
^{\frac{22^*_\mu-\gamma_p-\gamma_q}{22^*_\mu-2}}(a^2+b^2)^{\frac{p+q-\gamma_p-\gamma_q}{2}}},\end{equation}
where $C_{N,p,q}$ is given in (\ref{1.14}).
Then $\mathcal{P}^0_\nu(a,b)=\emptyset$ and $\mathcal{P}_\nu(a,b)$ is a submanifold of $H$.\end{lemma}
{\bf Proof}: Firstly we show $\mathcal{P}^0_\nu(a,b)=\emptyset$ and argue by contradiction we assume $(u,v)\in \mathcal{P}^0_\nu(a,b)$. In view of (\ref{1.8.0})
, $(\Psi^\nu_{(u,v)})''(0)=0$, (\ref{1.14})
and $\gamma_p+\gamma_q<2$, we obtain
\begin{equation}\label{2.14}|\nabla u|^2_2+|\nabla v|^2_2\leq\Bigl[\frac{22^*_\mu-\gamma_p-\gamma_q}{22^*_\mu-2}\nu
(\gamma_p+\gamma_q)C_{N,p,q}(a^2+b^2)^{\frac{p+q-\gamma_p-\gamma_q}{2}}\Bigr]
^{\frac{2}{2-\gamma_p-\gamma_q}}.\end{equation}
Moreover, as (\ref{3.5}) we have
\begin{equation*}\aligned&\frac{2-\gamma_p-\gamma_q}{22^*_\mu-\gamma_p-\gamma_q}(|\nabla u|^2_2+|\nabla v|^2_2)\leq S^{-2^*_\mu}_{H,L}(|\nabla u|^{2}_2+|\nabla v|^{2}_2)^{2^*_\mu},\endaligned\end{equation*} contradicts (\ref{2.14}) due to $\nu<\nu_0$.

Next we show $\mathcal{P}_\nu(a,b)$ is a submanifold of $H$. Note that $$\mathcal{P}_\nu(a,b)=\{(u,v)\in H:P(u,v)=0, G(u)=0, F(v)=0\},$$
where $G(u)=a^2-|u|^2_2$ and $F(v)=b^2-|v|^2_2$. It suffices to show
$d(P,G,F)$: $H\rightarrow\mathbb{R}^3$ is surjective.
Otherwise, by the  independence of $dG(u)$ and $dF(v)$, there must
be that $dP(u,v)$ is a linear combination of $dG(u)$ and $dF(v)$, i.e. there exist $\mu_1,\mu_2\in\mathbb{R}$, such that $(u,v)$ is a weak solution of
\begin{equation*}\aligned
\left\{ \begin{array}{lll}
-\Delta u+\mu_1 u=2^*_\mu(I_\mu\ast|u|^{2^*_\mu})
|u|^{2^*_\mu-2}u+\nu p\frac{\gamma_p+\gamma_q}{2}(I_\mu\ast|v|^{q})
|u|^{p-2}u,\\
-\Delta v+\mu_2 v=2^*_\mu(I_\mu\ast|v|^{2^*_\mu})
|v|^{2^*_\mu-2}v+\nu q\frac{\gamma_p+\gamma_q}{2}(I_\mu\ast|u|^{p})
|v|^{q-2}v, \\
|u|^2_2=a^2, |v|^2_2=b^2.
\end{array}\right.\endaligned
\end{equation*}
Testing the above system with $(u,v)$ and combining with Pohozaev identity we deduce $(\Psi^\nu_{(u,v)})''(0)=0$. Then $(u,v)\in \mathcal{P}^0_\nu({a,b})$, contradicts $\mathcal{P}^0_\nu(a,b)=\emptyset$. \ \ \ \ \ $\Box$

Using the definition of $S_{H,L}$ and (\ref{1.14}) we get
$$\aligned J_\nu(u,v)\geq\frac12(|\nabla u|^2_2+|\nabla v|^2_2)-A(|\nabla u|^2_2+|\nabla v|^2_2)^{\frac{\gamma_p+\gamma_q}{2}}-B(|\nabla u|^2_2+|\nabla v|^2_2)^{2^*_\mu},\endaligned$$
where
\begin{equation}\label{2.15}A=\nu C_{N,p,q}
(a^2+b^2)^{\frac{p+q-\gamma_p-\gamma_q}{2}},\quad B=({22^*_\mu})^{-1}S^{-2^*_\mu}_{H,L}.\end{equation}
For $\rho>0$, define
\begin{equation}\label{2.16}h(\rho)=\frac12\rho^2-A\rho^{\gamma_p+\gamma_q}-B\rho^{22^*_\mu}, \quad g(\rho)=\rho^{2-\gamma_p-\gamma_q}-22^*_\mu B\rho^{22^*_\mu -\gamma_p-\gamma_q}.\end{equation}
Setting
\begin{equation*}\bar{\rho}=\Bigl(\frac{2-\gamma_p-\gamma_q}{22^*_\mu-\gamma_p-\gamma_q}
S^{2^*_\mu}_{H,L}\Bigr)^{\frac1{22^*_\mu-2}},\end{equation*}
it is easy to check that $g(\rho)$ is increasing in $(0,\bar{\rho})$, and decreasing in
$(\bar{\rho},+\infty)$. By direct computations, the assumption $0<\nu<\nu_0$ with $\nu_0$ defined in (\ref{2.13})
implies $g(\bar{\rho})>(\gamma_p+\gamma_q)A$ and $h'(\bar{\rho})>0$. This means that $h(\bar{\rho})$ has exactly two critical points
$0<\rho_1<\bar{\rho}<\rho_2$ with
$h(\rho_1)=\min_{0<\rho<\bar{\rho}}h(\rho)<0$, $ h(\rho_2)=\max_{\rho>0}h(\rho)>0.$
Moreover, there exist $R_1>R_0>0$ such that $h(R_0)=h(R_1)=0$ and $h(\rho)>0$ if and only if $\rho\in (R_0,R_1)$.  Taking the standard argument as in \cite[Lemma 5.2]{SoaveJFA}, we have the following lemma.

\begin{lemma}\label{l2.2}Let $0<\nu<\nu_0$. For every $(u,v)\in \mathcal{S}(a,b)$, $\Psi^\nu_{(u,v)}(t)$ has exactly two critical
points $t_\nu(u,v)<s_\nu(u,v)$ and two zeros $c_\nu(u,v)<d_\nu(u,v)$ with $t_\nu(u,v)<c_\nu(u,v)<s_\nu(u,v)<d_\nu(u,v)$. Moreover,

\noindent(i)  $t_\nu\star(u,v)\in \mathcal{P}^+_\nu(a,b)$, $s_\nu\star(u,v)\in \mathcal{P}^-_\nu(a,b)$, and if $t\star(u,v)\in \mathcal{P}_\nu(a,b)$, then either
$t=t_\nu(u,v)$ or $t=s_\nu(u,v)$.\\
\noindent(ii) If $t\leq c_\nu(u,v)$, then $|\nabla (t\star u)|^2_2+|\nabla (t\star v)|^2_2\leq R^2_0$,  and
\begin{equation*}\aligned
\Psi^\nu_{(u,v)}(t_\nu(u,v))=\min\bigl\{\Psi^\nu_{(u,v)}(t):t\in\mathbb{R},
|\nabla (t\star u)|^2_2+|\nabla (t\star v)|^2_2\leq R^2_0\bigr\}<0.
\endaligned\end{equation*}
Moreover, $\Psi^\nu_{(u,v)}(s_\nu(u,v))=\max_{t\in\mathbb{R}}\Psi^\nu_{(u,v)}(t)$ and $\Psi^\nu_{(u,v)}(t)$ is decreasing on $(s_\nu(u,v),+\infty)$.\\
\noindent(iii)  The maps $(u,v)\mapsto t_\nu(u,v)$ and $(u,v)\mapsto s_\nu(u,v)$ are of class $C^1$.
\end{lemma}

For $R>0$, $a,b>0$, we define
\begin{equation}\label{2.2.1}A_R(a,b):=\bigl\{(u,v)\in \mathcal{S}(a,b):|\nabla u|^2_2+|\nabla v|^2_2<R^2\bigr\}.\end{equation}
\begin{lemma}\label{l2.3}If $0<\nu<\nu_0$, then  $m_\nu(a,b)=\inf_{\mathcal{P}^+_\nu(a,b)}J_\nu=\inf_{A_{R_0}(a,b)}J_\nu<0$. Moreover, for some $r>0$ small \begin{equation}\label{4.5.0}m_\nu(a,b)<\inf_{\overline{A_{R_0}(a,b)}\backslash A_{R_0-r}(a,b)}J_\nu.\end{equation}
In addition, $m_\nu(a,b)$ satisfies (\ref{2.1}). \end{lemma}
{\bf Proof}: From Lemma \ref{l2.2} we know $\mathcal{P}^+_\nu(a,b)\subset A_{R_0}(a,b)$ and
$m_\nu(a,b)
=\inf_{\mathcal{P}^+_\nu(a,b)}J_\nu<0.$
Then $m_\nu(a,b)\geq\inf_{A_{R_0}(a,b)}J_\nu$. On the other hand, for any $(u,v)\in A_{R_0}(a,b)$, there holds
$$m_\nu(a,b)\leq J_\nu(t_\nu(u,v)\star(u,v))\leq J_\nu(u,v).$$
Thus, $m_\nu(a,b)=\inf_{A_{R_0}(a,b)}J_\nu$. Moreover, by the continuity of $h$, there is $r>0$ such that $h(t)\geq\frac{m_\nu(a,b)}{2}$ if $t\in[R_0-r, R_0]$. Observe that
$J_\nu(u,v)\geq h((|\nabla u|^2_2+|\nabla v|^2_2)^{\frac12})$.
Then (\ref{4.5.0}) follows.
With minor modifications by following the strategy in \cite[Lemma 2.3]{Bartsch-Li-Zou}, where local Schr\"{o}dinger system was considered, we can infer that $m_\nu(a,b)$ satisfies (\ref{2.1}).\ \ \ \ $\Box$

{\bf Proof of Theorem \ref{th1.1} (i), (iii) and  Theorem \ref{th1.4}} Let $\{(\bar{u}_n,\bar{v}_n)\}\subset A_{R_0}(a,b)$ be a minimizing sequence of $m_\nu(a,b)$ and we may assume that $\bar{u}_n,\bar{v}_n$ are radially decreasing (if this is not the case, we can
replace $(\bar{u}_n,\bar{v}_n)$ with $(|\bar{u}_n|^*,|\bar{v}_n|^*)$, the Schwarz rearrangement of $(|\bar{u}_n|,|\bar{v}_n|)$, then $(|\bar{u}_n|^*,|\bar{v}_n|^*)\in A_{R_0}(a,b)$ and
satisfies $J_\nu(|\bar{u}_n|^*,|\bar{v}_n|^*)\leq J_\nu(\bar{u}_n,\bar{v}_n)$). By Lemma \ref{l2.2}, for  every $n$, there
exists $t_n$ such that $t_n\star(\bar{u}_n,\bar{v}_n)\in \mathcal{P}^+_\nu(a,b)\cap A_{R_0}(a,b)$ and
$$J_\nu(t_n\star(\bar{u}_n,\bar{v}_n))=\min\{J_\nu(t\star(\bar{u}_n,\bar{v}_n)):t\in\mathbb{R}, |\nabla \bar{u}_n|^2_2+|\nabla \bar{v}_n|^2_2<R^2_0\}\leq J_\nu(\bar{u}_n,\bar{v}_n).$$
Therefore, we can get a new minimizing sequence $\{(\hat{u}_n,\hat{v}_n)=t_n\star(\bar{u}_n,\bar{v}_n)\}\subset \mathcal{S}_r(a,b)\cap\mathcal{P}^+_\nu(a,b)$, which is radially
decreasing for every $n$. Using Lemma \ref{l2.3} we have $(\hat{u}_n,\hat{v}_n)\in A_{R_0-r}(a,b)$.
  Since $(|\hat{u}_n|, |\hat{v}_n|)\in \mathcal{P}_\nu(a,b)$ and $J_\nu(|\hat{v}_n|, |\hat{v}_n|)=J_\nu(\hat{u}_n, \hat{v}_n)$, we may also assume that $\hat{u}_n,\hat{v}_n$ are nonnegative. Thus, Ekeland varational principle implies that there is a radially symmetric Palais-Smale sequence $(u_n,v_n)$ for $J_\nu|_{\mathcal{S}(a,b)}$ satisfying $\|(u_n,v_n)-(\hat{u}_n,\hat{v}_n)\|_H\rightarrow0$ as $n\rightarrow\infty$. Then
$$P_\nu(u_n,v_n)=P_\nu(\hat{u}_n,\hat{v}_n)+o_n(1)=o_n(1)\quad\text{and}\ u^-_n,v^-_n\rightarrow0\ \text{a.e.\ in}\ \mathbb{R}^N.$$
Applying Lemma \ref{l3.4} with $c=m_\nu(a,b)$ we get $(u_n,v_n)\rightarrow(u,v)$ in $H$ and $(\lambda_{1,n},\lambda_{2,n})\rightarrow(\lambda_{1},\lambda_{2})$ with $\lambda_1,\lambda_2>0$. Hence, $(u,v)$ is a normalized ground state of (\ref{0.0}).\ \ \ \ $\Box$

\begin{lemma}\label{l3.3}For fixed $a,b>0$ and the given function $h$ in (\ref{2.16}), there exist $\nu_0=\nu_0(a,b)>0$ and $\rho_0=\rho_0(a,b)>0$ such that $h(\rho)>0$ if $\rho_0\leq\rho\leq2 \rho_0$
for all $0<\nu<\nu_0$ and so $R_0\leq\rho_0$ for all $0<\nu<\nu_0$.
\end{lemma}
{\bf Proof}:
For the function $h$ given in (\ref{2.16}) and constants $A,B$ given in (\ref{2.15}), we can choose $\rho_0$ small enough such that
$$B\rho^{22^*_\mu-2}\leq B(2\rho_0)^{22^*_\mu-2}\leq\frac18, \ \text{when}\ \rho\leq2\rho_0.$$
Note that $\gamma_p+\gamma_q<2$. Choose $\nu_0$ small enough such that for all $0<\nu<\nu_0$ there holds
$$A\rho^{\gamma_p+\gamma_q-2}
\leq\nu C_{N,p,q}(a^2+b^2)^{\frac{p+q-\gamma_p-\gamma_q}{2}}\rho^{\gamma_p+\gamma_q-2}_0\leq\frac18, \ \text{when}\ \rho_0\leq\rho.$$
Then $h(\rho)\geq\frac14\rho^2_0$ if $\rho_0\leq\rho\leq2\rho_0$ for $0<\nu<\nu_0$.
Note that in the paragraph before Lemma \ref{l2.2}, $h(R_0)=h(R_1)=0$ and $h(\rho)=0$ if and only if $\rho\in(R_0,R_1)$. Hence, $R_0\leq\rho_0$ for $0<\nu<\nu_0$.\ \ \ \ $\Box$

{\bf Proof of Theorem \ref{th1.1} (ii), (iii) and  Theorem \ref{th1.4}}: As in the proof of Theorem \ref{th1.1} (i), there is a radially symmetric Palais-Smale sequence $(u_n,v_n)$ for $J_\nu|_{\mathcal{S}(a,b)}$ satisfying $\|(u_n,v_n)-(\hat{u}_n,\hat{v}_n)\|_H\rightarrow0$ with $(\hat{u}_n,\hat{v}_n)\in A_{R_0-r}(a,b)$,
$$J_\nu(u_n,v_n)\rightarrow m_\nu(a,b), \ P_\nu(u_n,v_n)=o_n(1)\quad\text{and}\ u^-_n,v^-_n\rightarrow0\ \text{a.e.\ in}\ \mathbb{R}^N,$$
and
$$J'_\nu(u_n,v_n)+\lambda_{1,n}(u_n,0)+\lambda_{2,n}(0,v_n)\rightarrow0,
\quad\lambda_{1,n},\lambda_{2,n}\in\mathbb{R}.$$
Applying Lemma \ref{l3.4.0} we get $\{(u_n,v_n)\}$ is bounded in $H$ and so $\lambda_{1,n},\lambda_{2,n}$ are bounded.
Assume  $(u_n,v_n)\rightarrow(u,v)$ in $H$ and $(\lambda_{1,n},\lambda_{2,n})\rightarrow(\lambda_1,\lambda_2)$. Then $u,v\geq0$ a.e. in $\mathbb{R}^N$ and $(u,v)$ satisfies
\begin{equation}\label{3.4}J'_\nu(u,v)+\lambda_{1}(u,0)+\lambda_{2}(0,v)=0.\end{equation}
 Therefore $P_\nu(u,v)=0$. From Lemma \ref{l3.4.1} and weak lower semi-continuity of the norm it follows that
 $$\aligned
 J_\nu(u,v)-\frac1{22^*_\mu}P_\nu(u,v)=&
 \frac{2^*_\mu-1}{22^*_\mu}(|\nabla u|^2_2+|\nabla v|^2_2)+\frac{\gamma_p+\gamma_q-22^*_\mu}
 {22^*_\mu}\nu
 \int_{\mathbb{R}^N}(I_\mu\ast|u|^p)|v|^q\\
 \leq&\frac{2^*_\mu-1}{22^*_\mu}(|\nabla u_n|^2_2+|\nabla v_n|^2_2)+\frac{\gamma_p+\gamma_q-22^*_\mu}
 {22^*_\mu}\nu
 \int_{\mathbb{R}^N}(I_\mu\ast|u_n|^p)|v_n|^q\\
 =&J_\nu(u_n,v_n)-\frac1{22^*_\mu}P_\nu(u_n,v_n)=m_\nu(a,b)+o_n(1).\endaligned$$
Then $J_\nu(u,v)\leq m_\nu(a,b)<0$. Therefore, $(u,v)\neq(0,0)$.

{\bf Claim 1}. $u\neq0,v\neq0$. Otherwise, by above $(u,v)\neq(0,0)$, we may assume $u\neq0$ and $v=0$. Then $|\nabla u_n|^2_2\rightarrow l_1>0$ and $|\nabla v_n|^2_2\rightarrow l_2\geq0$. As in the proof of Lemma \ref{l3.4}, we get $l_1+l_2\geq S^{\frac{2^*_\mu}{2^*_\mu-1}}_{H,L}$ and so $m_\nu(a,b)\geq(\frac12-\frac1{22^*_\mu})S^{\frac{2^*_\mu}{2^*_\mu-1}}_{H,L}$. However, $m_\nu(a,b)<0$. This is a contradiction. Then $u\neq0,v\neq0$.

{\bf Claim 2}. $u_n\rightarrow u,v_n\rightarrow v$ in $H^1(\mathbb{R}^N)$.
 In fact, we firstly show that $\lambda_1,\lambda_2>0$. If $\lambda_1\leq0$, testing (\ref{3.4}) with $(u,0)$ we get
 $$|\nabla u|^2_2\geq \int_{\mathbb{R}^N}(I_\mu\ast|u|^{2^*_\mu})|u|^{2^*_\mu}+\nu p\int_{\mathbb{R}^N}(I_\mu\ast|u|^{p})|v|^{q}.$$
 Jointing the above inequality with $P_\nu(u,v)=0$ we obtain
\begin{equation}\label{3.7}|\nabla v|^2_2-\int_{\mathbb{R}^N}(I_\mu\ast|v|^{2^*_\mu})|v|^{2^*_\mu}\leq
\nu(\gamma_p+\gamma_q-p)\int_{\mathbb{R}^N}(I_\mu\ast|u|^{p})|v|^{q}.\end{equation}
 Since $(\frac{N}{2}-1)p+\frac{N}{2}q\leq 2N-\mu$, the coefficient of $\int_{\mathbb{R}^N}(I_\mu\ast|u|^{p})|v|^{q}$ is non-positive. By taking, if necessary $\nu_0$ and $\rho_0$ smaller (and thus $R_0$ smaller) than in Lemma \ref{l3.3} we get
 $$|\nabla v|^2_2-\int_{\mathbb{R}^N}(I_\mu\ast|v|^{2^*_\mu})|v|^{2^*_\mu}
 \geq|\nabla v|^2_2-C|\nabla v|^{22^*_\mu}_2>0,$$
 due to the definition of $S_{H,L}$. Therefore, there is a contradiction from (\ref{3.7}). Then $\lambda_1>0$. Similarly, $\lambda_2>0$. Below as in the proof of Lemma \ref{l3.4} we have $(u_n,v_n)\rightarrow(u,v)$ in $H$. Hence, $(u,v)$ is a normalized ground state of (\ref{0.0}).\ \ \ \ $\Box$

\subsection{$L^2$-supercritical case}

\begin{lemma}\label{l2.4}Let $\nu>0$. Then $\mathcal{P}_\nu(a,b)$ is a smooth manifold of $H$ and for every $(u,v)\in \mathcal{S}(a,b)$, there exists a unique $t_\nu(u,v)\in\mathbb{R}$ such that $t_\nu(u,v)\star(u,v)\in \mathcal{P}_\nu(a,b)$ and $\Psi^\nu_{(u,v)}(t_\nu(u,v))=
\max_{t\in\mathbb{R}}\Psi^\nu_{(u,v)}(t)>0$. Moreover
\\
\noindent(i) $\mathcal{P}_\nu(a,b)=\mathcal{P}^-_\nu(a,b)$, $t_\nu(u,v)<0$ if and only if
$P_\nu(u,v)<0$.\\
\noindent(ii) $\Psi^\nu_{(u,v)}(t)$ is decreasing and concave on $(t_\nu(u,v),+\infty)$.\\
\noindent(iii)  The map $(u,v)\mapsto t_\nu(u,v)$ is of class $C^1$.
\end{lemma}
{\bf Proof}: Taking similar arguments as in Lemma \ref{l2.2.0}, we know $\mathcal{P}_\nu(a,b)$ is a smooth manifold of $H$.  For any $(u,v)\in \mathcal{P}_\nu(a,b)\backslash{\mathcal{P}^-_\nu(a,b)}$ we deduce
$$(22^{*}_\mu-2)\int_{\mathbb{R}^N}[(I_\mu\ast|u|^{2^*_\mu})
|u|^{2^*_\mu}+(I_\mu\ast|v|^{2^*_\mu})|v|^{2^*_\mu}]\leq\nu (\gamma_p+\gamma_q)(2-\gamma_p-\gamma_q)\int_{\mathbb{R}^N}(I_\mu\ast|u|^{p})
|v|^{q},$$
which is impossible because $\gamma_p+\gamma_q>2$. Then $\mathcal{P}_\nu(a,b)={\mathcal{P}^-_\nu(a,b)}$. From the definition of $\Psi^\nu_{(u,v)}(t)$ and $\gamma_p+\gamma_q>2$, we get $\Psi^\nu_{(u,v)}(t)\rightarrow0^+$ as $t\rightarrow-\infty$ and $\Psi^\nu_{(u,v)}(t)\rightarrow-\infty$ as $t\rightarrow+\infty$. Hence $\Psi^\nu_{(u,v)}(t)$ has a global maximum point at positive level. Furthermore, $e^{-2t}(\Psi^\nu_{(u,v)})'(t)$ is decreasing in $t$, $(\Psi^\nu_{(u,v)})'(t)$ has exactly one zero point. Note that $t_\nu(u,v)\star(u,v)\in \mathcal{P}^-_\nu(a,b)$. Then $(\Psi^\nu_{(u,v)})''(t_\nu(u,v))=(\Psi^\nu_{t_\nu(u,v)\star(u,v)})''(0)<0$ and $(\Psi^\nu_{(u,v)})''(t)<0$ if $t>t_\nu(u,v)$ since  $\bigl(\Psi^\nu_{(u,v)}\bigr)''(t)e^{-2t}$ is decreasing in $t$. Let $\Phi(t,(u,v))=\bigl(\Psi^\nu_{(u,v)}\bigr)'(t)$. Then
$\Phi(t_\nu(u,v),(u,v))=0$ and $$\partial_t\Phi(t,(u,v))\bigl|_{t=t_\nu(u,v)}
=\bigl(\Psi^\nu_{(u,v)}\bigr)''(t_\nu(u,v))<0.$$  Applying the implicit function theorem, the map $(u,v)\rightarrow t_\nu(u,v)$ is of class $C^1$. \ \ \ \ $\Box$
\begin{Remark}\label{r3.1}  From Lemma \ref{l2.4}, it is easy to see that the minima $m_\nu(a,b)$ and $m_{r,\nu}(a,b)$ have minimax characterization:
\begin{equation}\label{4.11.0}m_\nu(a,b)
=\inf_{\mathcal{S}(a,b)}\max_{t\in\mathbb{R}}J_\nu(t\star(u,v)),\quad m_{r,\nu}(a,b)
=\inf_{\mathcal{S}_r(a,b)}\max_{t\in\mathbb{R}}J_\nu(t\star(u,v)).\end{equation}\end{Remark}
\begin{lemma}\label{l4.5.0}$m_{r,\nu}(a,b)\geq m_\nu(a,b)>0$.\end{lemma}
{\bf Proof}:
Since $p+q>4+\frac{4-2\mu}{N}$, we have $\gamma_p+\gamma_q>2$. For $(u,v)\in \mathcal{P}_\nu(a,b)$, using the definition of $S_{H,L}$ and (\ref{1.14}) there holds
$$\aligned |\nabla u|^2_2+|\nabla v|^2_2
&\leq C(|\nabla u|^2_2+|\nabla v|^2_2)^{2^*_\mu}+C(|\nabla u|^2_2+|\nabla v|^2_2)^{\frac{\gamma_p+\gamma_q}{2}},\endaligned$$
which implies $\inf_{\mathcal{P}_\nu(a,b)}(|\nabla u|^2_2+|\nabla v|^2_2)>0$ since $\gamma_p+\gamma_q>2$. So
\begin{equation}\label{4.6}\aligned m_\nu(a,b)
=&\inf_{\mathcal{P}_\nu(a,b)}\Bigl\{\frac{\gamma_p+\gamma_q-2}{2(\gamma_p+\gamma_q)}(|\nabla u|^2_2+|\nabla v|^2_2)\\&+\frac{22^*_\mu-\gamma_p-\gamma_q}{22^*_\mu(\gamma_p+\gamma_q)}\int_{\mathbb{R}^N}[(I_\mu\ast|u|^{2^*_\mu})
|u|^{2^*_\mu}+(I_\mu\ast|v|^{2^*_\mu})|v|^{2^*_\mu}]\Bigr\}>0.
\endaligned\end{equation}
This ends the proof.
\ \ \ \ $\Box$

As \cite[Lemma 4.5]{Yeyang}, we have the following result.
\begin{lemma}\label{l4.5}
 There exists $k>0$ small enough, such that $\inf_{\overline{A_k(a,b)}}J_\nu>0$, $\inf_{\overline{A_k(a,b)}}P_\nu>0$ and $\sup_{\overline{A_k(a,b)}}J_\nu<m_\nu(a,b)$, where $A_k(a,b)$ is defined in (\ref{2.2.1}).
\end{lemma}

Following the strategy introduced in \cite{Jeanjean}, we consider the functional $\hat{J}_\nu(t,(u,v)):\mathbb{R}\times H\rightarrow\mathbb{R}$ defined by $\hat{J}_\nu(t,(u,v)):=J_\nu(t\star(u,v))$. Clearly, $\hat{J}_\nu$ is of class $C^1$.  Let $J^c_\nu$ be the closed sublevel set
$\{(u,v)\in \mathcal{S}(a,b):J_\nu(u,v)\leq c\}$. Let
$$\sigma(a,b,\nu):=\inf_{\gamma\in \Gamma_\nu}\max_{(t,(u,v))\in \gamma([0,1])}\hat{J}_\nu(t,(u,v)),$$
where
$$\aligned\Gamma_\nu=&\bigl\{\gamma=(\alpha,\beta)\in C([0,1],\mathbb{R}\times\mathcal{S}_r(a,b)): \gamma(0)\in\{0\}\times A_k(a,b), \gamma(1)\in\{0\}\times J^0_\nu\bigr\},\endaligned$$
with $k$ given in Lemma \ref{l4.5}. Similar to \cite[Lemma 4.7]{Yeyang} and \cite[Lemma 4.8]{Yeyang}, we have $\sigma(a,b,\nu)=m_{r,\nu}(a,b),$ as well as the following result.
\begin{lemma}\label{l2.6}Let $\nu>0$. There exists a Palais-Smale sequence $(u_n,v_n)\subset \mathcal{S}_r(a,b)$ for $J_\nu|_{\mathcal{S}(a,b)}$ at the level $m_{r,\nu}(a,b)=\sigma(a,b,\nu)$ with $P_\nu(u_n,v_n)\rightarrow0$ and $u^-_n$, $v^-_n\rightarrow0$ a.e. in $\mathbb{R}^N$ as $n\rightarrow\infty$.
 \end{lemma}

Below we give some properties of $m_{r,\nu}(a,b)$ and $m_{\nu}(a,b)$ to recover the
compactness of Palais-Smale sequence.

\begin{lemma}\label{l2.7} (i) $m_{\nu}(a,b)$ satisfies (\ref{2.1}), $m_{r,\nu}(a,b)$ satisfies (\ref{2.1}) with $m_{\nu}(a,b)$ replacing by $m_{r,\nu}(a,b)$, $m_{r,\nu}(a,b)$ and $m_{\nu}(a,b)$ are decreasing in $\nu\in(0,+\infty)$.\\ \noindent(ii) $\lim_{\nu\rightarrow+\infty}m_{r,\nu}(a,b)=
 \lim_{\nu\rightarrow+\infty}m_{\nu}(a,b)=0^+$.
 \end{lemma}
{\bf Proof}: (i) With minor modifications of the argument in \cite[Lemma 3.2]{Bartsch-Li-Zou}, where Schr\"{o}dinger system with local nonlinearity was considered, we can show $m_{\nu}(a,b)$ satisfies (\ref{2.1}). Here we give the details for the reader's convenience. In fact, by the definition of $m_\nu(a_1,b_1)$, for any $\epsilon>0$, there exists $(u,v)\in \mathcal{P}_\nu(a_1,b_1)$ such that
 \begin{equation}\label{3.2}J_\nu(u,v)\leq m_\nu(a_1,b_1)+\epsilon.\end{equation}
 Let $\phi\in C^\infty_0(\mathbb{R}^N,[0,1])$ be radial and satisfy
$\phi(x)=0$ if $|x|\geq2$ and $\phi(x)=1$ if $|x|\leq1$.
Consider $u_\delta(x)=u(x)\phi(\delta x)$ and $v_\delta(x)=v(x)\phi(\delta x)$, where $\delta>0$. Clearly $({u}_\delta,{v}_\delta)\rightarrow(u,v)$ in $H$ as $\delta\rightarrow0^+$.
Let $\bar{u}_\delta=a_1|u_\delta|^{-1}_2u_\delta$ and $\bar{v}_\delta=b_1|v_\delta|^{-1}_2v_\delta$. Then $(\bar{u}_\delta,\bar{v}_\delta)\rightarrow (a_1|u|^{-1}_2u,b_1|v|^{-1}_2v)=(u,v)$ in $H$ as $\delta\rightarrow0^+$.
Lemma \ref{l2.4} (iii) implies $t_\nu(\bar{u}_\delta,\bar{v}_\delta)\star(\bar{u}_\delta,
\bar{v}_\delta)\rightarrow t_\nu(u,v)\star(u,v)$ in $H$ as $\delta\rightarrow0^+$.
Then by (\ref{3.2}) there exists $\delta>0$ small such that
\begin{equation}\label{3.3}J_\nu
(t_\nu(\bar{u}_\delta,\bar{v}_\delta)\star(\bar{u}_\delta,\bar{v}_\delta))\leq J_\nu(t_\nu(u,v)\star(u,v))+\epsilon\leq J_\nu(u,v)+\epsilon\leq m_\nu(a_1,b_1)+2\epsilon.\end{equation}
Consider
$$w_{a_2}=\sqrt{a^2_2-|\bar{u}_\delta|^2_2}|\varphi|^{-1}_2\varphi,\quad w_{b_2}=\sqrt{b^2_2-|\bar{v}_\delta|^2_2}|\varphi|^{-1}_2\varphi,$$
where $\varphi\in C^\infty_0(\mathbb{R}^N)$ satisfies $supp\varphi\subset \mathbb{R}^N\backslash{B_{\frac{4}{\delta}}}$.
For $\tau\leq0$, define
 $(\tilde{u}_\tau, \tilde{v}_\tau):=(\bar{u}_\delta+\tau\star w_{a_2},\bar{v}_\delta+\tau\star w_{b_2})$. Noting that
$$(supp(\bar{u}_\delta)\cup supp(\bar{v}_\delta))\cap(supp(\tau\star w_{a_2})\cup supp(\tau\star w_{b_2}))=\emptyset,$$
we know $(\tilde{u}_\tau, \tilde{v}_\tau)\in \mathcal{S}(a_2,b_2)$. Let $t_\tau=t_\nu(\tilde{u}_\tau, \tilde{v}_\tau)$ be such that ${P}_\nu(t_\tau\star(\tilde{u}_\tau, \tilde{v}_\tau))=0$. Then
$$\aligned&e^{-(22^*_\mu-2)t_\tau}(|\nabla\tilde{u}_\tau|^2_2+
|\nabla\tilde{v}_\tau|^2_2)+\nu(\gamma_p+\gamma_q)e^{-
(22^*_\mu-\gamma_p-\gamma_q)t_\tau}\int_{\mathbb{R}^N}
(I_\mu*|\tilde{u}_\tau|^p)|\tilde{v}_\tau|^q
\\=&\int_{\mathbb{R}^N}
[(I_\mu*|\tilde{u}_\tau|^{2^*_\mu})|\tilde{u}_\tau|^{2^*_\mu}+
(I_\mu*|\tilde{v}_\tau|^{2^*_\mu})|\tilde{v}_\tau|^{2^*_\mu}].\endaligned$$
Since $(\tilde{u}_\tau, \tilde{v}_\tau)\rightarrow(\bar{u}_\delta,\bar{v}_\delta)\neq(0,0)$ in $H$ as $\tau\rightarrow-\infty$, we have $\limsup_{\tau\rightarrow-\infty}t_\tau<+\infty$ and then $t_\tau+\tau\rightarrow-\infty$ as $\tau\rightarrow-\infty$.  Thus, for $\tau<-1$ small,
$J_\nu((t_\tau+\tau)\star(w_{a_2},w_{b_2}))<\epsilon$,
which combining with (\ref{3.3}) we obtain
$$\aligned J_\nu(t_\tau\star(\tilde{u}_\tau, \tilde{v}_\tau))=J_\nu(t_\tau\star(\bar{u}_\delta,\bar{v}_\delta))
+J_\nu((t_\tau+\tau)
\star (w_{a_2},w_{b_2}))
\leq m_\nu(a_1,b_1)+3\epsilon.\endaligned$$
Note that  $m_\nu(a_2,b_2)\leq J_\nu(t_\tau\star(\tilde{u}_\tau, \tilde{v}_\tau))$.
From the arbitrariness
 of $\epsilon>0$, we know $m_\nu(a,b)$ satisfies (\ref{2.1}).
Taking the above arguments with $m_{\nu}(a,b)$,  $\mathcal{S}(a,b)$ and $\mathcal{P}_\nu(a,b)$ replacing by $m_{r,\nu}(a,b)$, $\mathcal{S}_r(a,b)$ and $\mathcal{P}_\nu(a,b)\cap H_{rad}$ respectively, and with $\varphi\in C^\infty_0(\mathbb{R}^N)$ chosen to be radial, we can deduce that $m_{r,\nu}(a_2,b_2)\leq m_{r,\nu}(a_1,b_1)$. In addition, using (\ref{4.11.0}), $m_{r,\nu}(a,b)$ and $m_{\nu}(a,b)$ are decreasing with respect to $\nu\in(0,+\infty)$.

(ii) For simplicity, denote
\begin{equation}\label{3.6.2}K(u)=\frac{1}2|\nabla u|^2_2-\frac{1}{22^*_\mu}\int_{\mathbb{R}^N}
(I_\mu\ast|u|^{2^*_\mu})|u|^{2^*_\mu},\quad \forall u\in H^1(\mathbb{R}^N).\end{equation}
It is easy to see that
\begin{equation}\label{3.6.1}\max_{t\in\mathbb{R}}K(t\star u)=\frac{2^*_\mu-1}{22^*_\mu}\frac{|\nabla u|^{\frac{22^*_\mu}{2^*_\mu-1}}_2}{\bigl(\int_{\mathbb{R}^N}(I_\mu\ast|u|^{2^*_\mu})|u|^{2^*_\mu}\bigr)^{\frac1{2^*_\mu-1}}}
,\quad \forall u\in H^1(\mathbb{R}^N)\backslash\{0\}.\end{equation}
Choose $\phi\in C^\infty_0(\mathbb{R}^N)$ be radial with $|\phi|_2\leq\min\{a,b\}$. By (i) we have
\begin{equation*}\aligned
m_{r,\nu}(a,b)&\leq m_{r,\nu}(|\phi|_2,|\phi|_2)\leq\max_{t\in\mathbb{R}}J_\nu(t\star(\phi,\phi))
\\&=\max_{t\in\mathbb{R}}\Bigl[2K(t\star\phi)-\nu e^{(\gamma_p+\gamma_q)t}\int_{\mathbb{R}^N}(I_\mu\ast|\phi|^{p})
|\phi|^{q}\Bigr].\endaligned
\end{equation*}
Note that $K(t\star\phi)\rightarrow0^+$ as $t\rightarrow-\infty$. Then for any $\epsilon>0$, there exists $t_0>0$ such that $K(t\star\phi)<\epsilon$ for any $t<-t_0$. From (\ref{3.6.1}) it follows that for some $\bar{\nu}>0$ there holds
\begin{equation*}\aligned
&\max_{t\geq -t_0}\bigl[2K(t\star\phi)
-\nu e^{(\gamma_p+\gamma_q)t}\int_{\mathbb{R}^N}(I_\mu\ast|\phi|^p)|\phi|^q\bigr]
<\epsilon,\quad\text{if}\ \nu\geq\bar{\nu}.\endaligned\end{equation*}
 Then $m_{r,\nu}(a,b)<\epsilon$ when $\nu\geq\bar{\nu}$. The conclusion follows from the arbitrariness of $\epsilon$ and Lemma \ref{l4.5.0}.\ \ \ \ \ $\Box$

In the super-critical case, due to the non-homogeneous and non-locality of the nonlinear term, it is
difficult for us to use the Riesz's rearrangement inequality (see \cite{LL}) to obtain $m_\nu(a,b)=m_{r,\nu}(a,b)$.
Consequently, different from \cite{Bartsch-Li-Zou} where only one threshold was considered, it is necessary to introduce two thresholds as follows
\begin{equation*}\aligned&\nu_1:=\inf\bigl\{\nu>0:
{ m_{\nu}(a,b)}<\frac{2^*_\mu-1}{22^*_\mu}
S^{\frac{2^*_\mu}{2^*_\mu-1}}_{H,L}\bigr\},\\ &\nu_2:=\inf\bigl\{\nu>0:
{ m_{r,\nu}(a,b)}<\frac{2^*_\mu-1}{22^*_\mu}
S^{\frac{2^*_\mu}{2^*_\mu-1}}_{H,L}\bigr\}.\endaligned\end{equation*}
Clearly, $0\leq\nu_1\leq\nu_2$. From Lemma \ref{l2.7} (ii) we know $\nu_1,\nu_2<+\infty$. To investigate $\nu_1,\nu_2$, we shall use the test function $\eta_\epsilon:=\phi U_\epsilon$, where
$U_\epsilon:=U_{\epsilon,0}$ is defined by (\ref{1.16}) and $\phi\in C^\infty_0(\mathbb{R}^N,[0,1])$ is a radial cut-off function satisfying
$\phi(x)=1$ if $|x|\leq\delta$ and $\phi(x)=0$ if $|x|\geq2\delta$ for some $\delta>0$.

In the following lemma, we estimate the coupling term $\int_{\mathbb{R}^N}(I_\mu\ast|\eta_\epsilon|^{p})|\eta_\epsilon|^{q}$.
\begin{lemma}\label{l4.10}
Assume (H$_0$) and (\ref{bu1}) are satisfied. Then for $\eta_\epsilon=\phi U_\epsilon$, as $\epsilon\rightarrow0^+$ we have
\begin{equation*}
\int_{\mathbb{R}^N}(I_\mu\ast|\eta_\epsilon|^{p})|\eta_\epsilon|^{q}\geq O(\epsilon^{2N-\mu-\frac{N-2}{2}(p+q)})-
O(\epsilon^{-\frac{N-2}{2}|p-q|+\frac{2N-\mu}2}).
\end{equation*}
\end{lemma}
{\bf Proof}:
Observe that
$$\aligned&\int_{\mathbb{R}^N}(I_\mu\ast|\eta_\epsilon|^{p})
|\eta_\epsilon|^{q}\\ \geq&\Bigl[\int_{B_{2\delta}}\int_{B_{2\delta}}-\int_{B_{2\delta}\setminus{B_\delta}}
\int_{B_{\delta}}-\int_{B_{\delta}}\int_{B_{2\delta}\setminus{B_\delta}}
-\int_{B_{2\delta}\setminus{B_\delta}}
\int_{B_{2\delta}\setminus{B_\delta}}\Bigr]
\frac{|U_\epsilon(x)|^{p}
|U_\epsilon(y)|^{q}}{|x-y|^\mu}dydx\\
:=&\mathbb{A}-\mathbb{B}-\mathbb{C}-\mathbb{D}.\endaligned$$
Let $C=(N(N-2))^{\frac{N-2}{4}(p+q)}$. For $\mathbb{A}$, we have
$$\aligned
\mathbb{A}
\geq&O(\epsilon^{2N-\mu-\frac{N-2}{2}(p+q)})
\int_{B_{{2\delta}}}\int_{B_{{2\delta}}}
\frac{1}{(1+{|x|^2})
^{\frac{N-2}{2}p}}\frac{1}{|x-y|^\mu}
\frac{1}{(1+{|y|^2})^{\frac{N-2}{2}q}}dydx\\=&O(\epsilon^{2N-\mu-\frac{N-2}{2}(p+q)}).
\endaligned$$
Note that $2_{\mu,\ast}=\frac{2N-\mu}{N}=\frac{2N-\mu}{2(N-2)}$ if $N=4$, which together with (\ref{bu1}), we know $p,q>\frac{2N-\mu}{2(N-2)}$ for both $N=3$ and $N=4$. Then
\begin{equation*}\int^{\frac{\delta}{\epsilon}}_0
\frac{r^{N-1}}{(1+r^2)
^{\frac{N-2}{2}q\frac{2N}{2N-\mu}}}dr<+\infty.\end{equation*}
Thus the Hardy-Littlewood-Sobolev inequality implies
$$\aligned
\mathbb{B}
\leq&C\epsilon^{2N-\mu-\frac{N-2}{2}(p+q)}
\Bigl(\int_{\frac{2\delta}{\epsilon}\geq|x|\geq\frac{\delta}{\epsilon}}
\frac{1}{|x|
^{({N-2})p\frac{2N}{2N-\mu}}}dx\Bigr)^{\frac{2N-\mu}{2N}}
\Bigl(\int_{|y|\leq\frac{\delta}{\epsilon}}
\frac{1}{(1+{|y|^2})
^{\frac{N(N-2)q}{2N-\mu}}}dy\Bigr)^{\frac{2N-\mu}{2N}}
\\ \leq &C\epsilon^{2N-\mu-\frac{N-2}{2}(p+q)}
\Bigl(|x|^{-(N-2)p\frac{2N}{2N-\mu}+N}\Big|^{\frac{2\delta}{\epsilon}}
_{\frac{\delta}{\epsilon}}\Bigr)^{\frac{2N-\mu}{2N}}
\Bigl(\int^{\frac{\delta}{\epsilon}}_0\frac{r^{N-1}}{(1+r^2)
^{\frac{N-2}{2}q\frac{2N}{2N-\mu}}}dr\Bigr)^{\frac{2N-\mu}{2N}}\\
=&O(\epsilon^{\frac{2N-\mu}{2}+\frac{N-2}{2}(p-q)}).
\endaligned$$
Similarly, by $p>\frac{2N-\mu}{2(N-2)}$ we get
$\mathbb{C}\leq O(\epsilon^{\frac{N-2}{2}(q-p)+\frac{2N-\mu}{2}}).$
For $\mathbb{D}$, we infer
$$\aligned
\mathbb{D}
\leq&C\epsilon^{\frac{N-2}{2}(p+q)}\int_{B_{2\delta}\backslash{B_\delta}}
\int_{B_{2\delta}\backslash{B_\delta}}
\frac{1}{|x|
^{({N-2})p}}\frac{1}{|x-y|^\mu}
\frac{1}{|y|^{(N-2)q}}dydx=O(\epsilon^{\frac{N-2}{2}(p+q)}).
\endaligned$$
Then the conclusion holds true.
\ \ \ \ $\Box$

\begin{lemma}\label{l2.9} Assume (H$_0$) and (\ref{bu1}) are satisfied. Then $m_{\nu}(a,b)\leq m_{r,\nu}(a,b)\leq \frac{2^*_\mu-1}{22^*_\mu}
S^{\frac{2^*_\mu}{2^*_\mu-1}}_{H,L}$ for all $\nu>0$. If additionally \begin{equation*}\aligned
\left\{ \begin{array}{lll}
\min\{p,q\}<2\ \text{and} \ 2<\mu<3, \ & \text{if}\quad N=3,\\
\min\{p,q\}<2, \quad\quad\quad\quad \quad\quad   \ & \text{if}\quad
N=4,
\end{array}\right.\endaligned
\end{equation*} then $\nu_2=\nu_1=0$.
\end{lemma}
{\bf Proof}: We only consider $q<p$, since the case $q>p$ can be treated similarly. For the test function $\eta_\epsilon$ as in Lemma \ref{l4.10}, the following estimates hold, which can be found in \cite{Jeanjeanle,YangSCICHINA,SoaveJFA,WM}.
\begin{equation}\label{2.24}|\nabla \eta_\epsilon|^2_2=S^{\frac{N}{2}}+O(\epsilon^{N-2})
=C(N,\mu)^{\frac{N}{22^*_\mu}}
S^{\frac{N}{2}}_{H,L}+O(\epsilon^{N-2}),\end{equation}
where $C(N,\mu)$ is given in Proposition \ref{p1.1}. Moreover,
\begin{equation}\label{2.25}
|\eta_\epsilon|^2_2=
O(\epsilon^2|\ln\epsilon|)\  \text{if}\ N=4,\quad
|\eta_\epsilon|^2_2=O(\epsilon)  \  \text{if}\ N=3,
\end{equation}
\begin{equation}\label{2.26}\int_{\mathbb{R}^N}
(I_\mu\ast|\eta_\epsilon|^{2^*_\mu})|\eta_\epsilon|^{2^*_{\mu}}\geq
 C(N,\mu)^{\frac N{2}}S^{\frac{2N-\mu}{2}}_{H,L}
 -{O}(\epsilon^{\frac{2N-\mu}{2}}),\end{equation}
 \begin{equation}\label{2.27.0}
|\eta_\epsilon|^p_p=\left\{ \begin{array}{lll}
O(\epsilon^{N-\frac{N-2}{2}p}),\ & \text{if}\ N\geq4\ \text{and}\ p\in(2,2^*)\ \text{or}\ \text{if}\ N=3\ \text{and}\ p\in(3,6),\\
O(\epsilon^{\frac p2}),\ & \text{if}\ N=3\ \text{and}\ p\in(2,3),\\
O(\epsilon^{\frac{3}{2}}|\ln\epsilon|),\ & \text{if}\ N=3\ \text{and}\ p=3.
\end{array}\right.
\end{equation}
 By the Hardy-Littlewood-Sobolev inequality and Lemma \ref{l1.1} we get
\begin{equation}\label{2.28}\aligned
&\Bigl(\int_{\mathbb{R}^N}
(I_\mu\ast|\eta_\epsilon|^{2^*_\mu})
|\eta_\epsilon|^{2^*_{\mu}}\Bigr)^{\frac1{2^*_\mu}}
\leq C(N,\mu)^{\frac1{2^*_{\mu}}}|\eta_\epsilon|^2_{2^*}
=C(N,\mu)^{\frac{N}{22^*_\mu}}
S^{\frac{N}{2}}_{H,L}+O(\epsilon^{N-2}).\endaligned
\end{equation}

Set $u_\epsilon=\frac{a}{|\eta_\epsilon|_2}\eta_\epsilon$,
$v_\epsilon=\frac{\epsilon^\theta}{|\eta_\epsilon|_2}\eta_\epsilon$ with constant $\theta>0$ to be determined. For $t_\epsilon=t_\nu(u_\epsilon,v_\epsilon)$ given by Lemma \ref{l2.4}, $t_\epsilon\star(u_\epsilon,v_\epsilon)\in \mathcal{P}_\nu(a,\epsilon^\theta)\cap H_{rad}$.
 By the definition ${P}_\nu(t_\epsilon\star(u_\epsilon,v_\epsilon))=0$ we have $e^{t_\epsilon}\leq \mathbb{A}^{\frac{1}{22^*_\mu-2}},$
where $$\mathbb{A}:=\frac{|\nabla u_\epsilon|^2_2+|\nabla v_\epsilon|^2_2}{\int_{\mathbb{R}^N}
[(I_\mu\ast|u_\epsilon|^{2^*_\mu})|u_\epsilon|^{2^*_\mu}+
(I_\mu\ast|v_\epsilon|^{2^*_\mu})|v_\epsilon|^{2^*_\mu}]}=\frac{a^2+\epsilon^{2\theta}}{a^{22^*_\mu}+
\epsilon^{22^*_\mu\theta}}\frac{|\nabla \eta_\epsilon|^2_2|\eta_\epsilon|^{22^*_\mu-2}_2}{\int_{\mathbb{R}^N}
(I_\mu\ast|\eta_\epsilon|^{2^*_\mu})|\eta_\epsilon|^{2^*_\mu}}.$$

Also we set$$\mathbb{D}:= \frac{\int_{\mathbb{R}^N}(I_\mu\ast|u_\epsilon|^p)|v_\epsilon|^q}
{\int_{\mathbb{R}^N}[(I_\mu\ast|u_\epsilon|^{2^*_\mu})|u_\epsilon|^{2^*_\mu}+
(I_\mu\ast|v_\epsilon|^{2^*_\mu})|v_\epsilon|^{2^*_\mu}]}=\frac{\epsilon^{\theta q}|\eta_\epsilon|^{22^*_\mu-p-q}_2}{a^{22^*_\mu}+\epsilon^{22^*_\mu\theta}}\frac{\int_{\mathbb{R}^N}
(I_\mu\ast|\eta_\epsilon|^{p})|\eta_\epsilon|^{q}
}{\int_{\mathbb{R}^N}
(I_\mu\ast|\eta_\epsilon|^{2^*_\mu})|\eta_\epsilon|^{2^*_\mu}}.$$
In view of ${P}_\nu(t_\epsilon\star(u_\epsilon,v_\epsilon))=0$, $\gamma_p+\gamma_q>2$, and $e^{t_\epsilon}\leq \mathbb{A}^{\frac{1}{22^*_\mu-2}}$, for $\epsilon>0$ small there holds
\begin{equation*}\aligned
&e^{(22^*_\mu-2)t_\epsilon}=\mathbb{A}-(\gamma_p+\gamma_q)\nu e^{(\gamma_p+\gamma_q-2)t_\epsilon}\mathbb{D}
\geq\mathbb{A}\bigl[1-(\gamma_p+\gamma_q)\nu\mathbb{A}^{
\frac{\gamma_p+\gamma_q-22^*_\mu}{22^*_\mu-2}}\mathbb{D}\bigr]\\
\geq& \frac{|\nabla \eta_\epsilon|^2_2|\eta_\epsilon|^{22^*_\mu-2}_2}{\int_{\mathbb{R}^N}
(I_\mu\ast|\eta_\epsilon|^{2^*_\mu})|\eta_\epsilon|^{2^*_\mu}}
\biggl[C_1-C_2\nu\epsilon^{\theta q}\frac{\int_{\mathbb{R}^N}
(I_\mu\ast|\eta_\epsilon|^{p})|\eta_\epsilon|^{q}}
{|\eta_\epsilon|^{p+q-\gamma_p-\gamma_q}_2}\frac{|\nabla \eta_\epsilon|^{\frac{\gamma_p+\gamma_q-22^*_\mu}{2^*_\mu-1}}_2}{
\Bigl(\int_{\mathbb{R}^N}
(I_\mu\ast|\eta_\epsilon|^{2^*_\mu})|\eta_\epsilon|^{2^*_\mu}\Bigr)
^{\frac{\gamma_p+\gamma_q-2}{22^*_\mu-2}}}\biggr]
.\endaligned\end{equation*}
By (\ref{2.24}), (\ref{2.26}) and (\ref{2.28}),  it suffices to estimate
$\frac{\int_{\mathbb{R}^N}
(I_\mu\ast|\eta_\epsilon|^{p})|\eta_\epsilon|^{q}}
{|\eta_\epsilon|^{p+q-\gamma_p-\gamma_q}_2}.$
If $N=4$, note that $\frac{2N}{2N-\mu}p,\frac{2N}{2N-\mu}q\in(2,2^*)$.
If $N=3$, $p,q>\frac{6-\mu}{2}$ implies that $\frac{2N}{2N-\mu}p,\frac{2N}{2N-\mu}q\in(3,6)$.
From (\ref{2.27.0}), (\ref{2.25}) and the Hardy-Littlewood-Sobolev inequality it follows that
$$\aligned\frac{\int_{\mathbb{R}^N}
(I_\mu\ast|\eta_\epsilon|^{p})|\eta_\epsilon|^{q}}
{|\eta_\epsilon|^{p+q-\gamma_p-\gamma_q}_2}&\leq C_1\frac{|\eta_\epsilon|^p_{\frac{2N}{2N-\mu}p}|\eta_\epsilon|^q_{\frac{2N}{2N-\mu}q}}{|\eta_\epsilon|
^{p+q-\gamma_p-\gamma_q}_2}\leq\left\{ \begin{array}{lll}
O(\epsilon^{\frac{p+q-\gamma_p-\gamma_q}{2}}), \quad   \ &\text{if}\quad N=3,\\
\frac{C_1}{|\ln \epsilon|^{\frac{p+q-\gamma_p-\gamma_q}{2}}}=o_\epsilon(1),\quad \ &\text{if}\quad N=4,
\end{array}\right.\endaligned$$
where we have used the fact $2N-\mu-\frac{N-2}{2}(p+q)=p+q-\gamma_p-\gamma_q$.
Thus
$$e^{(22^*_\mu-2)t_\epsilon}\geq C\frac{|\nabla \eta_\epsilon|^2_2|\eta_\epsilon|^{22^*_\mu-2}_2}{\int_{\mathbb{R}^N}
(I_\mu\ast|\eta_\epsilon|^{2^*_\mu})|\eta_\epsilon|^{2^*_\mu}}.$$
Then
\begin{equation*} e^{t_\epsilon}\geq C|\eta_\epsilon|_2,\ \text{as\ } \epsilon\rightarrow0 \text{\ for some constant\ } C>0.\end{equation*} and so for $\epsilon>0$ small we obtain
\begin{equation}\label{2.30.0}\aligned
&m_{r,\nu}(a,b)\leq m_{r,\nu}(a,\epsilon^\theta)\leq J_\nu(t_\epsilon\star(u_\epsilon,v_\epsilon))\\
\leq&\max_{s>0}\Bigl(\frac{a^2+\epsilon^{2\theta}}{2}s^2
\frac{|\nabla \eta_\epsilon|^2_2}{|\eta_\epsilon|^2_2}-\frac{a^{22^*_\mu}+
\epsilon^{22^*_\mu\theta}}{22^*_\mu}s^{22^*_\mu}
\frac{\int_{\mathbb{R}^N}(I_\mu\ast|\eta_\epsilon|^{2^*_\mu})
|\eta_\epsilon|^{2^*_\mu}}{|\eta_\epsilon|^{22^*_\mu}_2}\Bigr)\\&-\nu \epsilon^{\theta q}e^{(\gamma_p+\gamma_q)t_\epsilon}
\frac{\int_{\mathbb{R}^N}(I_\mu\ast|\eta_\epsilon|^p)|\eta_\epsilon|^q}{
|\eta_\epsilon|^{p+q}_2}
\\=&\frac{2^*_\mu-1}{22^*_\mu}
\Bigl[\frac{(a^2+\epsilon^{2\theta})^{2^*_\mu}}{(a^{22^*_\mu}+
\epsilon^{22^*_\mu\theta})}\frac{|\nabla \eta_\epsilon|^{22^*_\mu}_2}
{\int_{\mathbb{R}^N}
(I_\mu*|\eta_\epsilon|^{2^*_\mu})|\eta_\epsilon|^{2^*_\mu}}\Bigr]^{\frac{1}{2^*_\mu-1}}-\nu \epsilon^{\theta q}e^{(\gamma_p+\gamma_q)t_\epsilon}
\frac{\int_{\mathbb{R}^N}(I_\mu*|\eta_\epsilon|^p)|\eta_\epsilon|^q}{
|\eta_\epsilon|^{p+q}_2}\\
\leq&\frac{2^*_\mu-1}{22^*_\mu}
S^{\frac{2^*_\mu}{2^*_\mu-1}}_{H,L}+O(\epsilon^{N-2})+O(\epsilon^{2\theta})-
\nu C
\epsilon^{\theta q}\frac{\int_{\mathbb{R}^N}(I_\mu\ast|\eta_\epsilon|^p)|\eta_\epsilon|^q}{
|\eta_\epsilon|^{p+q-\gamma_p-\gamma_q}_2}.\endaligned\end{equation}
Note that  $q<p$, by Lemma \ref{l4.10} and (\ref{2.25}) we get
\begin{equation*}\frac{\int_{\mathbb{R}^N}(I_\mu\ast|\eta_\epsilon|^p)
|\eta_\epsilon|^q}{
|\eta_\epsilon|^{p+q-\gamma_p-\gamma_q}_2}\geq\left\{ \begin{array}{lll}
O(\epsilon^{\frac{p+q-\gamma_p-\gamma_q}{2}})
-O(\epsilon^{\frac{2N-\mu}{2}-\frac{N-2}{2}(p-q)-\frac{p+q-\gamma_p-\gamma_q}{2}}),\ & N=3,\\
C|\ln\epsilon|^{\frac{\gamma_p+\gamma_q-p-q}{2}}(1-O(\epsilon^{\frac{2N-\mu}{2}-\frac{N-2}{2}(p-q)-p-q+\gamma_p+\gamma_q}),\ & N=4.
\end{array}\right.\end{equation*}
Then using (\ref{2.30.0}) there holds
\begin{equation}\label{2.34.0}\aligned m_{r,\nu}(a,b)\leq\frac{2^*_\mu-1}{22^*_\mu}
S^{\frac{2^*_\mu}{2^*_\mu-1}}_{H,L}+\left\{ \begin{array}{lll}
O(\epsilon)+O(\epsilon^{2\theta})-
\epsilon^{\theta q}O(\epsilon^{\frac{p+q-\gamma_p-\gamma_q}{2}}),\quad& N=3,\\
O(\epsilon^2)+O(\epsilon^{2\theta})-
C\epsilon^{\theta q}|\ln\epsilon|^{\frac{\gamma_p+\gamma_q-p-q}{2}},\quad& N=4.
\end{array}\right.\endaligned\end{equation}
 Letting $\epsilon\rightarrow0^+$, we know $m_\nu(a,b)\leq m_{r,\nu}(a,b)\leq\frac{2^*_\mu-1}{22^*_\mu}
S^{\frac{2^*_\mu}{2^*_\mu-1}}_{H,L}$. If $q<p$ and $N=3$, we can choose $\theta>0$ in the above argument such that $\theta q+\frac{p+q-\gamma_p-\gamma_q}{2}<\min\{1,2\theta\}$ since $\frac{6-\mu}{2}<\min\{p,q\}<2$ with $2<\mu<3$. It follows from (\ref{2.34.0}) that $m_{r,\nu}(a,b)<\frac{2^*_\mu-1}{22^*_\mu}
S^{\frac{2^*_\mu}{2^*_\mu-1}}_{H,L}$ for $\epsilon$ small enough.
If $q<p$ and $N=4$, we choose $\theta=1$ in the above argument and from (\ref{2.34.0}) we get $m_{r,\nu}(a,b)<\frac{2^*_\mu-1}{22^*_\mu}
S^{\frac{2^*_\mu}{2^*_\mu-1}}_{H,L}$ using $\min\{p,q\}<2$. In a word, $\nu_2=0$ and so $\nu_2=\nu_1=0$. \ \ \ $\Box$

\begin{lemma} \label{l2.9.0}Let constants $A,B>0$ and $\theta\geq1$ be fixed. For any $\nu>0$, let
$$\Theta_\nu=\{\sigma\in\mathbb{R}^+|A\sigma^2\leq S^{-2^*_\mu}_{H,L}\sigma^{22^*_\mu}+B\nu\sigma^{2\theta}\}.$$
Then for all $\epsilon>0$, there exists $\tilde{\nu}>0$ depending only on $\epsilon, A,B$ and $\theta$ such that
$$\inf \Theta_\nu\geq (1-\epsilon)S^{\frac{2^*_\mu}{2(2^*_\mu-1)}}_{H,L}A^{\frac1{2(2^*_\mu-1)}},\quad\text{for all}\ 0<\nu<\tilde{\nu}.$$
\end{lemma}
{\bf Proof}: If $\theta=1$, the conclusion easily yields. Below we assume $\theta>1$.  Consider a function $f(\sigma)=A\sigma^{2-2\theta}-S^{-2^*_\mu}_{H,L}\sigma^{22^*_\mu-2\theta}.$
Then $\sigma_0:=f^{-1}(0)=(S^{2^*_\mu}_{H,L}A)^{\frac{1}{22^*_\mu-2}}$.
 If $1<\theta\leq2^*_\mu$, then $f$ is strictly decreasing. If $\theta>2^*_\mu$, then $f$ has a global negative minimum at
$\sigma=\Bigl(\frac{2-2\theta}{22^*_\mu-2\theta}S^{2^*_\mu}_{H,L}A\Bigr)
^{\frac1{22^*_\mu-2}}>\sigma_0$, and $f$ tends to $0$ as $\sigma\rightarrow+\infty$. In any case, $f$ is decreasing in $(0, \sigma_0)$, has only one zero point  at $\sigma_0$, $\lim_{\sigma\rightarrow0^+}f(\sigma)=+\infty$, and $f(\sigma)<0$ in $(\sigma_0,+\infty)$. Hence, $\inf\Theta_\nu=f^{-1}(B\nu)\rightarrow \sigma_0$ as $\nu\rightarrow0^+$.\ \ \ \ $\Box$

\begin{lemma}\label{l2.9.1} Let (H$_0$) hold.
If $\min\{p,q\}\geq2+\frac{4-\mu}{N}$, then $\nu_2>0$.
\end{lemma}
{\bf Proof}: We argue by contradiction by assuming there exists a sequence $\nu_n\rightarrow 0^+$ such that
$0<m_{r,\nu_n}(a,b)<\frac{2^*_\mu-1}{22^*_\mu}
S^{\frac{2^*_\mu}{2^*_\mu-1}}_{H,L}$ {for all } $n.$
In view of Lemmas \ref{l2.6} and \ref{l2.7}, we can apply Lemma \ref{l3.4} with $c=m_{r,\nu_n}(a,b)$ to obtain a minimizer $(u_n,v_n)\in H_{rad}$ of $m_{r,\nu_n}(a,b)$ for any $n$. Then $(u_n,v_n)$ solves (\ref{0.0}), i.e.
\begin{equation}\label{2.35}
\aligned
\left\{ \begin{array}{lll}
-\Delta u_n+\lambda_{1,n}u_n=(I_\mu\ast |u_n|^{2^*_\mu})|u_n|^{2^*_\mu-2}u_n+\nu_n p(I_\mu\ast |v_n|^q)|u_n|^{p-2}u_n\ & \text{in}\quad \mathbb{R}^N,\\
-\Delta v_n+\lambda_{2,n}v_n=(I_\mu\ast |v_n|^{2^*_\mu})|v_n|^{2^*_\mu-2}v_n+\nu_n q(I_\mu\ast |u_n|^p)|v_n|^{q-2}v_n\ & \text{in}\quad \mathbb{R}^N,\\
\int_{\mathbb{R}^N}u^2_n=a^2,\quad\int_{\mathbb{R}^N}v^2_n=b^2,
\end{array}\right.\endaligned
\end{equation}
with $\lambda_{1,n},\lambda_{2,n}>0$. Taking similar arguments as in Lemma \ref{l3.4.0}, we get $(u_n,v_n)$ is bounded in $H$. The Hardy-Littlewood-Sobolev inequality implies $\Bigl|\int_{\mathbb{R}^N}(I_\mu\ast|u_n|^p)|v_n|^q\Bigr|\leq C$.
Then
\begin{equation}\label{2.36}\aligned&\frac{2^*_\mu-1}{22^*_\mu}
S^{\frac{2^*_\mu}{2^*_\mu-1}}_{H,L}>m_{r,\nu_n}(a,b)
=J_{\nu_n}(u_n,v_n)-
\frac{1}{22^*_\mu}P_{\nu_n}(u_n,v_n)\\
\geq&\bigl(\frac12-\frac1{22^*_\mu}\bigr)(|\nabla u_n|^2_2+|\nabla v_n|^2_2)-C\nu_n=\bigl(\frac12-\frac1{22^*_\mu}\bigr)(|\nabla u_n|^2_2+|\nabla v_n|^2_2)+o_n(1).
\endaligned\end{equation}
Using the definition of $S_{H,L}$, the inequality (3.3) in \cite{Gvans} and (\ref{1.14}), testing the first equation of  (\ref{2.35}) with $(u_n,0)$ we get
$$\aligned|\nabla u_n|^2_2&\leq S^{-2^*_\mu}_{H,L}|\nabla u_n|^{22^*_\mu}_2+\nu_np\Bigl(\int_{\mathbb{R}^N}
(I_\mu\ast|u_n|^p)|u_n|^p \Bigr)^{\frac12}\Bigl(\int_{\mathbb{R}^N}(I_\mu\ast|v_n|^q)|v_n|^q \Bigr)^{\frac12}\\
&\leq S^{-2^*_\mu}_{H,L}|\nabla u_n|^{22^*_\mu}_2+\nu_nC_1|\nabla u_n|^{\frac{N(p-2)+\mu}{2}}_2,  \endaligned$$
for some $C_1>0$ independent of $n$, where we have used $\{v_n\}$ is bounded in $H^1(\mathbb{R}^N)$.
Since $\min\{p,q\}\geq2+\frac{4-\mu}{N}$ we get $\frac{N(p-2)+\mu}{2}\geq2$. Applying Lemma \ref{l2.9.0} with $\theta=\frac{N(p-2)+\mu}{4}$, we infer
$\liminf_{n\rightarrow\infty}|\nabla u_n|^2_2\geq S^{\frac{2^*_\mu}{2^*_\mu-1}}_{H,L}$. Similarly, $\liminf_{n\rightarrow\infty}|\nabla v_n|^2_2\geq S^{\frac{2^*_\mu}{2^*_\mu-1}}_{H,L},$
contradicting (\ref{2.36}).\ \ \ \ \ \ \ $\Box$

Letting $e^t=s$ and  $s\circ(u,v)=(s\circ u,s\circ v)$ with $s\circ u=s^\frac{N}{2}u(sx)$, from (\ref{4.11.0}) we know \begin{equation}\label{4.3}m_\nu(a,b)=\inf_{\mathcal{S}(a,b)}
\max_{s>0}J_\nu(s\circ(u,v)).\end{equation}
For any $(u,v)\in \mathcal{S}(a,b)$, let $e^{t_\nu(u,v)}=s_\nu(u,v)$, where $t_\nu(u,v)$ is  defined in Lemma \ref{l2.4}, then $s_\nu(u,v)\circ(u,v)\in \mathcal{P}_\nu(a,b)$, $\max_{s>0}J_\nu(s\circ(u,v))
=J_\nu(s_\nu(u,v)\circ(u,v))$ and $f(\nu,(u,v),s_\nu(u,v))=0$, where
\begin{equation*}\aligned
f(\nu,(u,v),s)=&\nu(\gamma_p+\gamma_q)s^{\gamma_p+\gamma_q-2}
\int_{\mathbb{R}^N}(I_\mu\ast|u|^{p})|v|^{q}+
s^{22^*_\mu-2}\int_{\mathbb{R}^N}\bigl[(I_\mu\ast|u|^{2^*_\mu})|u|^{2^*_\mu}
\\&+
(I_\mu\ast|v|^{2^*_\mu})|v|^{2^*_\mu}\bigr]-(|\nabla u|^2_2+|\nabla v|^2_2), \quad s>0.
\endaligned\end{equation*}
Clearly, $f(\nu,(u,v),s)$ is increasing in $s>0$.

\begin{lemma}\label{l2.10.0} If $0<\nu<\nu_1$, then $m_\nu(a,b)=\frac{2^*_\mu-1}{22^*_\mu}
S^{\frac{2^*_\mu-1}{2^*_\mu}}_{H,L}$ is not achieved.\end{lemma}
{\bf Proof}: Clearly, if $0<\nu<\nu_1$, then $m_\nu(a,b)=\frac{2^*_\mu-1}{22^*_\mu}
S^{\frac{2^*_\mu-1}{2^*_\mu}}_{H,L}$ and it suffices to show $m_\nu(a,b)$ is not achieved.  Denote $a(\nu)=m_\nu(a,b)$. Argue  by contradiction we assume $a(\nu)$ is achieved by some $(u_\nu,v_\nu)\in H$.
By Lemma \ref{l2.7} we know $a'(\nu)=\frac{d{m_\nu(a,b)}}{d\nu}\leq0$ exists for almost every $\nu\in (0,+\infty)$.  Let $\nu\in(0,+\infty)$ be such that $a'(\nu)$ exists.  The following claim describes the relation between $a'(\nu)$ and $(u_\nu,v_\nu)$, which is inspired by
\cite{Chen-Zou}. We claim that
\begin{equation} \label{2.0.8} a'(\nu)=-\int_{\mathbb{R}^N}(I_\mu\ast|u_\nu|^p)|v_\nu|^q.\end{equation}
Indeed, for simplicity, denote
$$\aligned A&=|\nabla u_\nu|^2_2+|\nabla v_\nu|^2_2,\quad B=\int_{\mathbb{R}^N}(I_\mu\ast |u_\nu|^p)|v_\nu|^q,\\
C&=\int_{\mathbb{R}^N}[(I_\mu\ast |u_\nu|^{2^*_\mu})|u_\nu|^{2^*_\mu}+(I_\mu\ast |v_\nu|^{2^*_\mu})|v_\nu|^{2^*_\mu}],
\endaligned$$
 and consider
$$F(\mu,s):=f(\mu,(u_\nu,v_\nu),s)=\mu(\gamma_p+\gamma_q)
s^{\gamma_p+\gamma_q-2}B+s^{22^*_\mu-2}C-A, \ s>0.$$
Note that $F(\nu,1)=0$, $\frac{\partial}{\partial s}F(\nu,1)=\nu(\gamma_p+\gamma_q)(\gamma_p+\gamma_q-2)B+(22^*_\mu-2)C>0$ since $\gamma_p+\gamma_q>2$,
and $F(\mu,s_\mu(u_\nu,v_\nu))\equiv0$. By the implicit function theorem, there is a small neighborhood $(\nu-\delta,\nu+\delta)$ of $\nu$ in $(0,+\infty)$ such that $s(\mu):=s_\mu(u_\nu,v_\nu)\in C^\infty((\nu-\delta,\nu+\delta),\mathbb{R})$ as a function of $\mu$.
By $F(\mu,s(\mu))\equiv0$ we have
\begin{equation}\label{2.4.0}
s'(\nu)=-\frac{\frac{\partial F}{\partial\mu}}{\frac{\partial F}{\partial s}}\Bigl|_{\mu=\nu}=-\frac{(\gamma_p+\gamma_q)B}{
\nu(\gamma_p+\gamma_q)(\gamma_p+\gamma_q-2)B+(22^*_\mu-2)C},
\end{equation}
where we have used the fact that  $s_\nu(u_\nu,v_\nu)\equiv1$ since $(u_\nu,v_\nu)$ is a critical point of $J_\nu$.
The  Taylor expansion implies $s(\mu)=1+s'(\nu)(\mu-\nu)+O((\mu-\nu)^2)$
and so
$$s^2(\mu)=1+2s'(\nu)(\mu-\nu)+O((\mu-\nu)^2),\ s^{22^*_\mu}(\mu)=1+22^*_\mu s'(\nu)(\mu-\nu)+O((\mu-\nu)^2).$$
Taking into account (\ref{2.4.0}) and the fact that $A=\nu(\gamma_p+\gamma_q)B+C$ we deduce
$$\aligned
a(\mu)&\leq\max_{s>0}J_\mu(s\circ(u_\nu,v_\nu))=\bigl(\frac12-\frac{1}{\gamma_p+\gamma_q}\bigr)
s^2(\mu)A+\bigl(\frac{1}{\gamma_p+\gamma_q}-\frac{1}{22^*_\mu}\bigr)
s^{22^*_\mu}(\mu)C\\&=a(\nu)-B(\mu-\nu)+O((\mu-\nu)^2).\endaligned$$
Letting $\mu\nearrow\nu$ and  $\mu\searrow\nu$ respectively, we obtain
$a'(\nu)=-B$, that is, (\ref{2.0.8}) holds true.
So $a(\nu)$ is not achieved for $\nu\in(0,\nu_1)$. \ \ \ \ $\Box$

\begin{Remark}Taking similar arguments of Lemma \ref{l2.10.0}, we can infer that $m_{r,\nu}(a,b)=\frac{2^*_\mu-1}{22^*_\mu}
S^{\frac{2^*_\mu-1}{2^*_\mu}}_{H,L}$ is not achieved when $0<\nu<\nu_2$.\end{Remark}
{\bf Proof of Theorems \ref{th1.3} and \ref{th1.4}} Applying Lemma \ref{l3.4} with $c=m_{r,\nu}(a,b)$, and combining Lemmas \ref{l4.5.0}, \ref{l2.6}, \ref{l2.7}, \ref{l2.9}, \ref{l2.9.1} and \ref{l2.10.0}, the results follow directly. \ \ \ \ $\Box$
\subsection{$L^2$-critical case}
Define
$$\Omega_\nu(a,b):=\Bigl\{(u,v)\in \mathcal{S}(a,b):|\nabla u|^2_2+|\nabla v|^2_2>2\nu\int_{\mathbb{R}^N}(I_\mu\ast|u|^p)|v|^q\Bigr\}.$$
Using (\ref{1.14}), it is easy to see that $\Omega_\nu(a,b)=\mathcal{S}(a,b)\neq\emptyset$ when $0<\nu<\nu'_0$ with
\begin{equation}\label{niu0}\nu'_0:=\frac12
(a^2+b^2)^{-\frac{p+q-2}{2}}C^{-1}_{N,p,q},\end{equation}
where $C_{N,p,q}$ is given in (\ref{1.14}).
Then replacing the condition $\nu>0$ by $0<\nu<\nu'_0$, Lemma \ref{l2.4} still hold, and as (\ref{4.11.0}) the minimum $m_\nu(a,b)$ and $m_{r,\nu}(a,b)$ with $0<\nu<\nu'_0$ have minimax characterizations:
\begin{equation}\label{4.11.1}m_\nu(a,b)=\inf_{\Omega_\nu(a,b)}
\max_{t\in\mathbb{R}}J_\nu(t\star(u,v)),\quad m_{r,\nu}(a,b)=\inf_{\Omega_\nu(a,b)\cap H_{rad}}
\max_{t\in\mathbb{R}}J_\nu(t\star(u,v)).\end{equation}
 \begin{lemma}\label{l2.14}
\noindent(i) $m_{\nu}(a,b)=m_{r,\nu}(a,b)$ if $0<\nu<\nu'_0$.\\
 \noindent(ii) $m_\nu(a,b)$ satisfies (\ref{2.1}) if $0<\nu<\nu'_0$ and $m_\nu(a,b)$ is decreasing in $\nu\in(0,\nu'_0)$.
 \end{lemma}
{\bf Proof}: (i) Argue by contradiction, we assume that there is $(u,v)\in \mathcal{P}_\nu(a,b)\backslash \mathcal{S}_r(a,b)$ such that
\begin{equation}\label{b1}J_\nu(u,v)<m_{r,\nu}(a,b).\end{equation}
Let $(u^*,v^*)$ be the symmetric decreasing rearrangement of $(u,v)$. From the properties of symmetric decreasing rearrangement and \cite[Theorem 3.7]{LL} we get
\begin{equation}\label{b3}\aligned&|\nabla u^*|^2_2\leq|\nabla u|^2_2,\ |\nabla v^*|^2_2\leq|\nabla v|^2_2, \ \int_{\mathbb{R}^N}(I_\mu\ast|u^*|^{2^*_\mu})|u^*|^{2^*_\mu}\geq
\int_{\mathbb{R}^N}(I_\mu\ast|u|^{2^*_\mu})|u|^{2^*_\mu},\\ & \int_{\mathbb{R}^N}(I_\mu\ast|u^*|^{p})|v^*|^{q}\geq
\int_{\mathbb{R}^N}(I_\mu\ast|u|^{p})|v|^{q},  \ \int_{\mathbb{R}^N}(I_\mu\ast|v^*|^{2^*_\mu})|v^*|^{2^*_\mu}\geq
\int_{\mathbb{R}^N}(I_\mu\ast|v|^{2^*_\mu})|v|^{2^*_\mu}
.\endaligned\end{equation}
Therefore
\begin{equation}\label{b2} \ J_\nu(u^*,v^*)\leq J_\nu(u,v), \ P_\nu(u^*,v^*)\leq P_\nu(u,v)=0.\end{equation}
If $P_\nu(u^*,v^*)=0$, then $(u^*,v^*)\in \mathcal{P}_\nu(a,b)\cap H_{rad}$ and by (\ref{b1}) we get $J_\nu(u,v)<J_\nu(u^*,v^*)$. This contradicts (\ref{b2}). If $P_\nu(u^*,v^*)<0$, then Lemma \ref{l2.4} with $0<\nu<\nu'_0$ implies $t_\nu:=t_\nu(u^*,v^*)<0$. Since $\gamma_p+\gamma_q=2$, by (\ref{b3}) we obtain
$$\aligned &m_{r,\nu}(a,b)\leq J_\nu(t_\nu(u^*,v^*)\star(u^*,v^*))\\=&
\frac{2^*_\mu-1}{22^*_\mu}e^{2t_\nu}\bigl[|\nabla u^*|^2_2+|\nabla v^*|^2_2-2\nu\int_{\mathbb{R}^N}(I_\mu\ast|u^*|^{p})|v^*|^{q}\bigr]\\
\leq&\frac{2^*_\mu-1}{22^*_\mu}\bigl[|\nabla u|^2_2+|\nabla v|^2_2-2\nu\int_{\mathbb{R}^N}(I_\mu\ast|u|^{p})|v|^{q}\bigr]
=J_\nu(u,v),
\endaligned$$
which contradicts the hypothesis (\ref{b1}). Therefore, $m_\nu(a,b)=m_{r,\nu}(a,b)$.

(ii) By (\ref{4.11.1}) and the arguments as in the proof of Lemma \ref{l2.7} (i), the conclusion (ii) follows.\ \ \ \ $\Box$

\begin{lemma}\label{l2.15} Assume
 (\ref{1.11.0}) is satisfied. Then
$m_{\nu}(a,b)\in(0,\frac{2^*_\mu-1}{22^*_\mu}
S^{\frac{2^*_\mu}{2^*_\mu-1}}_{H,L})$ if $0<\nu<\nu'_0$.
\end{lemma}
{\bf Proof}: For $(u,v)\in \mathcal{P}_\nu(a,b)$, the definition of $S_{H,L}$, (\ref{1.14}) and $\gamma_p+\gamma_q=2$ imply
\begin{equation*}\aligned \frac{\nu_0-\nu}{\nu_0} (|\nabla u|^2_2+|\nabla v|^2_2)\leq\int_{\mathbb{R}^N}[(I_\mu\ast|u|^{2^*_\mu})|u|^{2^*_\mu}+
(I_\mu\ast|v|^{2^*_\mu})|v|^{2^*_\mu}]\leq C(\nabla u|^2_2+|\nabla v|^2_2)^{2^*_\mu}.
\endaligned\end{equation*}
Then $$\inf_{\mathcal{P}_\nu(a,b)}\int_{\mathbb{R}^N}[(I_\mu\ast|u|^{2^*_\mu})|u|^{2^*_\mu}+
(I_\mu\ast|v|^{2^*_\mu})|v|^{2^*_\mu}]>0.$$ Moreover, using (\ref{4.6}) with $\gamma_p+\gamma_q=2$ we infer $m_\nu(a,b)>0$.
Next we show $m_{\nu}(a,b)<\frac{2^*_\mu-1}{22^*_\mu}
S^{\frac{2^*_\mu}{2^*_\mu-1}}_{H,L}$ and below we only consider the case $q<p$ since the case $q>p$ can be treated similarly. Choose
$u_\epsilon=\frac{a}{|\eta_\epsilon|_2}\eta_\epsilon$,
$v_\epsilon=\frac{\epsilon^\theta}{|\eta_\epsilon|_2}\eta_\epsilon$ as in Lemma \ref{l2.9} with $\theta>0$ to be determined and let $t_\epsilon=t_\nu(u_\epsilon,v_\epsilon)$ be given as in Lemma \ref{l2.4} with $\nu>0$ replacing by $0<\nu<\nu'_0$. Then $t_\epsilon\star(u_\epsilon,v_\epsilon)\in \mathcal{P}_\nu(a,\epsilon^\theta)$. Set $\mathbb{A}$ and $\mathbb{D}$ as in the proof of Lemma \ref{l2.9}.
In view of ${P}_\nu(t_\epsilon\star(u_\epsilon,v_\epsilon))=0$, $\gamma_p+\gamma_q=2$, $0<\nu<\nu'_0$ and (\ref{1.14}) we deduce for $\epsilon>0$ small
\begin{equation*}\aligned
e^{(22^*_\mu-2)t_\epsilon}=\mathbb{A}-2\nu\mathbb{D}
\geq \bigl[1-2\nu C_{N,p,q}(a^2+b^2)^{\frac{p+q-2}2}\bigr]\mathbb{A}\geq C\frac{|\nabla \eta_\epsilon|^2_2|\eta_\epsilon|^{22^*_\mu-2}_2}{\int_{\mathbb{R}^N}
(I_\mu\ast|\eta_\epsilon|^{2^*_\mu})|\eta_\epsilon|^{2^*_\mu}}
.\endaligned\end{equation*}
Then (\ref{2.24}) and (\ref{2.28}) imply $e^{t_\epsilon}\geq C|\eta_\epsilon|_2$ for some $C>0$ as $\epsilon\rightarrow0$.
Therefore, (\ref{2.30.0}) and (\ref{2.34.0}) hold with $\gamma_p+\gamma_q=2$ and $m_{r,\nu}(a,b)$ replacing by $m_{\nu}(a,b)$.
When $N=3$, by the assumptions $p+q=4+\frac{4-2\mu}{3}$ and $2<\mu<3$, we know $\min\{p,q\}<2$ if $\min\{p,q\}>\frac{6-\mu}{2}$. Then we choose $\theta>0$ in (\ref{2.34.0}) such that $\theta q+\frac{p+q-\gamma_p-\gamma_q}{2}<\min\{1,2\theta\}$. When $N=4$, by the assumption $\min\{p,q\}<2$, we choose $\theta=1$ in (\ref{2.34.0}). In any case,
for  $\epsilon>0$ small enough, we always get $m_\nu(a,b)<\frac{2^*_\mu-1}{22^*_\mu}
S^{\frac{2^*_\mu}{2^*_\mu-1}}_{H,L}$. \ \ \ $\Box$

{\bf Proof of Theorems \ref{th1.2} and \ref{th1.4}} Taking into account of Lemma \ref{l3.4} with $c=m_{r,\nu}(a,b)$, Lemma \ref{l2.6} with the condition $\nu>0$ replacing by $0<\nu<\nu'_0$, Lemmas \ref{l2.14} (ii) and \ref{l2.15}, we infer that system (\ref{0.0}) has a positive radial normalized solution $(u,v)$. Furthermore, Lemma \ref{l2.14} (i) implies that $(u,v)$ is a normalized ground state.\ \ \ \ $\Box$
\section{The asymptotic behavior of normalized solutions}
\subsection{Preliminaries}
\renewcommand{\theequation}{4.\arabic{equation}}
In this subsection, we are devoted to studying two limit problems. One type of asymptotic behavior is relevant to the minimizer of $S_{H,L}$ and the following nonlocal critical compactness lemma, which can be viewed as the counterpart of \cite[Theorem 1.41]{WM} about a local critical compactness lemma. Let $u\rightarrow u_{r,x_0}=r^{\frac{N-2}{2}}u(rx+x_0)$ be the rescaling, where $r\in\mathbb{R}^+$ and $x_0\in\mathbb{R}^N$.
\begin{lemma}\label{l5.3.0} Let $\{u_n\}\subset D^{1,2}(\mathbb{R}^N)$ be a minimizing sequence satisfying
\begin{equation}\label{5.3.0}\int_{\mathbb{R}^N}(I_\mu\ast|u_n|^{2^*_\mu})|u_n|^{2^*_\mu}=1,\quad |\nabla u_n|^2_2\rightarrow S_{H,L},\quad\text{as}\ n\rightarrow\infty.\end{equation}
Then there exists a sequence $(r_n,y_n)\subset(0,\infty)\times \mathbb{R}^N$ such that $\{(u_n)_{r_n,y_n}\}$ contains a convergent subsequence in $D^{1,2}(\mathbb{R}^N)$.
\end{lemma}
{\bf Proof}: From (\ref{5.3.0}) we know
\begin{equation}\label{5.3.1}-\Delta u_n-S_{H,L}(I_\mu*|u_n|^{2^*_\mu})|u_n|^{2^*_\mu-2}u_n\rightarrow0.\end{equation}
Assume $u_n\rightharpoonup u$ in $D^{1,2}(\mathbb{R}^N)$. Set $v_n=u_n-u$. Below we distinguish two cases.

{\bf Case 1: $u\neq0$.} Then
\begin{equation}\label{5.3.4}\aligned
S_{H,L}&=|\nabla u_n|^2_2+o_n(1)=|\nabla u|^2_2+|\nabla v_n|^2_2+o_n(1)\\
&\geq S_{H,L}\Bigl(\int_{\mathbb{R}^N}(I_\mu*|u|^{2^*_\mu})
|u|^{2^*_\mu}\Bigr)^{\frac{1}{2^*_\mu}}+
S_{H,L}\Bigl(\int_{\mathbb{R}^N}(I_\mu*|v_n|^{2^*_\mu})
|v_n|^{2^*_\mu}\Bigr)^{\frac{1}{2^*_\mu}}+o_n(1)\\
&\geq S_{H,L}\Bigl[\int_{\mathbb{R}^N}(I_\mu*|u|^{2^*_\mu})
|u|^{2^*_\mu}+\int_{\mathbb{R}^N}(I_\mu*|v_n|^{2^*_\mu})
|v_n|^{2^*_\mu}\Bigr]^{\frac{1}{2^*_\mu}}+o_n(1)\\
&=S_{H,L}\Bigl(\int_{\mathbb{R}^N}(I_\mu*|u_n|^{2^*_\mu})
|u_n|^{2^*_\mu}+o_n(1)\Bigr)^{\frac{1}{2^*_\mu}}+o_n(1)=S_{H,L}+o_n(1).
\endaligned\end{equation}
Since $u\neq0$, the second inequality of (\ref{5.3.4}) implies that $\int_{\mathbb{R}^N}(I_\mu*|v_n|^{2^*_\mu})
|v_n|^{2^*_\mu}\rightarrow0$,  and then the first inequality of (\ref{5.3.4}) implies $|\nabla v_n|^2_2\rightarrow0$. Consequently, $u_n\rightarrow u$ in $D^{1,2}(\mathbb{R}^N)$.

{\bf Case 2: $u=0$.} Then $v_n=u_n$. We may suppose $|u_n|_{2^*}\not\rightarrow0$. Indeed, if $|u_n|_{2^*}\rightarrow0$, then Hardy-Littlewood-Sobolev inequality implies
$\int_{\mathbb{R}^N}(I_\mu*|u_n|^{2^*_\mu})|u_n|^{2^*_\mu}\leq C|u_n|^{22^*_\mu}_{2^*}\rightarrow0$,
contradicts (\ref{5.3.0}). Then there exist $c_0, C_0>0$ independent on $n$ such that $C_0\geq|u_n|_{2^*}\geq c_0$.
Next we claim that there exist $r_n>0$ and $y_n\in\mathbb{R}^N$ such that $(u_n)_{r_n,y_n}\rightarrow w\neq0$ in $D^{1,2}(\mathbb{R}^N)$ as $n\rightarrow\infty$. In fact, define the Levy concentration function
\begin{equation}\label{levy}Q_n(r)=\sup_{z\in\mathbb{R}^N}\int_{B_r(z)}|u_n|^{2^*}dx.\end{equation}
Since $\lim_{r\rightarrow0^+}Q_n(r)=0$ and $\lim_{r\rightarrow\infty}Q_n(r)\geq c^{2^*}_0$, there exist $r_n>0$ and $y_n\in\mathbb{R}^N$ such that
$$Q_n(r_n)=|u_n|^{2^*}_{2^*,B_{r_n}(y_n)}=b, \ \text{with}\ 0<b<\min\Bigl\{S^{\frac{2N}{4-\mu}}{[2S_{H,L}C(N,\mu)C^{2^*_\mu}_0]^{\frac{2N}{\mu-4}}}, \ c^{2^*}_0\Bigr\}.$$
Define $w_n=(u_n)_{r_n,y_n}$, we may assume $w_n\rightharpoonup w$ in $D^{1,2}(\mathbb{R}^N)$ and  $w_n\rightarrow w$ a.e. in $\mathbb{R}^N$. Moreover, $\sup_{z\in\mathbb{R}^N}\int_{B_{1}(z)}|w_n|^{2^*}dx=|w_n|^{2^*}_{2^*,B_{1}(0)}=b$.
It is easy to see that $|u_n|_{2^*}=|w_n|_{2^*}$, (\ref{5.3.1}) and (\ref{5.3.0}) hold with $u_n$ replacing by $w_n$.
If $w=0$, then $w_n\rightarrow0$ in $L^2_{loc}(\mathbb{R}^N)$. Choose $\phi\in C^\infty_0(\mathbb{R}^N)$ such that $supp \phi\subset B_1(y)$ for some $y\in\mathbb{R}^N$. From (\ref{5.3.1}) with $u_n$ replacing by $w_n$, Hardy-Littlewood-Sobolev inequality, Holder inequality and the fact that $|u_n|_{2^*}\leq C_0$ it follows that
$$\aligned\int_{\mathbb{R}^N}|\nabla(\phi w_n)|^2&=\int_{\mathbb{R}^N}\nabla w_n\nabla(\phi^2 w_n)+o_n(1)\leq S_{H,L}\int_{\mathbb{R}^N}(I_\mu\ast|w_n|^{2^*_\mu})\phi^2 |w_n|^{2^*_\mu}+o_n(1)\\&\leq S_{H,L}C(N,\mu)|w_n|^{2^*_\mu}_{2^*}\Bigl(\int_{\mathbb{R}^N} |\phi w_n|^{\frac{4N}{2N-\mu}}|w_n|^{(2^*_\mu-2)\frac{2N}{2N-\mu}}\Bigr)^
{\frac{2N-\mu}{2N}}+o_n(1)\\
&\leq S_{H,L}C(N,\mu)|w_n|^{2^*_\mu}_{2^*}
|w_n|^{2^*_\mu-2}_{2^*,{B_1(y)}}|\phi w_n|^2_{2^*}+o_n(1)\leq\frac12|\nabla(\phi w_n)|^2_2+o_n(1). \endaligned$$
Then $\nabla w_n\rightarrow0$ in $L^2_{loc}(\mathbb{R}^N)$ and $w_n\rightarrow0$ in $L^{2^*}_{loc}(\mathbb{R}^N)$, which contradicts $|w_n|^{2^*}_{2^*,B_{1}(0)}=b$. Hence, $w\neq0$. Taking the arguments in (\ref{5.3.4}) with $u_n$, $u$ and $v_n$ replacing by $w_n$, $w$ and $w_n-w$ respectively,  we infer $w_n\rightarrow w$ in $D^{1,2}(\mathbb{R}^N)$ since $w\neq0$.\ \ \ \ \ $\Box$

As a product, we have the following result.
\begin{lemma}\label{l5.3.1} Let $\{u_n\}\subset D^{1,2}(\mathbb{R}^N)$ be a radial minimizing sequence satisfying (\ref{5.3.0}).
Then there exists a sequence $r_n\in(0,\infty)$ such that $\{(u_n)_{r_n,0}\}$ contains a convergent subsequence in $D^{1,2}(\mathbb{R}^N)$.
\end{lemma}
{\bf Proof}: Note that $u_n$ is radial. Taking the arguments as in Lemma \ref{l5.3.0} with the Levy concentration function (\ref{levy}) replacing by
$$Q_n(r)=\int_{B_r(0)}|u_n|^{2^*}dx,$$ the conclusion follows.\ \ \ $\Box$

Another type of asymptotic behavior is relevant to  system (\ref{limit}) and next we study system (\ref{limit}). For simplicity, denote
\begin{equation}\label{5.5}D^{-1}_0=\inf_{(u,v)\in\mathcal{S}(a,b)}
\frac{(|\nabla u|^2_2+|\nabla v|^2_2)^{\frac{\gamma_p+\gamma_q}{2}}}{\int_{\mathbb{R}^N}(I_\mu*|u|^{p})|v|^{q}}
.\end{equation}
\begin{lemma}\label{l4.2} Let $p+q\neq4+\frac{4-2\mu}{N}$. Then for $\tilde{m}(a,b)$ given in (\ref{1.15}), there holds
$\tilde{m}(a,b)=\frac{\gamma_p+\gamma_q-2}{2(\gamma_p+\gamma_q)}
[D_0(\gamma_p+\gamma_q)]^{\frac{2}{2-\gamma_p-\gamma_q}}$,  $(\gamma_p+\gamma_q-2)\tilde{m}(a,b)>0$ and $D_0>0$.
\end{lemma}
{\bf Proof}: For any $(u,v)\in \mathcal{S}(a,b)$, there exists a unique $t(u,v)$ defined by \begin{equation}\label{4.4}e^{(\gamma_p+\gamma_q-2)t(u,v)}=\frac{|\nabla u|^2_2+|\nabla v|^2_2}{(\gamma_p+\gamma_q)\int_{\mathbb{R}^N}
(I_\mu*|u|^{p})|v|^{q}},\end{equation}
such that $t(u,v)\star(u,v)\in\mathcal{\tilde{P}}(a,b)$, where $\mathcal{\tilde{P}}(a,b)$ is given in (\ref{1.15}), and
$$\tilde{m}(a,b):=\inf_{\mathcal{S}(a,b)}
\tilde{J}(t(u,v)\star(u,v))=\frac{\gamma_p+\gamma_q-2}{2(\gamma_p+\gamma_q)}
[D_0(\gamma_p+\gamma_q)]^{\frac{2}{2-\gamma_p-\gamma_q}}.$$
where $\tilde{J}$ given in (\ref{limitfuctional}).
If $\gamma_p+\gamma_q>2$, as in Lemma \ref{l4.5.0} we get $\tilde{m}(a,b)>0$. If $\gamma_p+\gamma_q<2$, for $(u,v)\in \mathcal{S}(a,b)$, we have $\tilde{J}(t(u,v)\star(u,v))<0$ and then $\tilde{m}(a,b)<0$. Consequently, $D_0>0$. \ \ \  $\Box$
\begin{lemma}\label{l2.17} If $p+q>4+\frac{4-2\mu}{N}$, then $$\tilde{m}(a,b)=\tilde{m}_{r}(a,b):=
\inf_{(u,v)\in\mathcal{\tilde{P}}(a,b)\cap H_{rad}}\tilde{J}(u,v).$$\end{lemma}
{\bf Proof}:  Assume that there is $(u,v)\in \mathcal{\tilde{P}}(a,b)\backslash \mathcal{S}_r(a,b)$ such that
$\tilde{J}(u,v)<\tilde{m}_{r}(a,b).$
In the same way as the  arguments in Lemma \ref{l2.14}, we get $\tilde{P}(u^*,v^*)\leq0$ and exclude the case $\tilde{P}(u^*,v^*)=0$. If $\tilde{P}(u^*,v^*)<0$, by (\ref{4.4}) we have $t_*:=t(u^*,v^*)<0$. Then
$$\aligned \tilde{m}_{r}(a,b)\leq \tilde{J}(t_*\star(u^*,v^*))
&=\frac{\gamma_p+\gamma_q-2}{2(\gamma_p+\gamma_q)}e^{2t_*}(|\nabla u^*|^2_2+|\nabla v^*|^2_2)\\
&\leq\frac{\gamma_p+\gamma_q-2}{2(\gamma_p+\gamma_q)}(|\nabla u|^2_2+|\nabla v|^2_2)
=\tilde{J}(u,v),
\endaligned$$
 contradicts the above hypothesis $\tilde{J}(u,v)<\tilde{m}_{r}(a,b).$\ \ \ \ $\Box$

{\bf Proof of Theorem \ref{th1.5}}
In the case of $p+q>4+\frac{4-2\mu}{N}$, as Lemma \ref{l2.6} one can construct a Palais-Smale sequence
$(u_n,v_n)\subset \mathcal{S}_r(a,b)$ of $\tilde{J}|_{\mathcal{S}(a,b)}$ at the level $\tilde{m}_r(a,b)>0$ with $\tilde{P}(u_n,v_n)\rightarrow0$ and $u^-_n,v^-_n\rightarrow0$ a.e. in $\mathbb{R}^N$ as $n\rightarrow\infty$.
Then\begin{equation}\label{4.1}
\tilde{J}'(u_n,v_n)+\lambda_{1,n}(u_n,0)+\lambda_{2,n}(0,v_n)\rightarrow0,\quad \lambda_{1,n},\lambda_{2,n}\in\mathbb{R}.\end{equation}
As in the proof of Lemma \ref{l3.4.0}, we get $(u_n,v_n)$ is bounded in $H$.
Note that (\ref{3.6}) holds true with $J_\nu$ replacing by $\tilde{J}$. Then $\lambda_{1,n},\lambda_{2,n}$ are bounded. So there exist $(u,v)\in H_{rad}$  and $\lambda_1, \lambda_2$ such that
$(u_n,v_n)\rightharpoonup (u,v)$ in $H_{rad}$, $(u_n,v_n)\rightarrow(u,v)$ in $L^q(\mathbb{R}^N)\times L^q(\mathbb{R}^N)$, $(u_n,v_n)\rightarrow(u,v)$ a.e. in $\mathbb{R}^{2N}$, and $(\lambda_{1,n},\lambda_{2,n})\rightarrow(\lambda_1,\lambda_2)$ in $\mathbb{R}^2$. We claim that $u\neq0$ and $v\neq0$.
Indeed, if $u=0$ or $v=0$, then by $\tilde{P}(u_n,v_n)\rightarrow0$ and (\ref{2.11}) we infer $|\nabla u_n|^2_2+|\nabla v_n|^2_2\rightarrow0$.
Therefore, $\tilde{J}(u_n,v_n)
\rightarrow0,$ contradicts $\tilde{m}_r(a,b)=\tilde{m}(a,b)>0$ using Lemmas \ref{l2.17} and Lemma \ref{l4.2}. Then $0<|u|_2\leq a$ and $0<|v|_2\leq b$. Furthermore, \begin{equation}\label{5.2}
\tilde{J}'(u,v)+\lambda_1(u,0)+\lambda_2(0, v)=0,\
u\geq0,\ v\geq0,
\end{equation}
and so $\tilde{P}(u,v)=0$.

{\bf Subcase 1}. $N\in\{3,4\}$. In this subcase, as in the proof of Lemma \ref{l3.4}, we
get $\lambda_1,\lambda_2>0$. Testing (\ref{5.2}) and (\ref{4.1}) with $(u_n-u,0)$ respectively and subtracting, together with Hardy-Littlewood-Sobolev inequality we obtain
 $$\aligned
 &|\nabla(u_n-u)|^2_2+\lambda_1|u_n-u|^2_2\\
 =&p\int_{\mathbb{R}^N}[(I_\mu\ast|v_n|^q)|u_n|^{p-2}u_n(u_n-u)-
 (I_\mu\ast|v|^q)|u|^{p-2}u(u_n-u)]=o_n(1).\endaligned$$
Hence, $u_n\rightarrow u$ in $H^1(\mathbb{R}^N)$. Analogously, $v_n\rightarrow v$ in $H^1(\mathbb{R}^N)$. Then $(u,v)\in \mathcal{\tilde{P}}(a,b)$ and $\tilde{J}(u,v)=\tilde{m}_r(a,b)=\tilde{m}(a,b)$ using Lemma \ref{l2.17}.

{\bf Subcase 2}. $N\geq5$ with $(\frac N2-1)p+\frac {N}{2}q\leq 2N-\mu$ and $(\frac N2-1)q+\frac {N}{2}p\leq 2N-\mu$. If this subcase occurs, then we firstly claim that $\lambda_1,\lambda_2>0$. Actually, if $\lambda_1\leq0$, testing (\ref{5.2}) with $(u,0)$ we get
 $$|\nabla u|^2_2\geq p\int_{\mathbb{R}^N}(I_\mu\ast|u|^{p})|v|^{q}.$$
 Jointing the above inequality with $\tilde{P}(u,v)=0$ we obtain
\begin{equation*}|\nabla v|^2_2\leq (\gamma_p+\gamma_q-p)\int_{\mathbb{R}^N}(I_\mu\ast|u|^{p})|v|^{q}.\end{equation*}
 This is impossible since $v\neq0$ and the coefficient of $\int_{\mathbb{R}^N}(I_\mu\ast|u|^{p})|v|^{q}$ is non-positive due to $(\frac{N}{2}-1)p+\frac{N}{2}q\leq 2N-\mu$. Therefore, $\lambda_1>0$. Similarly, $\lambda_2>0$. Below as in subcase 1 we deduce that $(u,v)$ is a normalized ground state of (\ref{limit}).

If  $p+q<4+\frac{4-2\mu}{N}$, then we may assume $(u_n,v_n)\in \mathcal{S}(a,b)$ is
a minimizing sequence of $\tilde{m}(a,b)$. Using symmetric
decreasing rearrangement and passing to $(|u_n|,|v_n|)$ we may assume that $u_n$ and $v_n$ are non-negative and radial. By Ekeland's variational principle we may assume $(u_n,v_n)$ is a Palais-Smale sequence. Then taking similar arguments as in the case $p+q>4+\frac{4-2\mu}{N}$, we obtain $\tilde{m}(a,b)$ is achieved by a normalized ground state.\ \ \ \ $\Box$

\subsection{Proof of Theorems \ref{th1.6} and \ref{th1.7}}
By Theorems \ref{th1.1}, \ref{th1.2}, \ref{th1.3} and \ref{th1.4}, for $\nu>0$
small or large, there may exist a positive radial normalized solution or normalized ground state $(u_\nu,v_\nu)\in H$ and $(\lambda_{1,\nu},\lambda_{2,\nu})\in\mathbb{R}^2$ solving system  (\ref{0.0}). In this subsection, our goal
 is to investigate the asymptotic behavior of $(u_\nu,v_\nu)$ as $\nu\rightarrow0^+$ or $\nu\rightarrow+\infty$.

{\bf Proof of Theorem \ref{th1.6}}: (i) For any $(u_0,v_0)\in \mathcal{L}(a,b)$, we have $\tilde{J}(u_0,v_0)=\tilde{m}(a,b)$. Then
\begin{equation}\label{5.6}
\int_{\mathbb{R}^N}(I_\mu*|u_0|^{p})|v_0|^{q}=D_0(|\nabla u_0|^2_2+|\nabla v_0|^2_2)^{\frac{\gamma_p+\gamma_q}{2}},
\end{equation}
where $D_0$ is given in (\ref{5.5}).
In view of  Lemma \ref{l2.2}, there exists $t_\nu:=t_\nu(u_0,v_0)$ such that
$t_\nu\star(u_0,v_0)\in \mathcal{P}^+_\nu(a,b)$. Combining with (\ref{5.5}) we obtain
\begin{equation}\label{5.9}\aligned
&\Bigl[D^{-1}_0\int_{\mathbb{R}^N}(I_\mu*|t_\nu\star{u_0}|^{p})
|t_\nu\star{v_0}|^{q}\Bigr]^{\frac{2}{\gamma_p+\gamma_q}}\leq|\nabla (t_\nu\star{u_0})|^2_2+|\nabla (t_\nu\star{v_0})|^2_2\\ \leq&\nu(22^*_\mu-2)^{-1}(22^*_\mu-\gamma_p-\gamma_q)(\gamma_p+\gamma_q)
\int_{\mathbb{R}^N}(I_\mu*|t_\nu\star{u_0}|^{p})|t_\nu\star{v_0}|^{q}.
\endaligned\end{equation}
Then
\begin{equation}\label{5.15}
\int_{\mathbb{R}^N}(I_\mu*|t_\nu\star{u_0}|^{p})|t_\nu\star{v_0}|^{q}\leq
\mathbb{A},\ \text{with} \ \mathbb{A}=
\Bigl(\frac{22^*_\mu-\gamma_p-\gamma_q}{22^*_\mu-2}
\nu(\gamma_p+\gamma_q)D^{\frac{2}{\gamma_p+\gamma_q}}_0\Bigr)
^{\frac{\gamma_p+\gamma_q}{2-\gamma_p-\gamma_q}}.
\end{equation}
Since $(u_\nu,v_\nu)\in \mathcal{P}^+_\nu(a,b)$, as above we have
\begin{equation}\label{5.16}
\int_{\mathbb{R}^N}(I_\mu*|u_\nu|^{p})|v_\nu|^{q}\leq\mathbb{A}
.\end{equation}
Consider the function
$$\bar{f}(t)=\bigl(\frac12-\frac{1}{22^*_\mu}\bigr)D^{-\frac{2}{\gamma_p+\gamma_q}}_0
t^{\frac{2}{\gamma_p+\gamma_q}}-
\frac{22^*_\mu-\gamma_p-\gamma_q}{22^*_\mu}\nu t.$$
 In view of $t_\nu\star(u_0,v_0)\in \mathcal{P}^+_\nu(a,b)$, the definition of $m_\nu(a,b)$ and (\ref{5.6}) we get
\begin{equation}\label{5.17}\aligned
m_\nu(a,b)\leq J_\nu(t_\nu\star(u_0,v_0))=\bar{f}\Bigl(
\int_{\mathbb{R}^N}(I_\mu*|t_\nu\star{u_0}|^{p})|t_\nu\star{v_0}|^{q}\Bigr).
\endaligned\end{equation}
Moreover, (\ref{5.5}) implies
\begin{equation}\label{5.18}\aligned
m_\nu(a,b)=J_\nu(u_\nu,v_\nu)\geq\bar{f}\Bigl(
\int_{\mathbb{R}^N}(I_\mu*|u_\nu|^{p})|v_\nu|^{q}\Bigr).
\endaligned\end{equation}
Note that $\bar{f}(t)$ is decreasing in $(0,\mathbb{A})$, from (\ref{5.15})-(\ref{5.18}) we know
\begin{equation}\label{5.19}
\int_{\mathbb{R}^N}(I_\mu*|t_\nu\star{u_0}|^{p})|t_\nu\star{v_0}|^{q}
\leq\int_{\mathbb{R}^N}(I_\mu*|u_\nu|^{p})|v_\nu|^{q}.\end{equation}
Combining the second inequality of (\ref{5.9}) and $P_\nu(t_\nu\star(u_0,v_0))=0$ we get
\begin{equation}\label{5.2.2}
\nu\int_{\mathbb{R}^N}(I_\mu*|t_\nu\star{u_0}|^{p})|t_\nu\star{v_0}|^{q}
\sim|\nabla (t_\nu\star{u_0})|^2_2+|\nabla (t_\nu\star{v_0})|^2_2.
\end{equation}
Note that $\gamma_p+\gamma_q<2$. Then
\begin{equation}\label{5.2.1}e^{t_\nu}\thicksim
\nu^{\frac{1}{2-\gamma_p-\gamma_q}},\quad\text{as}\ \nu\rightarrow0^+.\end{equation}
Using  $P_\nu(t_\nu\star(u_0,v_0))=0$ again there holds
\begin{equation}\label{5.20}
e^{t_\nu}=(1+o_\nu(1))\Bigl(\frac{\nu(\gamma_p+\gamma_q)
\int_{\mathbb{R}^N}(I_\mu*|u_0|^{p})|v_0|^{q}
}{|\nabla u_0|^2_2+|\nabla v_0|^2_2}\Bigr)^{\frac{1}{2-\gamma_p-\gamma_q}}=(1+o_\nu(1))\nu^{\frac{1}{2-\gamma_p-\gamma_q}}.
\end{equation}

Since $(u_\nu,v_\nu)$ solves system (\ref{0.0}), $(\tilde{u}_\nu,\tilde{v}_\nu):=(-t_\nu)\star(u_\nu,v_\nu)$ solves
 \begin{equation*}
\aligned
\left\{ \begin{array}{lll}
-\Delta \tilde{u}_\nu+\lambda_{1,\nu}e^{-2t_\nu}\tilde{u}_\nu=&e^{(22^*_\mu-2)t_\nu}
(I_\mu\ast |\tilde{u}_\nu|^{2^*_\mu})|\tilde{u}_\nu|^{2^*_\mu-2}\tilde{u}_\nu\\&+\nu pe^{(\gamma_p+\gamma_q-2)t_\nu}(I_\mu\ast |\tilde{v}_\nu|^q)|\tilde{u}_\nu|^{p-2}\tilde{u}_\nu\  \text{in}\quad \mathbb{R}^N,\\
-\Delta \tilde{v}_\nu+\lambda_{2,\nu}e^{-2t_\nu}\tilde{v}_\nu=&e^{(22^*_\mu-2)t_\nu}(I_\mu\ast |\tilde{v}_\nu|^{2^*_\mu})|\tilde{v}_\nu|^{2^*_\mu-2}\tilde{v}_\nu\\&+\nu qe^{(\gamma_p+\gamma_q-2)t_\nu}(I_\mu\ast |\tilde{u}_\nu|^p)|\tilde{v}_\nu|^{q-2}\tilde{v}_\nu\  \text{in}\quad \mathbb{R}^N,\\
|\tilde{u}_\nu|^2_2=a^2,\quad|\tilde{v}_\nu|^2_2=b^2.
\end{array}\right.\endaligned
\end{equation*}
Below we claim that
\begin{equation}\label{claim}\lambda_{1,\nu}+\lambda_{2,\nu}\thicksim |\nabla u_\nu|^2_2+|\nabla v_\nu|^2_2\thicksim \nu^{\frac{2}{2-\gamma_p-\gamma_q}},\quad\text{as}\quad \nu\rightarrow0^+.\end{equation}
In fact, since $(u_\nu,v_\nu)\in \mathcal{P}^+_\nu(a,b)$, (\ref{5.2.2}) holds with $t_\nu\star(u_0,v_0)$ replacing by $(u_\nu,v_\nu)$. Then combining with (\ref{1.14}) we obtain
\begin{equation}\label{5.2.0}
\nu\int_{\mathbb{R}^N}(I_\mu*|u_\nu|^{p})|v_\nu|^{q}\sim|\nabla u_\nu|^2_2+|\nabla v_\nu|^2_2\lesssim\nu^{\frac{2}{2-r_p-r_q}}.
\end{equation}
For any $(u,v)\in \mathcal{S}(a,b)$ and $\nu>0$ small, by Lemma \ref{l2.2} there exists a unique $t'_\nu:=t_\nu(u,v)$ such that $t'_\nu\star(u,v)\in \mathcal{P}^+_\nu(a,b)$.
As (\ref{5.2.1}) we have
$e^{t'_\nu}\thicksim \nu^{\frac{1}{2-\gamma_p-\gamma_q}}$ as $\nu\rightarrow0^+.$ Thus
 $$\aligned
 J_\nu(t'_\nu\star(u,v))=&(\frac12-\frac1{\gamma_p+\gamma_q})e^{2t'_\nu}(|\nabla u|^2_2+|\nabla v|^2_2)+(\frac1{\gamma_p+\gamma_q}\\
&-\frac{1}{22^*_\mu})e^{22^*_\mu t'_\nu}\int_{\mathbb{R}^N}
[(I_\mu*|u|^{2^*_\mu})|u|^{2^*_\mu}+(I_\mu*|v|^{2^*_\mu})|v|^{2^*_\mu}]\sim-\nu^{\frac{2}{2-\gamma_p-\gamma_q}}.
\endaligned$$
Using $J_\nu(t'_\nu\star(u,v))\geq m_\nu(a,b)$ and $m_\nu(a,b)\gtrsim-\nu\int_{\mathbb{R}^N}(I_\mu*|u_\nu|^{p})|v_\nu|^{q}$. Then $\nu\int_{\mathbb{R}^N}(I_\mu*|u_\nu|^{p})
|v_\nu|^{q}\gtrsim\nu^{\frac{2}{2-\gamma_p-\gamma_q}},$ together with (\ref{5.2.0}) we deduce \begin{equation*}\nu\int_{\mathbb{R}^N}(I_\mu*|u_\nu|^{p})|v_\nu|^{q}\sim|\nabla u_\nu|^2_2+|\nabla v_\nu|^2_2\sim\nu^{\frac{2}{2-r_p-r_q}}.
\end{equation*}
From Pohozaev identity (\ref{1.8.0}) and testing (\ref{0.0}) with $(u_\nu,v_\nu)$ it follows that $$\lambda_{1,\nu}|u_\nu|^2_2+\lambda_{2,\nu}|v_\nu|^2_2=\nu(p+q-\gamma_p-\gamma_q)
\int_{\mathbb{R}^N}(I_\mu*|u_\nu|^{p})|v_\nu|^{q}.$$
Then (\ref{claim}) follows.

By (\ref{claim}) and (\ref{5.20}) we know
$|\nabla \tilde{u}_\nu|^2_2+|\nabla \tilde{v}_\nu|^2_2\thicksim1.$
Then $\{(\tilde{u}_\nu,\tilde{v}_\nu)\}$ is bounded in $H$. Since $({u}_\nu,{v}_\nu)$ is radial, we get $(\tilde{u}_\nu,\tilde{v}_\nu)$ is also radial. Passing to a subsequence, we assume $(\tilde{u}_\nu,\tilde{v}_\nu)\rightharpoonup (\tilde{u},\tilde{v})$ in $H_{rad}$, $(\tilde{u}_\nu,\tilde{v}_\nu)\rightarrow (\tilde{u},\tilde{v})$ in $L^q(\mathbb{R}^N)\times L^q(\mathbb{R}^N)$ for $2<q<2^*$ and $(\tilde{u}_\nu,\tilde{v}_\nu)\rightarrow (\tilde{u},\tilde{v})$ a.e. in $\mathbb{R}^{2N}$. Using (\ref{claim}) and (\ref{5.20}) again we have $\{\lambda_{1,\nu}e^{-2t_\nu}\}$ and $\{\lambda_{2,\nu}e^{-2t_\nu}\}$ are bounded in $\nu$. Hence we may assume $\lambda_{1,\nu}e^{-2t_\nu}\rightarrow\lambda_1$ and $\lambda_{2,\nu}e^{-2t_\nu}\rightarrow\lambda_2$ as $\nu\rightarrow0^+$.
On the other hand, by (\ref{5.20}) we get
$e^{(22^*_\mu-2)t_\nu}\rightarrow0$ and $\nu e^{(\gamma_p+\gamma_q-2)t_\nu}\rightarrow1$ as $\nu\rightarrow0^+$. Then
\begin{equation*}
\aligned
\left\{ \begin{array}{lll}
-\Delta \tilde{u}+\lambda_1\tilde{u}= p(I_\mu\ast |\tilde{v}|^q)|\tilde{u}|^{p-2}\tilde{u}\ & \text{in}\quad \mathbb{R}^N,\\
-\Delta \tilde{v}+\lambda_2\tilde{v}= q(I_\mu\ast |\tilde{u}|^p)|\tilde{v}|^{q-2}\tilde{v}\ & \text{in}\quad \mathbb{R}^N.
\end{array}\right.\endaligned
\end{equation*}
As (\ref{2.11}) we get
$$\int_{\mathbb{R}^N}(I_\mu*|\tilde{u}_\nu|^{p})|\tilde{v}_\nu|^{q}\rightarrow
\int_{\mathbb{R}^N}(I_\mu*|\tilde{u}|^{p})|\tilde{v}|^{q}.$$
Together with  (\ref{5.19}) and $(u_\nu,v_\nu)=t_\nu\star( \tilde{u}_\nu, \tilde{v}_\nu)$, we have
\begin{equation}\label{5.14}\int_{\mathbb{R}^N}(I_\mu*|\tilde{u}|^{p})|\tilde{v}|^{q}\geq\int_{\mathbb{R}^N}(I_\mu*|u_0|^{p})
|v_0|^{q},\end{equation}
which implies $\tilde{u}\neq0$ and $\tilde{v}\neq0$. As in the proof of Theorem \ref{th1.5} we deduce $\lambda_1,\lambda_2>0$, $\tilde{u}_\nu\rightarrow\tilde{u}$ and $\tilde{v}_\nu\rightarrow\tilde{v}$ in $H^1(\mathbb{R}^N)$ as $\nu\rightarrow0^+$. Therefore $\tilde{m}(a,b)\leq \tilde{J}(\tilde{u},\tilde{v})$. By (\ref{5.14}) we know
\begin{equation*}\aligned  \tilde{J}(\tilde{u},\tilde{v})=\frac{\gamma_p+\gamma_q-2}{2}
\int_{\mathbb{R}^N}(I_\mu*|\tilde{u}|^{p})|\tilde{v}|^{q}\leq\frac{\gamma_p+\gamma_q-2}{2}
\int_{\mathbb{R}^N}(I_\mu*|u_0|^{p})|v_0|^{q}=\tilde{J}(u_0,v_0).\endaligned\end{equation*}
Then $\tilde{J}(\tilde{u},\tilde{v})=\tilde{m}(a,b)$ since $\tilde{J}(u_0,v_0)=\tilde{m}(a,b)$. Thus, $(\tilde{u},\tilde{v})\in \mathcal{L}(a,b)\cap H_{rad}$.

(ii) For any $(u,v)\in \mathcal{S}_r(a,b)$, there exists a unique $t_\nu(u,v)>0$ such that $t_\nu(u,v)\star(u,v)\in \mathcal{P}^-_\nu(a,b)\cap H_{rad}$. Using similar arguments as that of (\ref{claim}), we can show that $e^{t_\nu(u,v)}\lesssim\nu^{\frac{1}{2-\gamma_p-\gamma_q}}$ and $m_{r,\nu}(a,b)\lesssim \nu^{\frac{2}{2-\gamma_p-\gamma_q}}$. Then
 $|\nabla u_\nu|^2_2+|\nabla v_\nu|^2_2\rightarrow0$
as $\nu\rightarrow+\infty$.
By Pohozaev identity (\ref{1.8.0}), (\ref{5.5}) and the definition of $S_{H,L}$ we get
\begin{equation*}\aligned
(\gamma_p+\gamma_q)\nu\int_{\mathbb{R}^N}
(I_\mu*|u_\nu|^{p})
|v_\nu|^{q}
\geq&(|\nabla u_\nu|^2_2+|\nabla v_\nu|^2_2)\bigl[1-S^{-2^*_\mu}_{H,L}(|\nabla u_\nu|^2_2+|\nabla v_\nu|^2_2)^{2^*_\mu-1}\bigr]\\
\geq& D^{-\frac{2}{\gamma_p+\gamma_q}}_0(1+o_\nu(1))\bigl(\int_{\mathbb{R}^N}
(I_\mu*|u_\nu|^{p})
|v_\nu|^{q}\bigr)^{\frac{2}{\gamma_p+\gamma_q}}.
\endaligned\end{equation*}
Then
\begin{equation}\label{5.23}\int_{\mathbb{R}^N}(I_\mu*|u_\nu|^{p})
|v_\nu|^{q}\geq(1+o_\nu(1))[(\gamma_p+\gamma_q)\nu D^{\frac{2}{\gamma_p+\gamma_q}}_0]^{\frac{\gamma_p+\gamma_q}
{2-\gamma_p-\gamma_q}}.\end{equation}
For any $(u_0,v_0)\in \mathcal{L}(a,b)\cap H_{rad}$, let $t_\nu:=t_\nu(u_0,v_0)$ be given in Lemma \ref{l2.4}. From $P_\nu(t_\nu\star(u_0,v_0))=0$ and (\ref{5.6}) it follows that
\begin{equation}\label{5.23.0}\aligned e^{(\gamma_p+\gamma_q)t_\nu}
\int_{\mathbb{R}^N}(I_\mu*|u_0|^{p})|v_0|^{q}&\leq [(\gamma_p+\gamma_q)\nu D^{\frac{2}{\gamma_p+\gamma_q}}_0]^{\frac{\gamma_p+\gamma_q}
{2-\gamma_p-\gamma_q}}.\endaligned\end{equation}
Therefore by (\ref{5.23.0}) we deduce
$$\aligned m_{r,\nu}(a,b)\leq& J_\nu(t_\nu\star(u_0,v_0))=\bigl(\frac{\gamma_p+\gamma_q}2-1\bigr)\nu
e^{(\gamma_p+\gamma_q)t_\nu}\int_{\mathbb{R}^N}(I_\mu*|u_0|^{p})|v_0|^{q}\\&
+\frac{2^*_\mu-1}{22^*_\mu}e^{22^*_\mu t_\nu}
\int_{\mathbb{R}^N}[(I_\mu*|{u_0}|^{2^*_\mu})|{u_0}|^{2^*_\mu}+(I_\mu*
|{v_0}|^{2^*_\mu})|{v_0}|^{2^*_\mu}]\\
\leq&(1+o_\nu(1))\nu\bigl(\frac{\gamma_p+\gamma_q}2-1\bigr)[(\gamma_p+\gamma_q)\nu D^{\frac{2}{\gamma_p+\gamma_q}}_0]^{\frac{\gamma_p+\gamma_q}
{2-\gamma_p-\gamma_q}}.
\endaligned$$
Observe that
\begin{equation*}\label{5.23.1}
m_{r,\nu}(a,b)=J_\nu(u_\nu,v_\nu)\geq\frac{\gamma_p+\gamma_q-2}{2}\nu\int_{\mathbb{R}^N}(I_\mu*|u_\nu|^{p})
|v_\nu|^{q}.\end{equation*}
Then $$\int_{\mathbb{R}^N}(I_\mu*|u_\nu|^{p})
|v_\nu|^{q}\leq(1+o_\nu(1))[(\gamma_p+\gamma_q)\nu D^{\frac{2}{\gamma_p+\gamma_q}}_0]^{\frac{\gamma_p+\gamma_q}
{2-\gamma_p-\gamma_q}},$$
together with (\ref{5.23}), (\ref{5.6}) and $\tilde{P}(u_0,v_0)=0$ it follows that
\begin{equation}\label{5.28}
\aligned\int_{\mathbb{R}^N}(I_\mu*|{u_\nu}|^{p})|{v_\nu}|^{q}
=&(1+o_\nu(1))[(\gamma_p+\gamma_q)\nu D^{\frac{2}{\gamma_p+\gamma_q}}_0]^{\frac{\gamma_p+\gamma_q}
{2-\gamma_p-\gamma_q}}\\=&(1+o_\nu(1))\nu^{\frac{\gamma_p+\gamma_q}{2-\gamma_p-\gamma_q}}
\int_{\mathbb{R}^N}(I_\mu*|{u_0}|^{p})|{v_0}|^{q}.\endaligned\end{equation}
Similarly
\begin{equation}\label{5.24.0}|\nabla u_\nu|^2_2+|\nabla v_\nu|^2_2=(1+o_\nu(1))[(\gamma_p+\gamma_q)\nu D_0]^{\frac{2}
{2-\gamma_p-\gamma_q}}.\end{equation}

Let $(\tilde{u}_\nu,\tilde{v}_\nu):=s_\nu\star({u}_\nu,{v}_\nu)$
with $e^{s_\nu}=\nu^{\frac{1}{\gamma_p+\gamma_q-2}}$. Then
$|\tilde{u}_\nu|^2_2=a^2$, $|\tilde{v}_\nu|^2_2=b^2$ and from (\ref{5.24.0}) we get $\{(\tilde{u}_\nu,\tilde{v}_\nu)\}$ is bounded in $H$ and assume $(\tilde{u}_\nu,\tilde{v}_\nu)\rightharpoonup (\tilde{u},\tilde{v})$ in $H_{rad}$. Using (\ref{5.28}) and arguing as (\ref{2.11}) there holds
\begin{equation*}\int_{\mathbb{R}^N}(I_\mu\ast |\tilde{u}|^p)|\tilde{v}|^{q}=\lim_{\nu\rightarrow+\infty}\int_{\mathbb{R}^N}(I_\mu\ast |\tilde{u}_\nu|^p)|\tilde{v}_\nu|^{q}=\int_{\mathbb{R}^N}(I_\mu\ast |{u}_0|^p)|{v}_0|^{q}.\end{equation*}
Below taking the argument as in the proof of (i), we have $(\tilde{u}_\nu,\tilde{v}_\nu)\rightarrow(\tilde{u},\tilde{v})$ in $H$ and $(\tilde{u},\tilde{v})\in \mathcal{L}(a,b)\cap H_{rad}$. \ \ \ \ $\Box$

{\bf Proof of Theorem \ref{th1.7}}:
(i) If $p+q= 4+\frac{4-2\mu}{N}$, then $m_{\nu}(a,b)$ is achieved by a positive radial normalized ground state $(u_\nu,v_\nu)$.
Since $m_{\nu}(a,b)$ is bounded with respect to $\nu$ small, as in Lemma \ref{l3.4.0} or Remark \ref{r1} one easily has that $\{(u_\nu,v_\nu)\}$ is bounded in $H$. Let $|\nabla u_\nu|^2_2\rightarrow l_1$ and $|\nabla v_\nu|^2_2\rightarrow l_2$. If $l_1=l_2=0$, then $\lim_{\nu\rightarrow0^+}m_{\nu}(a,b)=\lim_{\nu\rightarrow0^+}J_\nu(u_\nu,v_\nu)=0$. However, in view of Lemmas \ref{l2.15} and \ref{l2.14} we get
$\lim_{\nu\rightarrow0^+}m_{\nu}(a,b)>0$. This is a contradiction. So $l_1+l_2>0$. By Pohozaev identity (\ref{1.8.0}) and $\nu\rightarrow0^+$, as (\ref{3.5}) we infer
$l_1+l_2\geq S^{\frac{2^*_\mu}{2^*_\mu-1}}_{H,L}$ and equality holds if and only if $l_1=0$ or $l_2=0$.
Moreover
$$\aligned
\frac{2^*_\mu-1}{22^*_\mu}S^{\frac{2^*_\mu}{2^*_\mu-1}}_{H,L}
\geq\lim_{\nu\rightarrow0^+}m_{\nu}(a,b)=
\lim_{\nu\rightarrow0^+}J_\nu(u_\nu,v_\nu)=\frac{2^*_\mu-1}{22^*_\mu}
(l_1+l_2)\geq
\frac{2^*_\mu-1}{22^*_\mu}S^{\frac{2^*_\mu}{2^*_\mu-1}}_{H,L},
\endaligned$$
and therefore either $l_1=0$ or $l_2=0$. We may assume
$l_1=0$ since the case $l_2=0$ can be treated similarly. Then $l_2=S^{\frac{2^*_\mu}{2^*_\mu-1}}_{H,L}$. By Pohozaev identity (\ref{1.8.0}) and $\nu\rightarrow0^+$ again we get
$$|\nabla u_\nu|^2_2\rightarrow0,\ |\nabla v_\nu|^2_2\rightarrow S^{\frac{2^*_\mu}{2^*_\mu-1}}_{H,L},\ \int_{\mathbb{R}^N}(I_\mu\ast|v_\nu|^{2^*_\mu})|v_\nu|^{2^*_\mu}\rightarrow S^{\frac{2^*_\mu}{2^*_\mu-1}}_{H,L}.$$
Let $$w_\nu:=
\bigl(\int_{\mathbb{R}^N}
(I_\mu\ast|v_\nu|^{2^*_\mu})
|v_\nu|^{2^*_\mu}\bigr)^{-\frac{1}{22^*_\mu}}v_\nu.$$
Since $v_\nu$ is radial, applying Lemma \ref{l5.3.1}, there exists $r_\nu>0$ such that $(w_\nu)_{r_\nu,0}\rightarrow w$ in $D^{1,2}(\mathbb{R}^N)$. Moreover, $w$ is a radial minimizer of $S_{H,L}$ and satisfies
$$-\Delta w=S_{H,L}(I_\mu\ast|w|^{2^*_\mu})|w|^{2^*_\mu-2}w.$$
Letting $v=S^{\frac1{2(2^*_\mu-1)}}_{H,L}w$, we have $(v_\nu)_{r_\nu,0}\rightarrow v$ in $D^{1,2}(\mathbb{R}^N)$ and
 Lemma \ref{l1.1} implies that $v=\widetilde{U}_{\epsilon,0}$ for some
 $\epsilon>0$.

 (ii) If $p+q>4+\frac{4-2\mu}{N}$, then taking some arguments as the proof of (i) with $m_{\nu}(a,b)$ replacing by $m_{r,\nu}(a,b)$, the conclusion (ii) follows.\ \ \ \ $\Box$

\section{Nonexistence results in the case $\nu<0$}
\renewcommand{\theequation}{5.\arabic{equation}}
In this section, we show there is no normalized ground state of (\ref{0.0}) when $\nu<0$.
\begin{lemma}\label{l3.2}Let $\nu\leq0$. Then the following results hold.\\
\noindent(i) $m_\nu(a,b)=\inf_{\Lambda_\nu(a,b)\cup\mathcal{P}_\nu(a,b)}\bar{J}_\nu
=\inf_{\Lambda_\nu(a,b)}\bar{J}_\nu$, where
$$\Lambda_\nu(a,b)=\{(u,v)\in \mathcal{S}{(a,b)}: P_\nu(u,v)<0\},\  \bar{J}_\nu(u,v)
:={J}_\nu(u,v)-\frac{1}{22^*_\mu}P_\nu(u,v),$$
and $m_\nu(a,b)$ satisfies (\ref{2.1}).\\
\noindent(ii)  $m_\nu(a,b)=m_0(a,b)$.
\end{lemma}
{\bf Proof}: (i) It is easy to see that Lemma \ref{l2.4} holds with $0<\nu<\nu_0$ replacing by $\nu\leq0$. Then $\Lambda_\nu(a,b)$ is dense in $\Lambda_\nu(a,b)\cup\mathcal{P}_\nu(a,b)$, and so
$\inf_{\Lambda_\nu(a,b)}\bar{J}_\nu=
\inf_{\Lambda_\nu(a,b)\cup\mathcal{P}_\nu(a,b)}\bar{J}_\nu\leq m_\nu(a,b).$
Moreover, for $(u,v)\in \Lambda_\nu(a,b)$, Lemma \ref{l2.4} with $0<\nu<\nu_0$ replacing by $\nu\leq0$ implies $t_\nu(u,v)<0$. From the definition of $\bar{J}_\nu$ it follows that
$\bar{J}_\nu(u,v)>\bar{J}_\nu(t_\nu(u,v)\star(u,v))
\geq\inf_{\mathcal{P}_\nu(a,b)}\bar{J}_\nu.$
Then $\inf_{\Lambda_\nu(a,b)}\bar{J}_\nu= \inf_{\Lambda_\nu(a,b)\cup\mathcal{P}_\nu(a,b)}\bar{J}_\nu=m_\nu(a,b)$. Taking the argument as in Lemma \ref{l2.7}, we have $m_\nu(a,b)$ satisfies (\ref{2.1}).

(ii) The fact that $\inf_{\Lambda_0(a,b)}\bar{J}_0\leq
\inf_{\Lambda_\nu(a,b)}\bar{J}_\nu$ and the conclusion (i) imply
$m_0(a,b)\leq m_\nu(a,b)$.
On the other hand, for any $(u,v)\in\mathcal{P}_0(a,b)$, set $u_t=t^{\frac{N-2}{2}}u(tx)$ and $v_t=t^{\frac{N-2}{2}}v(tx)$ with $t>0$. Then for $t>1$ large enough we infer
$$\aligned P_\nu(u_t,v_t)
=&P_0(u,v)+\frac{(\gamma_p+\gamma_q)|\nu|}{t^{p+q-\gamma_p-\gamma_q}}
\int_{\mathbb{R}^N}(I_\mu\ast|u|^p)|v|^q<0.\endaligned$$
Using the conclusion (i) again, for $t\rightarrow+\infty$ there holds
$$\aligned
m_\nu(a,b)&\leq m_\nu(t^{-1}a,t^{-1}b)\leq \bar{J}_\nu(u_t,v_t)\\
&=\bar{J}_0(u,v)+\nu\bigl(\frac{\gamma_p+\gamma_q}{22^*_\mu}-1\bigr)t^{\gamma_p+\gamma_q-p-q}
\int_{\mathbb{R}^N}(I_\mu\ast|u|^p)|v|^q\rightarrow \bar{J}_0(u,v).
\endaligned$$
Then $m_\nu(a,b)\leq m_0(a,b)$ and so $m_\nu(a,b)=m_0(a,b)$.\ \ \ \ $\Box$

{\bf Proof of Theorem \ref{th1.8}} For $K$ given in (\ref{3.6.2}), by (\ref{3.6.1}) one easily has
\begin{equation*}
\inf_{u\in \mathcal{S}_{a}}\max_{t\in\mathbb{R}}K(t\star u)
=\inf_{u\in H^1(\mathbb{R}^N)\backslash\{0\}}\frac{2^*_\mu-1}{22^*_\mu}\frac{|\nabla u|^{\frac{22^*_\mu}{2^*_\mu-1}}_2}{\bigl(\int_{\mathbb{R}^N}(I_\mu\ast|u|^{2^*_\mu})|u|^{2^*_\mu}\bigr)^{\frac1{2^*_\mu-1}}}
=\frac{2^*_\mu-1}{22^*_\mu} S^{\frac{2^*_\mu}{2^*_\mu-1}}_{H,L},\end{equation*}
where $\mathcal{S}_a:=\{u\in H^1(\mathbb{R}^N):|u|_2=a\}$. Then for any $\epsilon>0$, there exists $u\in \mathcal{S}_a$ such that
$$\max_{t\in\mathbb{R}}K(t\star u)<\frac{2^*_\mu-1}{22^*_\mu}
S^{\frac{2^*_\mu}{2^*_\mu-1}}_{H,L}+\epsilon.$$
Set $v\in \mathcal{S}_b$ and let $t(s):=t_0(u,s\star v)$ be such that $P_0(t(s)\star(u,s\star v))=0$.
Then $\limsup_{s\rightarrow-\infty}t(s)<+\infty$ as in Lemma \ref{l2.7}. Thus
$$\aligned J_0(t(s)\star(u,s\star v))=K(t(s)\star u)+K((t(s)+s)\star v)\leq\frac{2^*_\mu-1}{22^*_\mu}
S^{\frac{2^*_\mu}{2^*_\mu-1}}_{H,L}+\epsilon+o(1),\endaligned$$
as $s\rightarrow-\infty$. Then $m_0(a,b)\leq \frac{2^*_\mu-1}{22^*_\mu}
S^{\frac{2^*_\mu}{2^*_\mu-1}}_{H,L}$. On the other hand,
for any $(u,v)\in \mathcal{P}_0(a,b)$, we have
\begin{equation*}\aligned
|\nabla u|^2_2+|\nabla v|^2_2
&\leq S^{-2^*_\mu}_{H,L}(|\nabla u|^{22^*_\mu}_2+|\nabla v|^{22^*_\mu}_2)< S^{-2^*_\mu}_{H,L}(|\nabla u|^{2}_2+|\nabla v|^{2}_2)^{2^*_\mu}.
\endaligned\end{equation*}
Then \begin{equation}\label{3.5.0}
|\nabla u|^2_2+|\nabla v|^2_2> S^{\frac{2^*_\mu}{2^*_\mu-1}}_{H,L},\end{equation}
and so $m_0(a,b)\geq\frac{2^*_\mu-1}{22^*_\mu}
S^{\frac{2^*_\mu}{2^*_\mu-1}}_{H,L}$.
Combining with Lemma \ref{l3.2} (ii) it follows that $m_\nu(a,b)=m_0(a,b)=\frac{2^*_\mu-1}{22^*_\mu}
S^{\frac{2^*_\mu}{2^*_\mu-1}}_{H,L}$.
Observe that for any $(u,v)\in \mathcal{P}_\nu(a,b)$ with $\nu\leq0$, (\ref{3.5.0})  still holds and
$$J_\nu(u,v)\geq
\frac{2^*_\mu-1}{22^*_\mu}(|\nabla u|^2_2+|\nabla v|^2_2).$$
Consequently, $m_\nu(a,b)$ with $\nu\leq0$ is not achieved.\ \ \ \ $\Box$

\end{document}